\documentclass[reqno,11pt, a4paper]{amsart}
\usepackage[utf8]{inputenc}
\usepackage[T1,T5]{fontenc}
\usepackage{amsmath,amsfonts,amssymb,amsthm}
\usepackage{appendix}
\usepackage[left=2cm,right=2cm,top=2cm,bottom=2.5cm]{geometry}
\usepackage{stmaryrd} 
\usepackage{mathrsfs}
\usepackage{subfiles} 
\usepackage{yhmath} 
\usepackage{cases} 
\usepackage{ulem} 
\usepackage{mathtools,color,url} 
\usepackage{abraces}
\mathtoolsset{showonlyrefs}  
\numberwithin{equation}{section}  

\newtheorem{theo}{Theorem}
\newtheorem{lemm}[theo]{Lemma}
\newtheorem{coro}[theo]{Corollary}
\newtheorem{defi}[theo]{Definition}
\newtheorem{propo}[theo]{Proposition}

\newtheorem{claim}[theo]{Claim}
\newtheorem{rema}[theo]{Remark}

\newcommand\R{{\ensuremath {\mathbb R} }}

\newcommand\N{{\ensuremath {\mathbb N} }}

\newcommand{\cE}{\mathcal{E}}

\newcommand{\eps}{\epsilon}

\title[Strong interactions for fmKdV]{Strongly interacting solitary waves for the fractional modified Korteweg-de Vries equation}
\author[A.Eychenne]{Arnaud Eychenne}
\address{Département of Mathematics, University of Bergen, Allégaten 41
Realfagbygget
5007 Bergen, Norway}
\email{arnaud.eychenne.waxweiler@gmail.com}
\author[F.Valet]{Frédéric Valet}
\address{Département of Mathematics, University of Bergen, Allégaten 41
Realfagbygget
5007 Bergen, Norway}
\email{frederic.valet@uib.no}
\date{\today}

\begin{document}

\maketitle

\begin{abstract}
We study one particular asymptotic behaviour of a solution of the fractional modified Korteweg-de Vries equation (also known as the dispersion generalised modified Benjamin-Ono equation):
\begin{align}\tag{fmKdV}
    \partial_t u + \partial_x (-\vert D \vert^\alpha u + u^3)=0.
\end{align}
The dipole solution is a solution behaving in large time as a sum of two strongly interacting solitary waves with different signs. We prove the existence of a dipole for fmKdV. A novelty of this article is the construction of accurate profiles. Moreover, to deal with the non-local operator $\vert D \vert^\alpha$, we refine some weighted commutator estimates.
\end{abstract}


\section{Introduction}

\subsection{Introduction of the equation}

This article is dedicated to the fractional-modified Korteweg-de Vries equation (also known as the dispersion generalised modified Benjamin-Ono equation):
\begin{align}\tag{fmKdV}\label{mBOt}
    \partial_t u + \partial_x \left( -\vert D \vert^\alpha u + u^3 \right)=0, \quad u: I_t \times \mathbb{R}_x \rightarrow \mathbb{R},  \quad 1 < \alpha<2,
\end{align}
where $I_t$ is a time interval, $\partial_x$ (respectively $\partial_t$) denotes the space (respectively time) derivative, and the symbol $\vert D \vert^\alpha$ is defined by the Fourier transform as an operator acting on the space of distributions:
\begin{align*}
    \mathcal{F}(\vert D \vert^\alpha u) (\xi) := \vert \xi \vert^\alpha \mathcal{F}(u)(\xi).
\end{align*}

For the purposes of motivating the equation, let us introduce the more generalised equation:
\begin{align}\label{eq:general_dispersive}
    \partial_t u + \mathcal{L} \partial_x u +\partial_x(f(u))=0.
\end{align}
The operator $L$ represents the dispersion of the equation, and $f(u)$  stands for the non-linearity. 

In the case of a quadratic non-linearity $f(u)=u^2$ and a dispersion $\mathcal{L}=-\vert D \vert^\alpha$, we get respectively the Benjamin-Ono equation (BO) and the Korteweg-de Vries equation (KdV) for $\alpha=1$ and $\alpha=2$. Shrira and Voronovich, in \cite{SV96}, introduced the equation of coastal waves, where the parameter is the evolution of the depth of the coast. If the evolution of the depth is algebraic and given by $-(1+X)^{\alpha-1}$, for $\alpha\in (1,2)$, then the dispersion operator is approximated, for waves with a small wave number, by $\displaystyle -c \vert D \vert ^{\alpha}$. Notice that other dispersions are justified by Klein, Linares, Pilod and Saut \cite{KLPS18}.

While the change of dispersion in the quadratic case models different phenomena, the change of non-linearity helps to understand the balance between non-linearity and dispersion. Indeed, studying equations with a cubic non-linearity $f(u)=u^3$ and different dispersions give new insights of the competition between those two terms. The case $\mathcal{L}=\partial_x^2=-\vert D \vert^2$ corresponds to the modified Korteweg-de Vries equation (mKdV), while the case $\mathcal{L}=-\vert D \vert$ corresponds to the modified Benjamin-Ono equation (mBO). We chose in this article to focus on the case of a non-local dispersion $\mathcal{L}=-\vert D \vert^\alpha$, with $1 < \alpha<2$.

Since, for $1<\alpha<2$, \ref{mBOt} does not enjoy a Lax pair as KdV, BO or mKdV, no tools from complete integrability can be applied to this equation. On the other hand, \ref{mBOt} possesses 3 conserved quantities (at least formally):
\begin{align*}
    \int_{\mathbb{R}} u(t,x) dx, \quad \frac{1}{2} \int_{\mathbb{R}} u^2(t,x) dx, \quad \int_{\mathbb{R}} \left( \frac{( \vert D \vert^\alpha u (t,x) )^2}{2} - \frac{u^4(t,x)}{4}\right)dx.
\end{align*}

We define the scaling operators by:
\begin{equation} \label{op:scaling}
    \forall \lambda \in \mathbb{R}_+^*, \quad u \mapsto u_\lambda, \quad \text{with} \quad u_{\lambda}(t,x):=\lambda^{\frac{\alpha}{2(1+\alpha)}} u(\lambda t,\lambda^{\frac{1}{1+\alpha}} x ).
\end{equation}
The set of solutions of \ref{mBOt} is fixed under the scaling operations. The mBO equation is mass-critical in the sense that the $L^2$-norm is preserved under any scaling operation. Meanwhile, \ref{mBOt} is mass-subcritical since the conserved space under the operator of scaling is the homogeneous Sobolev space $\dot{H}^s(\R)$ with $s=\frac{1-\alpha}{2}<0$ as soon as $\alpha>1$. The equation \ref{mBOt} has been proved to be locally well-posed in $H^{s}(\R)$ for $s \geq \frac{3-\alpha}{4}$ by Guo \cite{Guo12}, and the flow is locally continuous on that space. As a consequence, the equation is globally well-posed in the energy space $H^{\frac{\alpha}{2}}(\R)$ (see Appendix \ref{appendix:LWP}). We also refer to Guo and Huang \cite{GH22}, Kim and Schippa \cite{KS21}, Molinet and Tanaka \cite{MT22} for other well-posedness results. Moreover, in the case $\alpha=1$ the problem is locally well-posed in the energy space, see Kenig and Takaoka \cite{KT06}.

\subsection{Ground states and solitary waves.}

Different coherent structures may appear in the study of non-linear dispersive equations, and solitary waves are one of them. A solitary wave is a solution $u(t,x)=Q_c(x-ct)$ moving at a velocity $c$ in one direction, decaying at infinity and keeping its form along the time. The function $Q_c$ satisfies the elliptic equation:
\begin{align}
    -\vert D \vert^\alpha Q_c - c Q_c + Q_c^3=0.\label{eq:elliptic_GS}
\end{align}

A remarkable point is the existence of those objects for any velocity $c>0$. Unlike the mKdV equation the solutions $Q_c$ of \eqref{eq:elliptic_GS} are not explicit. The existence of a such  solution of the elliptic problem \eqref{eq:elliptic_GS} is related to the existence of a minimizer of an adequate functional.
Such a minimizer is called a ground state, and the existence of a ground state has been proved by Weinstein in  \cite{weinstein1987existence} and Albert-Bona-Saut in \cite{albert1997model}. Moreover, the ground state is positive. For now, the notation $Q_c$ will refer to the ground-state of the functional.

If we denote by $Q$ the positive ground state associated to $c=1$, all the other ground states $Q_c$ associated to the different values $c>0$ can be expressed in terms of the ground state $Q$ by the operation of scaling \eqref{op:scaling}:
\begin{align*}
    Q_c(x) = (Q_1)_c(x).
\end{align*}

The question of the uniqueness of the ground state of \eqref{eq:elliptic_GS} is difficult and has been solved by Frank-Lenzmann in \cite{FL13}. Note however that no result seems to be known for the uniqueness of solutions to \eqref{eq:elliptic_GS} which do not minimize the Euler-Lagrange functional. The non-locality of the operator $|D|^{\alpha}$ does not allow to use classical ODE's tools for this equation. The uniqueness of the solution of the non-local elliptic problem  \eqref{eq:elliptic_GS} is derived from the non-degerenency of the linearized operator 
\begin{align*}
    L=|D|^{\alpha}+1-3Q^2,
\end{align*}
by proving that $\ker(L)=\text{span}(Q')$.
This result was obtained by Frank-Lenzmann in \cite{FL13}. The proof is based on an extension process to the upper half-plane, introduced by Caffarelli-Silvestre \cite{CS07}, which allows to look at the operator $|D|^{\alpha}$ as a Dirichlet-Neumann operator.

Furthermore, as soon as $\alpha<2$, the function $Q$ has a algebraic decay (see \eqref{asympt:Q} for a more precise expansion):
\begin{align*}
    Q(x) \simeq_{+\infty} \frac{1}{x^{1+\alpha}}.
\end{align*}

The question of stability of a solitary wave in this case has been done by Angulo Pava \cite{Ang18}, see also \cite{LNP22}.

\vspace{0.5cm}

One conjecture in the field of dispersive equation states that any solution decomposes, at large time, into different dispersive objects (such as the solitary waves) plus a radiation term. Whereas the solitary waves move to the right, the radiation term moves to the left. This conjecture has been proved for the KdV equation using the tools of complete integrability, but remains open in the non-integrable cases. It is then natural to introduce multi-solitary waves, which are solutions $u$ that in large time $[T_0,+\infty)$ are close to a sum of $K$ decoupled solitary waves:
\begin{defi}
Let $K>0$, and $K$ different velocities $0<c_1<\cdots<c_K$. A function $u$ is called a multi-solitary waves associated to the previous velocities (or pure multi-solitary waves) if there exist $T_0>0$, $K$ functions $v_k : (T_0,+\infty)\rightarrow \mathbb{R}$ such that:
\begin{align*}
    \lim_{t \rightarrow +\infty} \left\| u(t) - \sum_{k=1}^K Q_{c_k}(\cdot - v_k(t)) \right\|_{H^{\frac{\alpha}{2}}} =0 \quad \text{and} \quad  \forall k \in (1,K), \quad \left\vert v_k(t) -c_kt \right\vert =o_{+\infty}(t).
\end{align*}
\end{defi}

Notice that the definition of a multi-solitary waves may depend on the information one can get from those objects. For example, in a recent result by the first author \cite{eychenne2021asymptotic}, the proof of the existence of the multi-solitary waves has been established for the equation fKdV with a dispersion $\alpha\in (\frac{1}{2},2)$ and an explicit rate of convergence of the solution to the sum of the $K$-decoupled solitary waves. Notice that the proof can easily be adapted to \ref{mBOt}, establishing then the existence of multi-solitary waves for this equation for $1<\alpha<2$. The proof of existence of those objects is a first step toward the soliton resolution conjecture for this equation. 

\subsection{Dipoles and main theorem}

Notice that in the previous definition of multi-solitary waves, all the velocities are distinct. One can wonder if there exist solutions $u$ behaving at infinity as a sum of two solitary waves with the same velocity $c$ and different signs. A solution satisfying this definition is called a dipole. 
In particular, if the two solitary waves have the same velocity, they interact in large time one with each other, and the velocity of the different solitary waves is thus expected to be of the form $ v_k(t) \sim_{+\infty} ct- g_k(t)$, with $g_k(t) = o_{+\infty}(t)$.

This object has first been observed on the mKdV equation using the complete integrability of the equation \cite{WO82}. For an odd non-linearity $f(u) = \vert u \vert^{p-1}u$, $p\in (2,5)$ and a dispersion $\mathcal{L}=\partial_x^2$, Nguyen in \cite{Ngu17} proved the existence of dipoles for those equations that are not completely integrable. 

In this paper, we prove the existence of a dipole for the \ref{mBOt} in the $L^2$-subcritical case:

\begin{theo}\label{maint_theo:version1}
Let $\alpha\in (1,2)$. There exist some constant $T_0>0,C>0$ and $U\in C^{0}([T_0,+\infty):H^{\frac{\alpha}{2}}(\R))$ solution of \eqref{mBOt} such that, for all $t\geq T_0$:
\begin{align*}
    \left\lVert U(t,\cdot) + Q\left(\cdot -t- \frac{a}{2}t^{\frac{2}{\alpha+3}}\right) - Q\left(\cdot -t+ \frac{a}{2}t^{\frac{2}{\alpha+3}}\right) \right\rVert_{H^{\frac{\alpha}{2}}}\leq C t^{-\frac{\alpha-1}{4(\alpha+3)}},
\end{align*}
where 
\begin{align}\label{defi:a_b_1}
a:=\left(\frac{\alpha+3}{2}\sqrt{\frac{-4b_1}{\alpha+1}} \right)^{\frac{2}{\alpha+3}}\quad \text{and} \quad \displaystyle b_1:=-2\frac{(\alpha+1)^2}{\alpha-1}\frac{\sin(\frac{\pi}{2}\alpha)}{\pi}\int_{0}^{+\infty}e^{-\frac{1}{r^{\alpha}}}dr\frac{\|Q\|^6_{L^{3}}}{\|Q\|^2_{L^2}}<0.
\end{align}
\end{theo}

This result sheds new light on the relation between the dispersion $\mathcal{L}$ and the distance between two solitary waves of a dipole. Indeed, Nguyen in \cite{Ngu17,Ngu19} studied the case of a dispersion $\mathcal{L}= -\vert D \vert^2=\partial_x^2$ and different non-linearities, which corresponds to the generalized Korteweg-de Vries equation. Since the ground states $Q$ have an exponential decay $e^{-\vert x \vert}$, the distance between the two solitary waves of a dipole is logarithmic in time $2\ln(tc)$, with $c$ depending on the non-linearity. A second example is the recent preprint of Lan and Wang \cite{LW22}, where they studied the generalized Benjamin-Ono equation with a dispersion $\mathcal{L}=-\vert D \vert=-\mathcal{H}\partial_x$ with $\mathcal{H}$ the Hilbert transform and different non-linearities. For this equation, since the ground states have a prescribed algebraic decay $ x^{-2}$, the solitary waves of the dipoles they studied have a distance $\alpha \sqrt{t} + \beta \ln(t) + \gamma$, where $\alpha$, $\beta$ and $\gamma$ are constants dependent only on the non-linearity. Theorem \ref{maint_theo:version1} emphasises how the dispersion influences the distance between the two solitary waves, that is $at^{\frac{2}{\alpha+3}}$. One can conjecture that the dipoles for an equation $\mathcal{L}=-\vert D \vert^\alpha$, for $\alpha\in (1,2)$ and a non-linearity $f(u)=\vert u \vert^{p-1}u$ with various values of $p$, are composed of two solitary waves at a distance $ct^{\frac{2}{\alpha+3}}$, with a constant $c$ dependent on $p$.

\subsection{Related results}

As explained in the introduction, the behaviour of a solution of \eqref{mBO} is determined by the balance between the non-linearity and the dispersion, therefore blow-ups are expected in the critical and super-critical cases. An important result for blow-up, in finite or infinite time, in a non-local setting has been obtained by Kenig-Martel-Robbiano in \cite{KMR09} for:
\begin{align*}
    \partial_tu -\partial_x|D|^{\alpha}u+|u|^{2\alpha}u=0.
\end{align*}
This equation is critical for all the values of $\alpha$. For $\alpha=2$ in the former equation, which corresponds to the critical general Korteweg-de Vries equation, Merle \cite{Mer01} proved the existence of blow-up solutions in finite or infinite time. Using this result, \cite{KMR09} proved by a perturbative argument the existence of blow-up for all $\alpha\in(\alpha_1,2]$, for some $1<\alpha_1<2$. The proof is based on the existence of a Liouville property and localized energy estimates. Those localized estimates generalize the pioneering work of Kenig and Martel \cite{KM09} for the asymptotic stability of the soliton of the Benjamin-Ono equation.

In the case $\alpha=1$ in fmKdV, the equation is $L^2$-critical and blow-up phenomena occur. Bona-Kalisch \cite{BK04}, and Klein-Saut-Wang \cite{KSW22} studied numerically the critical fmKdV and conjecture a blow-up in finite time for this equation. In \cite{MP17} Martel-Pilod proved rigorously the existence of minmial mass blow-up solution for mBO. We mention also the result by Kalisch-Moldabayev-Verdier in \cite{KMV17}, where they observed that two solitary waves may interact in such a way that the smaller wave is annihilated. 

For the super-critical case we refer to the work of Saut-Wang in \cite{SW21}, where they proved the global well-posedness for small initial data and \cite{KSW22} for numerical simulation of blow-up in finite time.

The phenomenon of strong interaction between two different objects also occurs in different situations. Let us enumerate different families of equations and results (this list may not be exhaustive) by beginning with the KdV family. By using the integrable structure of mKdV, Wadati and Ohkuma \cite{WO82} exhibited the existence of a dipole. More recently, Koch and Tataru \cite{KT20} characterized the set of complex two-solitons as a 8-dimensional symplectic submanifold of $H^s$ for $s>-\frac{1}{2}$. The explicit formula of a dipole holds for the mKdV equation only. In the non-integrable case Nguyen \cite{Ngu17} proved the existence of a dipole for \eqref{eq:general_dispersive} for a dispersion $\mathcal{L}=\partial_x^2$ and a non-linearity $f(u)=\vert u \vert^{p-1}u$, with $p \in (2,5)$. Moreover, he discovered that for each super-critical equation with a non-linearity $p>5$, there exists a dipole formed by two solitary waves with same signs, and the distance between the two objects is also logarithmic in time. Inspired by this result, Lan and Wang \cite{LW22} looked for the phenomenon of dipoles for a dispersion $\mathcal{L}=-\vert D \vert$ and a non-linearity $f(u)=\vert u \vert^{p-1}u$, with various values of $p\neq 3$. We also list some results in the setting of the strong interaction of two non-linear objects in the non-linear Schrödinger setting. Ovchinnikov and Sigal \cite{OS98} for the time-dependent Ginzburg Landau equation, with two vortices with different signs; Krieger, Martel and Raphaël \cite{KMR09} for the three dimensional gravitational Hartree equation with two solitons; Nguyen \cite{Ngu19} for the subcritical non-linear Schrödinger with two solitary waves  with different signs, and the same signs for the super-critical case; Nguyen and Martel \cite{MN20} for coupled non-linear Schrödinger, for two solitary waves with different velocities. The phenomenon of dipole also appears in the family of wave equations: Gerard, Lenzmann, Pocovnicu and Raphaël \cite{GLPR18} for the cubic half-wave equation; Côte, Martel, Yuan and Zhao \cite{CMYZ21} for the damped Klein-Gordon equation; Aryan \cite{Ary22} for the Klein-Gordon equation; Jendrej and Lawrie \cite{JL22} for the wave maps equation.

The strong interaction between different objects also gives rise to exotic behaviours. For example, the existence of strongly interacting objects has been proved with multi-solitary waves for the mass-critical non-linear Schrödinger equation by Martel and Raphaël \cite{MR18} and with bubbles for the critical gKdV equation by Combet and Martel \cite{CM18}.

Even if the question of dipoles occur at infinity, one can wonder what happens on the real line to a solution that behaves like a two soliton at $-\infty$. The problem of inelastic collision of two solitary waves has been investigated by Mizumachi \cite{Miz03}, Martel and Merle \cite{MM11,MM11inelastic} and Mu\~noz \cite{Mun10} for non-integrable equations in the KdV family. Indeed, only the completely integrable equations exhibit an elastic collision, that is a solution that can be decomposed at $+\infty$ with the same decomposition as at $-\infty$ (up to phase shift).

We end this part with open questions related to the dipoles of \ref{mBOt}. We begin with the particular case of the critical equation mBO: we do not know if the dipole phenomenon exists for this equation. For a fixed dispersion $\mathcal{L}=-\vert D \vert^\alpha$, one can also wonder about the importance of the non-linearity $f(u)=\vert u \vert^{p-1} u$ : if $p$ is close to $1$, does the structure of a dipole still make sense, or does the non-linearity breaks the structure? Concerning the \ref{mBOt} equation, if a solution behaves at time $-\infty$ as a sum of two different solitary waves, what will be the behaviour of this solution at $+\infty$? Even though this article does not answer those questions, it gives insights and tools to tackle those problems with non-local dispersion.

\subsection{Ideas of the proof}

Let us perform the following change of variables. Let $y:=x-t$, then $v(t,y):=u(t,x)$ verifies
\begin{align}\label{mBO}
    \partial_tv+\partial_y\left( -v-|D|^{\alpha}v+v^3 \right)=0.
\end{align}
This equation is better suited than \ref{mBOt} for the phenomenon of strong interaction, since most of the objects considered here are moving at a velocity close to $1$. Theorem \ref{maint_theo:version1} can be rewritten in this new setting:

\begin{theo}\label{main_theo}
Let $\alpha\in (1,2)$. There exist some constant $T_0>0,C>0$ and $w\in C^{0}([T_0,+\infty):H^{\frac{\alpha}{2}}(\R))$ solution of \eqref{mBO} such that, for all $t\geq T_0$:
\begin{align}
    \left\lVert w(t,\cdot) + Q(\cdot - \frac{a}{2}t^{\frac{2}{\alpha+3}}) - Q(\cdot + \frac{a}{2}t^{\frac{2}{\alpha+3}}) \right\rVert_{H^{\frac{\alpha}{2}}}\leq C t^{-\frac{\alpha-1}{4(\alpha+3)}},
\end{align}
with the constant $a$ defined in \eqref{defi:a_b_1}.
\end{theo}
From now on, we focus on proving the existence of the function $w$. We provide some ideas for the proof of Theorem \ref{main_theo}.

The first important point is the construction of a good approximation. We look for a solution closed to the sum of two solitary waves $-R_1 +R_2$ modulated by a set of parameters $\Gamma=(z_1,z_2,\mu_1,\mu_2)$, where $z_i(t)$ correspond to the centres of the solitary waves moving along the time, whereas $1+\mu_i(t)$ correspond to their size. To this aim, we search for an accurate description of $w+R_1-R_2$, and we introduce the approximation $V$ of the form $V(t,x)= -R_1(t,x)+R_2(t,x)+b(t)W(t,x)-P_1(t,x)+P_2(t,x)$. The goal is to adapt the four other functions such that $V$ almost solves \eqref{mBOt}, in the sense that the quantity $\mathcal{E}_V$ is close to $0$, with:
\begin{align*}
    \mathcal{E}_V:= \partial_t V+ \partial_y (- \vert D \vert^\alpha V -V +V^3).
\end{align*}
By computing the time derivative of $R_1$ and $R_2$, four intrinsic directions appear: $\partial_y R_1$, $\partial_y R_2$, $\Lambda R_1$ and $\Lambda R_2$. For convenience, we will write them under a vector form by $\overrightarrow{MV}$. They go hand in hand with the derivatives of the modulation parameters $\dot{z}_1$, $\dot{z}_2$, $\dot{\mu}_1$ and $\dot{\mu}_2$. Then, the function $W$ is inherent to the problem : it compensates two of those specific directions, and has a plateau between $z_2$ and $z_1$. Even if the previous constructions of strong interactions (\cite{MM11, Ngu17,Ngu19,MN20}) used this function, it seems to be the first time that it  is understood as an intrinsic part of the evolution of the solitary waves, and not only as a part of the profiles $P_i$. With this function we understand how the dispersion of the first solitary wave $-R_1$ on the front influences the second solitary wave in the back, and vice-versa. Once this function $W$ is defined, we fix the functions $P_1$ and $P_2$ with algebraic decay to cancel the remainder terms with algebraic decay too, concentrated around the solitary waves. As a conclusion of this construction, the error can be decomposed into:
\begin{align*}
    \mathcal{E}_V = \overrightarrow{m} \cdot \overrightarrow{MV} + \partial_y S + T,
\end{align*}
with $\overrightarrow{MV}$ containing the four peculiar directions cited above, $\overrightarrow{m}$ gives a system of ODEs that is satisfied by $\Gamma$ and adapted from the interaction terms. The two other source terms, $S$ and $T$ are error terms coming from the rough approximation and are bounded by functions depending on $\Gamma$. If one wants to go further in the development of the approximation, it suffices to extract from $S$ and/or $T$ the terms at the next order to build more precise profiles.

Once the approximation $V$ is constructed, the second step is to estimate the error between the approximation and a solution, and to find a set of equations satisfied by $\mu:=\mu_1-\mu_2$ and $z:=z_1-z_2$. Fix $S_n>>0$, and $v_n$ the solution of \ref{mBOt} with final condition $v_n(S_n)=V(S_n)$. We estimate the $H^{\frac{\alpha}{2}}$-norm of the error backward in time by using an adequate weighted functional, mostly composed of quadratic terms in the error. Whereas studying the error by the energy is quite classic, we adapt in this article the energy functional used by Nguyen \cite{Ngu17} by adding a source term $\int S \epsilon$, linear in $\epsilon$. This trick has been used by Martel and Nguyen \cite{MN20}, by mixing the source term $S$ in the functional, and allows to get rid of the term $\int \partial_y L S \epsilon$ in the functional. It generally helps to get a better approximation of the functional, but in our case, the use of the modified energy enables us not to compute the high Sobolev norms of the source term $S$. It means in particular that the influence of $S$ on the error of the approximation is lower than the one of $T$. 

One technical issue of this functional, as opposed to the ones previously used in this context, is the appearance of the non-local operator $\vert D \vert^\alpha$: two of the difficulties are the singularity of this operator for low frequencies, and the lack of an explicit Leibniz rule for this operator and the weight $\phi$. To bypass those difficulties, we generalize the weighted commutator estimates given in Lemma 6 and Lemma 7 of Kenig-Martel-Robbiano \cite{KMR09} and of the first author \cite{eychenne2021asymptotic}.

These estimates rely on the understanding of the operator $|D|^{\alpha}$. Since the operator is singular at frequency $0$, we need to localize in high and low frequencies : for the high frequencies, we use the pseudo-differential calculus, and the low frequency part is dealt with the theory of bounded operators on $L^2$. In particular, this method implies important restrictions on the choice of the weight.

When orthogonality conditions are imposed to the error, we get a system of ODEs ruling the behaviour of $z$ and $\mu$ in $\overrightarrow{m}$. Roughly speaking, the system is the following:
\begin{align*}
    \dot{\mu}(t) \sim \frac{2b_1}{z^{\alpha+2}(t)}, \quad \text{and} \quad \mu(t) \sim \dot{z}(t).
\end{align*}
Notice that it is the solution of this system that gives the distance between the two solitary waves in Theorem \ref{main_theo}.

To obtain a suitable bound on the different unknowns, we use a bootstrap argument. The more important ones are the error, the parameters $z$ and $\mu$. The error is dealt with the previous functional and $\mu$ by the bootstrap argument. Notice that a bootstrap argument alone would not have been sufficient to close the estimates: because of the algebraic decay in time of the different parameters, several integrations in time can not close the estimates. A topological argument, as introduced by C\^ote, Martel and Martel in \cite{CMM11}, is necessary to conclude the estimate on $z$: roughly speaking, this argument of connectedness asserts that there exists at least one initial data $z^{in}$, chosen in a fixed interval of initial data, such that the estimates hold on the all time interval. Once this initial data is chosen, the all set of estimates is proved to hold on $[T_0,S_n]$.

With these estimates in hand, a classical argument of extraction by compactness allows to get an adequate initial data. By weak-continuity of the flow, we prove that the chosen initial data is close at any time to the sum of the two decoupled solitary waves. Furthermore, we obtain the algebraic decay in time of the error between the final solution and the two solitary waves.

\subsection{Outline of the paper}

The paper is organised as follow. Section \ref{sec:ground_state} is dedicated to the properties related to the ground-state $Q$. It contains in particular the more recent results on those objects, the properties on the linearized operator and various lemmas related to this operator. Section \ref{sec:approximation} contains the construction of an approximation of the solution. Notice that the proof of the main theorem of this part can be skipped at first lecture. In section \ref{sec:modulation}, we give the modulation theorem to describe a solution close to the multi-solitary waves with strong interaction. Section \ref{sec:proof_main_theorem} provides the proof of the existence of the solution. The appendices recall satellite results used in this article : well-posedness, the pseudo-differential calculus, proofs of various lemmas based on pseudo-differential calculus, and the coercivity of the localised linearized operator.

\subsection{Notations}

Throughout the article, we use the following notations.

We denote by $C$ a positive constant, changing from lines to lines independent of the different parameters.

We say $x\sim y$ if there exists $0<c_1<c_2<+\infty$ such that $c_1 x\leq y\leq c_2 y$.

The japanese bracket $\langle \cdot \rangle$ is defined on $\mathbb{R}$ by $\left\langle x \right\rangle := (1+ x^2)^{\frac{1}{2}}$.

$L^2(\R)$ is the set of square integrable functions. We denote the scalar product on $L^{2}(\R)$ by $\langle u,v\rangle:=\int_{\R}u(x)v(x)dx$ with $u,v\in L^2(\R)$. The Fourier transform is defined by: 
\begin{align*}
    \forall f \in L^2(\R), \quad \hat{f}(\xi) := \int_{\R} e^{i x\xi} f(x) dx.
\end{align*}
We define the following spaces:
\begin{itemize}
    \item the Sobolev space, for $s\in \mathbb{R}$ : $\displaystyle H^s(\mathbb{R}):=\left\{f\in L^{2}(\mathbb{R}): \displaystyle\int_{\R}(1+|\xi|^2)^{\frac{s}{2}} \hat{f}(\xi)d\xi<+\infty\right\}$,
    \item the Schwartz space : $\displaystyle \mathcal{S}(\mathbb{R})= \left\{ f \in \mathcal{C}^\infty(\mathbb{R}); \forall \alpha \in \mathbb{N}, \forall \beta \in \mathbb{N}, \exists C_{\alpha,\beta}, \vert f^{\alpha}(x) \vert \leq C_{\alpha, \beta} \left\langle x \right\rangle^{-\beta}\right\}$,
    \item the set of functions with enough decay:
\begin{align}\label{defi:X}
    X^{s}(\mathbb{R}):=\left\{f\in H^{s}(\R): \exists C>0, \forall x\in\R, \quad  |f(x)|\leq \frac{C}{\langle x\rangle^{1+\alpha}}\right\}, \quad \text{and} \quad X^\infty(\mathbb{R}) = \bigcap_{s\in \mathbb{N}} X^s(\mathbb{R}).
\end{align} 
\end{itemize}

Let $f,g\in L^{2}(\R)$. We say that $f$ is orthogonal to $g$ if $\displaystyle \int_{\R}f(x)g(x)dx=0$, and is sometimes shortened by $f\perp g$.

$Q$ is the ground-state associated to the elliptic problem \eqref{eq:Q}, and for $c>0$, we set $Q_c(x) := c^{\frac{\alpha}{2(\alpha+1)}} Q(c^{\frac{1}{1+\alpha}}x)$. Moreover, let us define:
\begin{align}\label{defi:lambda_lambda_Q}
    &\Lambda Q_c := \frac{d}{dc'}Q_{c'_{|c'=c}}= \frac{1}{c}\bigg[\frac{\alpha}{2(\alpha+1)}Q + \frac{1}{\alpha+1}xQ'\bigg]_{c}, \\ &\Lambda^2 Q_c:= \frac{d^2}{dc'^2} Q_{\vert c'=c} = \frac{1}{c^2} \left(- \frac{\alpha(\alpha+2)}{4(\alpha+1)^2}Q +\frac{x^2 Q''}{(\alpha+1)^2} \right)_c.
\end{align}
    
The parameters of the approximation are $z_1$, $z_2$, $\mu_1$ and $\mu_2$. We denote by $\Gamma=(z_1,z_2,\mu_1,\mu_2)$ the set of those parameters. $z$ and $\mu$ are defined in \eqref{defi:mu_z}, and $\bar{z}$ and $\bar{\mu}$ in \eqref{defi:mu_bar}.
The two solitary waves are defined by:
\begin{align}
    R_1(\Gamma,y):=Q_{1+\mu_1}(y-z_1), \quad R_2(\Gamma,y):=Q_{1+\mu_2}(y-z_2), \quad \text{and} \quad \Lambda R_i(t,y) := (\Lambda Q_{1+\mu_i(t)}) (y-z_i(t)).
\end{align}
Along the article, the functions $z_1$, $z_2$, $\mu_1$, $\mu_2$ and $\Gamma$ can depend on the time, and it is precised when needed. The asset of this notation is to remark that the two solitary waves depend on the time through the parameter $\Gamma$. For purposes of notations, we can denote the solitary waves by $R_i(t)$ to emphasize on the time dependency. The solitary waves dependent only on the translation parameters are denoted by:
\begin{align} \label{def:R_tilde}
    \tilde{R}_i(t,y) :=Q(y-z_i(t)) \quad \text{and} \quad \Lambda \tilde{R}_i(t,y) := (\Lambda Q) (y-z_i(t)).
\end{align}

The derivatives are denoted by $\partial_y$ and $\partial_t$. The notation $\nabla_\Gamma$ holds for the gradient along the four directions of $\Gamma$. When no confusion is possible, we denote by prime (as in $Q'$) the space derivative, and by a dot (as in $\dot{\mu}$) the time derivative.


\section{Ground state}\label{sec:ground_state}

This part recalls the properties known on $Q$: existence, uniqueness and the recently proved asymptotic expansion. We emphasize that the asymptotic expansion is composed of terms with algebraic decay, and is thus different from the one of the (gKdV) family $-cQ_c+\Delta Q_c + Q_c^p =0$, with exponential decay. Next, we focus our attention on the linearized operator $L$.

\subsection{Ground state properties}

Considering the equation \eqref{mBO}, the existence of solitary waves is related to the existence of solutions to the following elliptic time-independent equation:
\begin{align}\label{eq:Q}
    - \vert D \vert^\alpha Q - Q +Q^3=0, \quad 1 < \alpha < 2.
\end{align}

The previous elliptic equation is related to a calculus of variation problem. If $Q$ is a minimizer of the following functional $J^\alpha$:
\begin{align}\label{GagliardoNirenberg}
		\displaystyle J^{\alpha}(v)=\frac{\left(\displaystyle\int ||D|^{\frac{\alpha}{2}}v|^2\right)^{\frac{1}{\alpha}}\left(\displaystyle\int |v|^2 \right)^{2-\frac{1}{\alpha}}}{\displaystyle\int \vert v\vert ^4}.
\end{align}
then it is a solution to the elliptic problem.

We now sum up the previous known results on the ground states, which are the minimizer of $J^\alpha$.

\begin{theo}[\cite{albert1997model,FL13,FLS16,weinstein1987existence}]\label{theo:resum_Q}
	Let $\alpha\in(1,2)$. There exists $Q\in H^{s}(\R)$ for all $s\geq 0$ such that  
	\begin{enumerate}
	\item (\textit{Existence}) The function $Q$ solves \eqref{eq:Q} and  $Q=Q(|x|)>0$ is even, positive and strictly decreasing in $|x|$. Moreover, the function $Q$ is a minimizer of $J^{\alpha}$ in the sense that: 
		\begin{equation} 
		J^{\alpha}(Q)=\inf_{v \in H^{\frac{\alpha}2}(\mathbb R)}J^{\alpha}(v).
		\end{equation}	
	\item (\textit{Uniqueness}) The even ground state solution $Q=Q(|x|)>0$ of \eqref{GagliardoNirenberg} is unique, up to the multiplication by a constant,  scaling and translation.
	\item (\textit{Decay}) The function $Q$ verifies the following decay estimate:
		\begin{align}
			\frac{C_1}{(1+|x|)^{1+\alpha}} \leq Q(x)\leq \frac{C_2}{(1+|x|)^{1+\alpha}},
		\end{align} 
     for some $C_1,C_2>0$.
     \item(Gagliardo-Niremberg inequality) There exists a constant $C=C(\alpha)$ such that:
\begin{align}\label{eq:GN1}
    \left\| v \right\|_{L^4} \leq C \left\| v \right\|_{L^2}^{1- \frac{1}{2\alpha}} \left\| v \right\|_{H^\frac{\alpha}{2}}^{\frac{1}{2\alpha}}.
\end{align}
	\end{enumerate}
\end{theo}

\begin{rema}
Notice that since the non-linearity is cubic, the function $Q$ in the theorem and $-Q$ are both solutions of the elliptic equation \eqref{eq:Q}.
\end{rema}

\begin{proof}
    We give some classic ideas to prove the Gagliardo-Niremberg inequality. The proof of this inequality relies on finding a universal constant $C$ which bounds $(J^\alpha)^{-1}$. Indeed, by denoting the following $2$-parameters transformation for $(\lambda, \gamma)\in \mathbb{R}_+^*\times \mathbb{R}_+^*$:
    \begin{align*}
        v_{\lambda,\gamma}(x):= \lambda v\left( \frac{x}{\gamma}\right),
    \end{align*}
    we notice that $J^\alpha(v_{\lambda,\gamma}) = J^\alpha(v)$. As a consequence, if for any $v$, the inequality is proved for some $v_{\lambda,\gamma}$ for some values of $\lambda$ and $\gamma$, then the inequality is proved for any function. In particular, if we choose $\lambda = \| v \|_{L^2}^{-1}$ and $\gamma = \| v \|_{\dot{H}^\frac{\alpha}{2}}^{\frac{2}{\alpha-1}}$, we have $\| v_{\lambda, \gamma} \|_{L^2}= \| v_{\lambda, \gamma} \|_{\dot{H}^\frac{\alpha}{2}}= 1$. Thus it suffices to prove that $(J^\alpha)^{-1}(v_{\lambda,\gamma})= \| v_{\lambda,\gamma} \|_{L^4}^4$ is uniformly bounded with the constraints on $\lambda$ and $\gamma$. By the Sobolev embedding for a certain constant $C$ (see for example \cite{DD07}):
    \begin{align*}
        \| v \|_{L^4} \leq C \| v \|_{H^\frac{\alpha}{2}},
    \end{align*}
    and thus we have $(J^\alpha)^{-1}(v_{\lambda,\gamma})\leq 2C$ independently of $v_{\lambda, \gamma}$. This last inequality concludes the proof of the Gagliardo-Niremberg inequality.
\end{proof}

\begin{rema}
As from \cite{weinstein1987existence,albert1997model,kenig2011local}, the optimal constant in the Gagliardo-Niremberg inequality can be given explicit in terms of $Q$.
\end{rema}

Recently, the asymptotic expansions of the ground states have been improved, see \cite{EV22}. We recall the results applied to our case:

\begin{theo}[\cite{EV22}]\label{thm:asympQ}
Let $\alpha\in(1,2)$ and $x>1$. The positive, even function $Q$ defined in Theorem \ref{theo:resum_Q} verifies:
\begin{enumerate}
    \item(\textit{First-order expansion}) The function $Q$ verifies the following decay estimate:
		\begin{align}
	\left\lvert Q^{(j)}(x) -(-1)^j \frac{(\alpha+j)!}{\alpha!}\frac{a_1}{x^{1+{\alpha}+j}}\right\rvert \leq \frac{C_j}{(1+x)^{2+\alpha+j}}, \quad j\in\N,\label{eq:decay_Q_j}
		\end{align} 
		for some $C_j>0$, with $a_1:=k_1\|Q\|^{3}_{L^3}>0$ and $k_1:= \displaystyle\frac{\sin\left(\frac{\pi}{2}\alpha\right)}{\pi} \displaystyle\int_{0}^{+\infty}e^{- r^{\frac{1}{\alpha}}} dr$ .
    \item (\textit{Higher order expansion})
        There exists $C>0$ such that:
    \begin{align}
    \left\vert Q(x) -\left( \frac{a_1}{x^{\alpha+1}} + \frac{a_2}{x^{2\alpha+1}} + \frac{a_3}{x^{\alpha+3}} \right)\right\vert & \leq \frac{C}{x^{3\alpha+1}} ,\label{asympt:Q}\\
    \left\vert Q'(x) +(\alpha+1)\frac{a_1}{x^{\alpha+2}} + (2\alpha+1) \frac{a_2}{x^{2\alpha+2}} \right\vert &  \leq \frac{C}{x^{3\alpha+1}}  , \label{asympt:Q'}\\
    \left\vert \Lambda Q(x) +  \frac{a_1(\alpha+2)}{2(\alpha+1)}\frac{1}{x^{\alpha+1}} + \frac{a_2(3\alpha+2)}{2(\alpha+1)}\frac{1}{x^{2\alpha+1}} \right\vert & \leq \frac{C}{x^{\alpha+3}} \label{asympt:LambdaQ}.
    \end{align}
    with 
    $a_2:=k_2\|Q\|^{3}_{L^{3}}$, $ k_2:=-\displaystyle\frac{2\sin\left(\pi\alpha \right)}{\pi} \displaystyle\int_{0}^{+\infty}re^{- r^{\frac{1}{\alpha}}}dr$, and $a_3\in\R$. 
\end{enumerate}
\end{theo}

We also recall some results of regularity given by convolution with the kernel $k$ associated to the dispersion:
\begin{align*}
    k(x):= \int_\R \frac{e^{ix\xi}}{1+\vert \xi \vert^\alpha} d\xi. 
\end{align*}
\begin{lemm}[\cite{EV22}]\label{lemma:est:kernel}
Let $g\in X^{0}(\R)$. There exists $C=C(g)$ such that: 
\begin{align*}
    |k\ast g|(x)\leq \frac{C}{\langle x \rangle^{1+\alpha}}.
\end{align*}

Furthermore, if $g\in \mathcal{C}^1(\R)$, and $\vert g'(x) \vert \leq C \langle x \rangle^{-2-\alpha}$, then there exists $C=C(g,g')$ such that: 
\begin{align*}
    \vert \partial_x \left( k \star g\right)\vert (x) \leq \frac{C}{\langle x \rangle^{2+\alpha}}.
\end{align*}
\end{lemm}

We set the expansion of the translated ground state $Q(x+z)$ at $+\infty$ in $x$ by:
\begin{align*}
    Q_{app}(x,z):=\frac{a_1}{z^{\alpha+1}}- (\alpha+1)a_1\frac{x}{z^{\alpha+2}} +\frac{a_2}{z^{2\alpha+1}} +\left( a_1 \frac{(\alpha+1)(\alpha+2)}{2}x^2 + a_3 \right) \frac{1}{z^{\alpha+3}}.
\end{align*} 

\begin{lemm}
Let $z$ be large enough. We have for all $|x|\leq \frac{z}{2}$:  
\begin{align}\label{eq:estimate_Qapp}
    |Q_{app}(x,z)-Q(x+z)|+|Q_{app}(-x,z)-Q(x-z)| & 
    \leq C \left( \frac{\vert x \vert ^3}{z^{\alpha+4}}+ \frac{\vert x \vert }{z^{2\alpha+2}}+ \frac{1}{z^{3\alpha+1}} \right), \\
    \vert \partial_x Q_{app}(x,z) - Q'(x+z) \vert + \vert \partial_x Q_{app} (-x,z) -Q'(x-z)\vert & \leq C \left( \frac{x^2}{z^{\alpha+4}}+ \frac{1}{z^{2\alpha+2}} \right) \label{eq:estimate_Qapp'} \\
    \left\vert \Lambda Q(x+z) + \frac{a_0(\alpha+2)}{2(\alpha+1)}\frac{1}{z^{\alpha+1}} \right\vert + \left\vert \partial_x \Lambda Q(x+z) \right\vert & \leq C \left( \frac{\vert x \vert}{z^{2+\alpha}} +\frac{1}{z^{2\alpha+1}} \right)
\end{align}
\end{lemm}

\begin{proof}
From the asymptotic of $Q$ in \eqref{asympt:Q} and the asymptotic expansions:
\begin{align*}
    \left\vert \frac{a_1}{\vert x- z \vert^{\alpha+1}} -\left( \frac{a_1}{z^{\alpha+1}} - a_1(\alpha+1) \frac{x}{z^{\alpha+2}} + a_1\frac{(\alpha+1)(\alpha+2)}{2} \frac{x^2}{z^{\alpha+3}} \right) \right\vert\leq C \frac{\vert x \vert^{3}}{z^{\alpha+4}}
\end{align*}
 and the ones of $\displaystyle \frac{a_2}{\vert x-z\vert^{2\alpha+1}}$ and $\displaystyle \frac{a_3}{\vert x - z\vert ^{\alpha+3}}$, we get the development of $Q(x+z)$.
The proof is similar for $Q'$ with \eqref{asympt:Q'}.

The proof of $\Lambda Q$ is a combination of the two previous asymptotic expansions.
\end{proof}

\begin{propo}\label{propo:DL_Q_mu}
Let $\mu^*>0$ be small enough. There exists a constant $C>0$, such that for any $\mu \leq \mu^*$, we have:
\begin{align}\label{eq:lambda_Q_DL}
    \left\vert Q_{1+\mu} - Q - \mu \Lambda Q\right\vert +
    \left\vert Q_{1+\mu}^2 -Q^2 - 2\mu Q \Lambda Q \right\vert \leq C\frac{\mu^2}{\langle x \rangle^{1+\alpha}}.
\end{align}
The following terms are also bounded in terms of $\mu$:
\begin{align}\label{eq:lambda_Q_H_1}
    \left\| Q_{1+\mu} -Q - \mu \Lambda Q \right\|_{H^2} \leq C \mu^2 \quad \text{and} \quad \left\| \Lambda Q_{1+\mu} - \Lambda Q \right\|_{H^1} \leq C\mu.
\end{align}
Moreover, the scalar product of $Q$ with $\Lambda Q$ is:
    \begin{align}\label{eq:Q_scalaire_LambdaQ}
        \langle Q, \Lambda Q \rangle = \frac{\alpha-1}{2(\alpha+1)} \| Q \|_{L^2}^2.
    \end{align}
\end{propo}

\begin{proof}
By the Taylor formula in $\mu$, we have, with \eqref{defi:lambda_lambda_Q}, \eqref{asympt:Q} and \eqref{eq:decay_Q_j} for the second derivative:
\begin{align*}
    Q_{1+\mu} - Q - \mu \Lambda Q & = \int_{1}^{1+\mu} (1+\mu-s) \Lambda^2 Q_s ds, \\
    \left\vert Q_{1+\mu} - Q - \mu \Lambda Q \right\vert & \leq \int_{0}^\mu \frac{\mu-s}{(1+s)^2} \frac{1}{\langle x \rangle^{\alpha+1}_{1+s}} ds \leq C\frac{ \mu^2}{\langle x \rangle^{\alpha+1} }.
\end{align*}
The proof is similar for $Q_{1+\mu}^2$.

Notice that the previous bound still holds for two more derivatives, and the integral gives the first part of \eqref{eq:lambda_Q_H_1}. The second part is similar.
\end{proof}

\subsection{Properties of the linearized operator}

We recall some results on the spectrum of the linearized operator $L$ and establish new inversion lemma on $L$.

\begin{theo}[\cite{weinstein1987existence,albert1997model,kenig2011local,FL13}]\label{theo:L}
	Let $\alpha\in]1,2[$. There exists $Q\in H^{\frac{\alpha}{2}}(\R)\cap C^{\infty}(\R)$ such that  
	\begin{enumerate}
		\item (\textit{Linearized operator}) Let $L$ be the unbounded operator defined on $L^{2}(\R)$ by:
		\begin{align}
		Lv=|D|^{\alpha}v+v-3Q^2v.\label{oplinea}
		\end{align}
		Then, the continuous spectrum of $L$ is $[1,+\infty[$, $L$ has one negative eigenvalue $\mu_0$, associated to an even eigenfunction $v_0>0$, and $\text{ker} \, L= \text{span} \, \{ Q'\}$.
		\item (Invertibility) For any $g\in L^2(\R)$ orthogonal to $v_0$ and $Q'$, there exists a unique $f\in L^2(\R)$ such that $Lf=g$ and $f\perp Q'$. Furthermore, if $g\in H^k(\R)$, then $f\in H^{k+\alpha}(\R)$.
	\end{enumerate}
	\end{theo} 
	
\begin{proof}
We give the proof of the second point. By the Lax-Milgram theorem on $H^{\frac{\alpha}{2}}(\R)$, we obtain the existence of $f$ in the same space. Because $f$ satisfies $\vert D \vert^\alpha f = g -f+3Q^2 f$, we have $f\in H^\alpha(\R)$.

Concerning the higher regularity of $g$, if $f$ is solution of $Lf=g$ with $g\in H^{k}(\R)$, then, since $[\partial_y,L]v =3\partial_y(Q^2)v$ for all $v\in \mathcal{S}(\R)$, we obtain that $f\in H^{k+\alpha}(\R)$.
\end{proof}

\begin{rema}\label{remark:inv_L}
    From Theorem \ref{theo:L} we have the operator $L$ verifies that there exists $\kappa>0$ such that for all $f\in H^{\frac{\alpha}{2}}(\R)$, with $f\perp v_0,Q'$ then: 
    \begin{align*}
        \langle Lf,f \rangle\geq \kappa\|f\|_{H^{\frac{\alpha}{2}}}^2.
    \end{align*}
    However, it is not convenient to work with $v_0$. An argument of Weinstein will allow us to replace the orthogonality on $v_0$ by an orthogonality on $Q$ to get the coercivity. Indeed, from Lemma $E.1$ in \cite{weinstein1985modulation} and since 
    \begin{align*}
        \langle L^{-1}Q,Q\rangle=-\langle \Lambda Q,Q\rangle = -\frac{\alpha-1}{2(\alpha+1)} \| Q \|_{L^2}^2<0. 
    \end{align*}
    we obtain the coercivity of $L$ up to the orthogonality condition on $Q$ and $Q'$:
    \begin{align}\label{eq:coercivite}
        \forall f \in H^{\frac{\alpha}{2}}(\R), \quad f\perp Q, Q' \quad \text{implies} \quad \langle Lf, f \rangle \geq \kappa \| f \|_{H^\frac{\alpha}{2}}^.
    \end{align}
\end{rema}

We continue this section with two lemmas on the characterisations on the inverse of particular functions by $L$ on specific directions.

\begin{lemm}\label{lemm:ant1}
Let $k>0$ and $g\in X^k(\R)$ with $g\perp Q'$, then there exist a unique $f\in X^{k+\alpha}(\R)$, $a\in \R$ such that: $$\begin{cases}
 Lf= g+ a Q\\
f\perp Q,\quad f\perp Q'
\end{cases}.$$
\end{lemm}

\begin{proof}
Since $g+ a Q\perp Q'$, we apply the invertibility property of Theorem \ref{theo:L} and there exists a unique $f\in H^{k+\alpha}(\R)$ such that:
$$\begin{cases}
Lf= g+ a Q\\
f\perp \ker(L)=\text{span}(Q')
\end{cases}.
$$

To obtain the second orthogonality condition, since $L\Lambda Q=-Q$, with \eqref{eq:Q_scalaire_LambdaQ} we deduce that: 
\begin{align}
    \langle f,Q \rangle=0 \iff \langle g+ aQ, \Lambda Q\rangle=0 \iff a=-\frac{\langle g,\Lambda Q\rangle}{\langle  Q,\Lambda Q\rangle}=-\frac{2(\alpha+1)}{\alpha-1}\frac{\langle g,\Lambda Q\rangle}{\|Q\|_{L^2}^2},
\end{align}

We finish with the decay in $\langle x \rangle^{-1-\alpha}$ from the definition \eqref{defi:X} of $X^{k+\alpha}(\R)$. Since $g+ a Q+3Q^2 f\in X^{k}(\R)$, we obtain by Lemma \ref{lemma:est:kernel} that $f=(|D|^{\alpha}+1)^{-1}(g+ a Q+3Q^2 f)\in X^{k+\alpha}(\R)$. This concludes the proof of  Lemma \ref{lemm:ant1}.
\end{proof}

We define a function $S_0$ such that the $\partial_y L S_0$ is close to $\Lambda Q$, in the sense that the remaining terms are of the form $\partial_y (g)$ for some function $g$:
\begin{align}\label{definition:S_0}
    S_0(y) := \int_y^{+\infty} \left( \vert D \vert^\alpha+1 \right)^{-1} \Lambda Q( \tilde{y} ) d\tilde{y}.
\end{align}
$S_0$ is a well-defined function. It has a limit at $-\infty$, which may be different from $0$ and is denoted by $l$:
\begin{align}\label{defi:l}
    l:= \lim_{y \rightarrow -\infty} S_0(y).
\end{align}
See Appendix \ref{proof:lemma:S0} for the justification of $S_0$.

\begin{lemm}\label{lemm:ant2}
Let $g\in X^k(\R)$. There exist a unique $a,\tilde{a}\in \R$ and a unique function $f\in X^{k+\alpha}(\R)$ such that: $$\begin{cases}
\partial_yL(f-\tilde{a}S_0)=  \partial_y g + a  Q' + \tilde{a}\Lambda Q\\
f-\tilde{a}S_0\perp Q,\quad f-\tilde{a}S_0\perp Q'
\end{cases}$$
with
\begin{align}\label{defi:tilde_a}
a=-\displaystyle\frac{2(\alpha+1)}{\alpha-1}\frac{\langle g-\tilde{a}(\vert D \vert^\alpha +1) S_0 ) ,\Lambda Q\rangle}{\|Q\|_{L^2}^2} \quad \text{ and } \quad \tilde{a}=\displaystyle\frac{2(\alpha+1)}{\alpha-1} \frac{\langle g, Q'\rangle}{\|Q\|_{L^2}^2}.
\end{align}

Similarly, there exist a unique $a$, $\tilde{a} \in \mathbb{R}$ and a unique function $f\in X^{k+\alpha}(\R)$ such that: $$\begin{cases}
\partial_yL\left(f+\tilde{a}(l-S_0)\right)=  \partial_y g + a  Q' + \tilde{a}\Lambda Q\\
f+\tilde{a}(l-S_0)\perp Q,\quad f+\tilde{a}(l-S_0)\perp Q'
\end{cases}$$
with
\begin{align*}
a=-\displaystyle\frac{2(\alpha+1)}{\alpha-1}\frac{\langle g+\tilde{a}(\vert D \vert^\alpha+1)(l-S_0)) ,\Lambda Q\rangle}{\|Q\|_{L^2}^2} \quad \text{ and } \quad \tilde{a}=\displaystyle\frac{2(\alpha+1)}{\alpha-1} \frac{\langle g, Q'\rangle}{\|Q\|_{L^2}^2}.
\end{align*}

\end{lemm}

\begin{proof}

We denote by $\mathcal{H}$ the Hilbert transform. Since $|D|^{\alpha}=|D|^{\alpha-1}\mathcal{H}\partial_y$, we deduce that:
\begin{align*}
    |D|^{\alpha}S_0=|D|^{\alpha-1}\mathcal{H}(|D|^{\alpha}+1)^{-1}\Lambda Q= \int_{y}^{+\infty} |D|^{\alpha}(|D|^{\alpha}+1)^{-1}\Lambda Q.
\end{align*}
Then, we get that: 
\begin{align}
    \partial_y LS_0(y)= -\Lambda Q(y) -3\partial_y\left(Q^2(y)S_0(y)\right).
\end{align}
Therefore, it is enough to prove that the following problem has a unique solution:
$$\begin{cases}
Lf=   g + a  Q - \tilde{a} 3Q^2S_0 \\
f-\tilde{a}S_0\perp Q,\quad f-\tilde{a}S_0\perp Q'
\end{cases}.$$
We choose $\tilde{a}$ such that $g+aQ+\tilde{a}3Q^2S_0$ is orthogonal to $Q'$, and then arguing as in the proof of Lemma \ref{lemm:ant1}, we conclude the proof the first identity of Lemma \ref{lemm:ant2}. The second identity is similar.

\end{proof}

\section{Construction of the approximation}\label{sec:approximation}

The approximation $V$ of the expected solution $u$ is built in this section. The purpose is to minimise the flow $\cE_V$ associated to the approximation, by detailing $V$. By taking the time derivative of the sum of two solitary waves $-R_1+R_2$, a particular direction intrinsic to the problem appears and is compensated by the use of a function $W$. This term possesses a tail at $-\infty$. We also define a time-dependent variable $b(z(t))$. We then minimise the flow associated to $-R_1+R_2+bW$ by adding localised profiles $-P_1$ and $P_2$ in the approximation to cancel the source term coming from the non-linearity.

\subsection{Notation}

Let us consider four $\mathcal{C}^1$ functions $\mu_1$, $\mu_2$, $z_1$ and $z_2$ on a time interval $I \subset \mathbb{R}$, and
\begin{align}
    \Gamma(t):=(\mu_1(t),\mu_2(t),z_1(t),z_2(t)).
\end{align}
We define the distance between the different functions by:
\begin{align}\label{defi:mu_z}
    \mu (t) := \mu_1(t)-\mu_2(t), \quad z(t):=z_1(t)-z_2(t).
\end{align}
For a fixed constant $C_0>0$, we use the following set of assumptions on the interval $I$:
\begin{align}
    -z(t) \leq  z_2(t) \leq -\frac{1}{8} z(t),& \quad  
    \frac{1}{8} z(t) \leq  z_1(t) \leq z(t), \label{hypothese:I} \\
    \left\vert \mu_1(t) \right\vert+\left\vert \mu_2(t) \right\vert+\left\vert \mu(t) \right\vert+ \left\vert \dot{z}_1(t) \right\vert + \left\vert \dot{z}_2(t) \right\vert & + \left\vert \dot{z}(t) \right\vert \leq  \frac{C_0}{z(t)^{\frac{\alpha+1}{2}}}, \label{hypothese:II} \\
    |\dot{\mu}_1(t)|+|\dot{\mu}_2(t)| & \leq \frac{C_0}{z(t)^{\alpha+2}}. \label{hypothese:III}
\end{align}

\begin{rema}
The constant $C_0$ is used to fix the set of assumptions on $\Gamma$. The computations of this section involve the constant $C_0$, but it does not have any influence on the final constant $C$ in Theorem \ref{main_theo}. For the sake of simplicity, we omit the presence of this constant in the computations. To close the bootstrap in subsection \eqref{subsec:bootstrap_setting}, we will fix the constant $C_0$ to be large enough so that the set assumptions on $\Gamma$ is satisfied.
\end{rema}
We define a function 
\begin{align}\label{defi:b}
    b(z(t)):= \frac{b_1}{z^{\alpha+2}(t)},
\end{align} 
with $\displaystyle b_1=-2a_1\frac{(\alpha+1)^2}{\alpha-1}\frac{\|Q\|_{L^3}^3}{\|Q\|_{L^2}^2}<0$ and $a_1>0$ defined in Theorem \ref{thm:asympQ}.

\subsection{Approximate solution}

\paragraph{}

To quantify the interaction between the two solitary waves, we set the bump function $W$ by:
\begin{align}\label{definition:W}
     W(\Gamma(t),y):= S_0(y-z_1(t))-S_0(y-z_2(t)),
\end{align}
where $S_0$ has been defined in \eqref{definition:S_0}.

We use the convenient notation of an index $i\in\{1,2\}$ to underline that $R_i$ and $P_i$ are functions centred at $z_i$.

\begin{theo}\label{theo:construction}
    Let $I\subset\R$ an interval such that the assumptions \eqref{hypothese:I}-\eqref{hypothese:III} on $\Gamma$ are satisfied.

There exist two constants $\beta_0$ and $\delta_0$ in $\mathbb{R}$, two functions $\beta(\Gamma)$ and $\delta(\Gamma)$ and two functions $P_1(\Gamma,y)$ and $P_2(\Gamma,y)$ such that the following holds:

\begin{itemize}
    \item Asymptotic of $\beta$ and $\delta$. The functions $\beta$ and $\delta$ have the following expansion:
\begin{align}\label{ineq:beta_delta}
    \left\vert \beta(\Gamma) - \frac{\beta_0}{z^{1+\alpha}}\right\vert + \left\vert \delta(\Gamma) -\frac{\delta_0}{z^{1+\alpha}}\right\vert \leq \frac{C}{z^{2+\alpha}}.
\end{align}
    \item Orthogonality conditions and limits. The profiles $P_i(\Gamma) \in \mathcal{C}(I, X^{2+\alpha}(\R))$ satisfy:
\begin{align}
    -P_1 +b(z)S_0(\cdot-z_1) \perp \tilde{R}_1, \partial_y \tilde{R}_1, \quad 
    P_2 +b(z)(l-S_0(\cdot-z_2)) \perp \tilde{R}_2, \partial_y \tilde{R}_2.
    \end{align}
    We then define the approximation $V$ of a solution by:
\begin{align}\label{defi:V}
    V(\Gamma,y):=\displaystyle\sum_{i=1}^2(-1)^{i}\left( R_i(\Gamma,y)+P_i(\Gamma,y)\right) + b(z)W(\Gamma,y),
\end{align}
and for simplicity we will write $V(t,y):=V(\Gamma(t),y)$.

\item Decomposition and estimate of the flow.
The flow $\cE_V$ of the approximation 
\begin{align}\label{defi:E_V}
    \cE_{V} := \partial_t V + \partial_y \left( - \vert D \vert^\alpha V -V +V^3 \right) 
\end{align}
can be decomposed into:
\begin{align}\label{definition:energy}
    \cE_V =\overrightarrow{m}\cdot\overrightarrow{MV} + \partial_yS + T 
\end{align}
with 
\begin{align}\label{defi:m_MV}
    \overrightarrow{m}(t)=\begin{pmatrix}
    -\dot{\mu}_1(t) +b(z(t)) \\
    \dot{z_1}(t) -\mu_1(t) +\beta(\Gamma(t)) \\
    \dot{\mu_2}(t) +b(z(t)) \\
    -\dot{z_2}(t)+\mu_2(t) -\delta(\Gamma(t))
    \end{pmatrix}, \quad 
    \overrightarrow{MV}(t,y)=\begin{pmatrix}
    \Lambda R_1(t,y)\\
    \partial_y R_1(t,y)\\
    \Lambda R_2(t,y)\\
    \partial_y R_2(t,y)
    \end{pmatrix},
\end{align}
 and the source term $S$ and the approximation due to the flow $T$ are in $\mathcal{C}^1(I, X^{2+\alpha}(\R))$ and satisfy the set of inequalities :
\begin{align}
    \|S\|_{H^1} & \leq \frac{C}{z^{\frac{5+3\alpha}{2}}},\label{est:S}\\
    \|\partial_tS\|_{L^2} & \leq \frac{C}{z^{3+2\alpha}}\label{est:dtS}, \\
    \|T\|_{H^1} & \leq \frac{C}{z^{\frac{5+3\alpha}{2}}} + \frac{C}{z^{1+\alpha}}\sum_{i=1}^2 \vert \dot{z}_i -\mu_i \vert. \label{est:T}
\end{align}
\end{itemize}

\end{theo}

We add some estimates related to the previously defined functions. We define two functions $\phi$ and $\Phi$ (that correspond to \eqref{defi:phi} and \eqref{Phi'}) by:
\begin{align}
    \phi(y)=\left(\int_{-\infty}^{+\infty} \frac{ds}{\langle s \rangle^{1+\alpha}}\right)^{-1}&\int_{y}^{+\infty} \frac{ds}{\langle s \rangle^{1+\alpha}} \quad \text{and} \quad  \Phi(y)=\sqrt{|\phi'(y)|}.
\end{align}

\begin{propo}\label{propo:estimates_flow}
With the previous notations, the following estimates hold:
\begin{itemize}
    \item Estimates on the solitary waves:
    \begin{align}
    \|R_i(\phi-\delta_{2i})\|_{H^1} + \left\|\partial_y R_i  \left(\phi-\delta_{2i} \right) \right\|_{H^1}
        +\| (1- \sqrt{\vert \delta_{1i}-\phi\vert})R_i) \|_{L^2} & \leq \frac{C}{z^{\alpha}} ,\quad i=1,2\label{eq:partial_y_Ri},\\
        \|\partial_yR_i\Phi\|_{L^2} + \|\Lambda R_i\Phi\|_{L^2}
        & \leq \frac{C}{z^{\frac{1+\alpha}{2}}},\quad i=1,2\label{eq:Ri_Phi},
    \end{align}
where $\delta_{ij}$ holds for Kronecker delta.
    \item Estimates on the profiles:
\begin{align}
    \|\left(P_1-P_2-bW \right)\|_{L^{\infty}} +
    \|\partial_y\left(P_1-P_2-bW \right)\|_{L^{\infty}}
        & \leq \frac{C}{z^{1+\alpha}},\label{partialyPibWchi}\\
    \|\partial_t\left(P_1-P_2-bW \right)\|_{L^{\infty}}
        & \leq \frac{C}{z^{\frac{3+3\alpha}{2}}},\label{partialtPibWchi}
\end{align}
    \item Estimates on the approximation:
\begin{align}
    \| V\|_{L^{\infty}}+\|\partial_y V\|_{L^{\infty}}
        & \leq C \label{V},\\
    \| V^{k} \Phi^2\|_{L^{\infty}} + \| (V^2 -R_1^2) \partial_y R_1 \|_{L^2}
        & \leq \frac{C}{z^{1+\alpha}} \quad k\in\N, \label{Vphi}\\
    \|\partial_t V \|_{L^{\infty}}
        & \leq \frac{C}{z^{\frac{1+\alpha}{2}}}. \label{partialtV}
\end{align}
\end{itemize}
\end{propo}

The next subsections are dedicated to the proof of the theorem on the approximation $V$. We begin with the expansion of $\cE_V$ defined in \eqref{defi:E_V}. Let us first compute the different time derivatives:
\begin{align*}
    \partial_t (- R_1) = \dot{z}_1 \partial_y R_1 - \dot{\mu}_1 \Lambda R_1 \quad \text{and} \quad \partial_t R_2 = -\dot{z}_2 \partial_y R_2 + \dot{\mu}_2 \Lambda R_2.
\end{align*}
By the definition of $V$ in \eqref{defi:V}, we get the development:
\begin{align}
    \cE_V
        & = \sum_{i=1}^2 (-1)^i \left( \dot{\mu}_i \Lambda R_i - \dot{z}_i \partial_y R_i \right) +\sum_{i=1}^2 (-1)^i \partial_y \left( -\vert D \vert^\alpha R_i-R_i +R_i^3\right)  \label{eq:decomposition_E_V_MV} \\
        & \quad +\sum_{i=1}^2 \partial_y \left( \left(- \vert D \vert^ \alpha-1 +3R_i^2 \right) ((-1)^i P_i) \right) + \partial_y \left( (-\vert D \vert^{\alpha} -1)(bW)\right) \notag \\
        & \quad+ \partial_y (V^3 + R_1^3 - R_2^3 +3 R_1^2 P_1- 3R_2^2 P_2) \label{eq:decomposition_E_V_S}\\
        & \quad + \frac{d}{dt}(-P_1+P_2+bW). \label{eq:decomposition_E_V_T}
\end{align}

The objective is to decompose $\cE_V$ in order to get the decomposition \eqref{definition:energy}, and to justify the definition of $S$ and $T$.

We now decompose the interaction term \eqref{eq:decomposition_E_V_S} coming from the non-linearity in order to get the first part $S_{V^3}$ of $S$, keeping in mind that the objective is to achieve the bound \eqref{est:S} on $S$ to obtain an accurate approximation of the $2$-solitary waves.

Since the natural space to study the approximation is the energy space, we look at the $L^2$-norm of the different terms, and extract from it the larger ones in terms of $z$. The next computations give insights of the decomposition made throughout the article and are formal, but the estimates are rigorously proved along the next subsections. Let us quantify the different terms involved in $V^3$:
\begin{itemize}
    \item The terms at the main order:
    \begin{align*}
        \left\| R_1^3\right\|_{L^2} +\left\| R_2^3 \right\|_{L^2} \simeq 1.
    \end{align*}
    \item The terms coming from the interaction. Note that the objects indexed by $i\in\{1,2\}$ are localised in $z_i$. Therefore, the product of two functions indexed by different values is "small" due to the different support of the functions:
    \begin{align*}
        \left\|3R_1^2R_2\right\|_{L^2} + \left\|3R_1R_2^2\right\|_{L^2} \simeq \frac{1}{z^{\alpha+1}}.
    \end{align*}
    We want to construct the profiles $P_1,P_2$,  respectively centred at $z_1,z_2$, with at least a decay of $y^{-\alpha-1}$ and that compensates the interaction terms above: 
    \begin{align*}
        \left\| \langle y -z_1 \rangle ^{\alpha+1} P_1 \right\|_{L^\infty} \simeq \frac{1}{z^{\alpha+1}}, \quad \left\| \langle y -z_2 \rangle ^{\alpha+1} P_2 \right\|_{L^\infty}\simeq \frac{1}{z^{\alpha+1}}.
    \end{align*}
    The process of construction of the profiles $P_1$ and $P_2$ involves an orthogonality condition that is not always satisfied. To this aim, we use the non-localized function $W$ and a coefficient $b(z)$ to get these orthogonality conditions. Consequently, the function $bW$ is expected to be at the same order as $P_1$ and $P_2$.
    
    From the decay property of the profiles $P_1,P_2$, we expect:  
    \begin{align*}
        \left\| 3R_1^2 P_1\right\|_{L^2} + \left\| 3R_2^2 P_2\right\|_{L^2} \simeq \frac{1}{z^{\alpha+1}},
    \end{align*}
    and by the definition \eqref{definition:W} of $W$, we have:
    \begin{align*}
        \left\| 3 R_1^2 b(z) S_0(\cdot-z_1)\right\|_{L^2} + \left\| 3 R_2^2 b(z)\left( l- S_0(\cdot-z_2) \right) \right\|_{L^2} \simeq b(z).
    \end{align*}
    \item The remaining terms. The other terms are small enough to get the condition \eqref{est:S}. Using the support of the two functions combined with decay of $P_1,P_2$: 
    \begin{align*}
        \| 3R_1^2 P_2\|_{L^2} + \| 3R_2^2 P_1\|_{L^2} + \|3R_1 P_1^2\|_{L^2}+ \|3R_2 P_2^2\|_{L^2}\simeq \frac{1}{z^{2\alpha+2}},
    \end{align*}
\end{itemize}

It is thus natural to define:
\begin{align}
    S_{V^3} := & V^3+R_1^3-R_2^3+3R_1^2P_1-3R_2^2P_2 -3R_1^2bS_0(y-z_1) - 3R_2^2b(l-S_0(y-z_2))-3R_1^2R_2+3R_1R_2^2 \notag\\
    = & 3R_1^2(P_2-bS_0(y-z_2))+3R_2^2(-P_1+b(S_0(y-z_1)-l))-6R_1R_2(-P_1+P_2+bW) \notag\\
    &  \quad +3(-R_1+R_2)(-P_1+P_2+bW)^2 + (-P_1+P_2+bW)^3. \label{defi:S_V_3}
\end{align}

Notice from the following identities:
\begin{align*}
    \partial_y (-\vert D \vert^\alpha(-R_1)-(-R_1) -R_1^3 ) = -\mu_1 \partial_y R_1 \quad \text{and} \quad \partial_y (-\vert D \vert^\alpha R_2-R_2+ R_2^3 ) = \mu_2 \partial_y R_2,
\end{align*}
and the definition of $W$ in \eqref{definition:W}, we thus rewrite the energy $\cE_V$ as:
\begin{align}
    \cE_V
        & = \sum_{i=1}^2 (-1)^i \left( \dot{\mu}_i \Lambda R_i + \left( -\dot{z}_i + \mu_i \right) \partial_y R_i \right)  \nonumber \\
        & \quad + \partial_y \left( \left(- \vert D \vert^ \alpha-1 +3R_1^2 \right) ( -P_1 +bS_0(\cdot - z_1) ) + 3R_1^2R_2 \right) \label{S_1_app} \\
        & \quad +\partial_y \left( \left(- \vert D \vert^ \alpha-1 +3R_2^2 \right) \left( P_2 +b(l-S_0(\cdot - z_1))\right) - 3R_1R_2^2 \right)  \label{S_2_app} \\
        & \quad + \frac{d}{dt}(-P_1+P_2+bW) + \partial_y (S_{V^3}). \nonumber 
\end{align}

We continue with a second term in the decomposition of $S$, given by the terms \eqref{S_1_app} and \eqref{S_2_app}. 

The estimates involved in the construction are dependent on the parameters $\mu_i$ and $z_i$. In order to explicit this dependency,  we need to "separate" the parameters by applying an asymptotic development in $\mu_i$ of the function $\mu_i\mapsto Q_{1+{\mu_i}}(\cdot - z_i)$ as in Proposition \ref{propo:DL_Q_mu}. With the bound on $\mu_i$ in \eqref{hypothese:II}, we have:   

\begin{align*}
    \left\| 3R_1^2 (-P_1) - 3 \tilde{R}_1^2 (-P_1) - 6\mu_1 \tilde{R}_1 \Lambda \tilde{R}_1 (-P_1) \right\|_{L^2} & \simeq \mu_1^2 \| P_1 \|_{L^2} \simeq \frac{1}{z^{2\alpha+2}}, \\
    \left\| 3R_1^2 bS_0(\cdot - z_1) -3 \tilde{R}_1^2 bS_0(\cdot- z_1) \right\|_{L^2} & \simeq \mu b(z) \simeq \frac{1}{z^{\frac{3\alpha+5}{2}}}, \\
    \left\| 3 R_1^2 R_2 - 3 \tilde{R}_1^2 \tilde{R}_2 -6\mu_1 \tilde{R}_1 \Lambda \tilde{R}_1 \tilde{R}_2 - 3\mu_2 \tilde{R}_1^2 \Lambda \tilde{R}_2 \right\|_{L^2} & \simeq (\mu_1^2 +\mu_2^2) \frac{1}{z^{\alpha+1}} \simeq \frac{1}{z^{2\alpha+2}},
\end{align*}
with $\tilde{R}_i$ defined in \eqref{def:R_tilde}. The computations are similar for $\tilde{R}_2$. We thus set:
\begin{align}\label{defi:S_cal_tilde}
    \forall i\neq j \in \{1,2\}, \quad \widetilde{\mathscr{S}}(i,j) := 6\mu_i \tilde{R}_i \Lambda \tilde{R}_i P_i + 3 R_i^2 R_j - 3 \tilde{R}_i^2 \tilde{R}_j - 6\mu_i \tilde{R}_i \Lambda \tilde{R}_i \tilde{R}_j - 3 \mu_j \tilde{R}_i^2 \Lambda \tilde{R}_j,
\end{align}
and:
\begin{align}\label{defi:S_tilde}
    \widetilde{S}
        & := 3 (R_1^2-\tilde{R}_1^2)(-P_1 + bS_0(\cdot-z_1)) + \widetilde{\mathscr{S}}(1,2) \nonumber\\
            & \quad + 3(R_2^2 - \tilde{R}_2^2)(P_2 + b(l-S_0(\cdot-z_2)))- \widetilde{\mathscr{S}}(2,1).
\end{align}
The energy $\cE_V$ can then be rewritten as:
\begin{align}
    \cE_V
        & = \sum_{i=1}^2 (-1)^i \left( \dot{\mu}_i \Lambda R_i + \left( -\dot{z}_i + \mu_i \right) \partial_y R_i \right) \label{m_dot_MV_app_bis} \\
        & \quad + \partial_y \left( \left(- \vert D \vert^ \alpha-1 +3\tilde{R}_1^2 \right) ( -P_1 +bS_0(\cdot - z_1) ) \right. \nonumber \\
            & \quad \quad \quad \left. + 3\tilde{R}_1^2\tilde{R}_2+6\mu_1 \tilde{R}_1 \Lambda \tilde{R}_1 (-P_1) + 6\mu_1 \tilde{R}_1 \Lambda \tilde{R}_1 \tilde{R}_2 + 3\mu_2 \tilde{R}_1^2 \Lambda \tilde{R}_2 \right) \label{S_A_1_app_bis}\\
        & \quad +\partial_y \left( \left(- \vert D \vert^ \alpha-1 +3\tilde{R}_2^2 \right) \left( P_2 +b(l-S_0(\cdot - z_1))\right) \right. \\
            & \quad \quad \quad \left. - 3\tilde{R}_1\tilde{R}_2^2 +6\mu_2 \tilde{R}_2 \Lambda \tilde{R}_2 P_2 - 6\mu_2 \tilde{R}_1 \tilde{R}_2  \Lambda \tilde{R}_2- 3\mu_1 \Lambda \tilde{R}_1 \tilde{R}_2^2 \right) \label{S_A_2_app_bis}\\
        & \quad + \frac{d}{dt}(-P_1+P_2+bW) + \partial_y (S_{V^3} + \tilde{S}). \label{d_dt_P1_app_bis} 
\end{align}

To apply Lemma \ref{lemm:ant2} on \eqref{S_A_1_app_bis} and \eqref{S_A_2_app_bis}, we need to adjust the directions $\partial_y \tilde{R}_1$, $\partial_y \tilde{R}_2$, $\Lambda \tilde{R}_1$ and $\Lambda \tilde{R}_2$ with coefficients $b(z)$, $\beta(\Gamma)$ and $\delta(\Gamma)$. We add the quantity $-b(z)\Lambda \tilde{R}_1-\beta(\Gamma) \partial_y \tilde{R}_1$ on the line \eqref{S_A_1_app_bis}, the quantity $-b(z)\Lambda \tilde{R}_2+\delta(\Gamma) \partial_y \tilde{R}_2$ on the line \eqref{S_A_2_app_bis}, and subtract those quantities in \eqref{m_dot_MV_app_bis}. However, in this new \eqref{m_dot_MV_app_bis}, we once again compare $R_i$ with $\tilde{R}_i$; we need to use the asymptotic development of $R_i$ in terms of $\mu_i$ to quantify this error term. For the direction $\Lambda R_i$ with the definition of $b$ in \eqref{defi:b}, we have:
\begin{align*}
    \left\| b(z) \Lambda \tilde{R}_i - b(z) \Lambda R_i \right\|_{L^2} \simeq b(z) \vert \mu_i \vert \simeq \frac{1}{z^{\frac{3\alpha+5}{2}}}.
\end{align*}
If we apply the same computation for the direction $\partial_yR_1$ with the estimate on $\beta$ in \eqref{ineq:beta_delta}, we get:
\begin{align*}
    \left\| \beta(\Gamma) \partial_y \tilde{R}_1 - \beta(\Gamma) \partial_y R_1 \right\|_{L^2} \simeq \beta(\Gamma) \vert \mu_i \vert \simeq \frac{1}{z^{\frac{3\alpha+3}{2}}}.
\end{align*}
However, this quantity is too large compared to \eqref{est:S}. Therefore, we need to adapt the profiles $P_1$ and $P_2$ to control this term. Using a further asymptotic development in $\mu_i$ and \eqref{ineq:beta_delta}, we obtain:
\begin{align*}
    \left\| \beta(\Gamma) \partial_y \tilde{R}_1 - \beta(\Gamma) \partial_y R_1 + \frac{\beta_0\mu_1}{z^{\alpha+1}} \partial_{y}\Lambda \tilde{R}_1 \right\|_{L^2}\simeq \frac{1}{z^{\frac{3\alpha+5}{2}}} + \frac{1}{z^{2\alpha+2}}.
\end{align*}
By applying the same arguments for the direction $\partial_yR_2$, we get:
\begin{align*}
    \left\| \delta(\Gamma) \partial_y \tilde{R}_2 - \delta(\Gamma) \partial_y R_2 + \frac{\delta_0\mu_2}{z^{\alpha+1}} \partial_{y}\Lambda \tilde{R}_2 \right\|_{L^2}\simeq \frac{1}{z^{\frac{3\alpha+5}{2}}} + \frac{1}{z^{2\alpha+2}}.
\end{align*}

The two new terms depending on $\beta_0$ and $\delta_0$ are compensated by adjusting $P_1$ and $P_2$, and thus we add $\displaystyle \partial_y \left( \beta_0 \frac{\mu_i}{z^{\alpha+1}} \Lambda \tilde{R}_1 \right)$ in \eqref{S_A_1_app_bis} and $\displaystyle \partial_y \left( \delta_0 \frac{\mu_2}{z^{\alpha+1}} \Lambda \tilde{R}_2 \right)$ in \eqref{S_A_2_app_bis}.

The last approximation comes from the time derivative of the profiles in \eqref{d_dt_P1_app_bis}. To detail the situation, let us give an idea of the construction of $P_1$. By the method of separation of variables, $P_1$ will be decomposed into a function $\overrightarrow{f}$ depending on $\Gamma$ only, multiplied by a function $\overrightarrow{B}$ depending on $y$ only and translated by $z_1$:
\begin{align*}
    P_1(\Gamma,y) = \overrightarrow{f}(\Gamma) \cdot \overrightarrow{B}(y-z_1).
\end{align*}
Taking the time derivative of $P_1$ yields:
\begin{align*}
    \frac{d}{dt} P_1(\Gamma(t),y)= \left( \frac{d}{dt} \overrightarrow{f}(\Gamma(t)) \right) \cdot \overrightarrow{B}(y-z_1(t)) - \dot{z}_1(t) \partial_y P_1(\Gamma(t),y),
\end{align*}
with $\dot{z}_1 \sim \mu_1 \sim z^{-\frac{\alpha+1}{2}}$. Since $P_1$, or equivalently $\overrightarrow{f}$ is of order $z^{-\alpha-1}$ we formally obtain that $\frac{d}{dt}(\overrightarrow{f}) \simeq z^{\frac{-3\alpha-5}{2}}$, which is convenient in \eqref{est:T}. For the second term, we get:
\begin{align*}
    \left\|\dot{z}_1 \partial_y P_1\right\|_{L^2} \simeq \frac{1}{z^{\frac{3\alpha+3}{2}}}.
\end{align*}
Therefore, we also need to compensate this term. However, we do not know if $z$ is a function of class $C^2$. Since:
\begin{align*}
    \left\| -\dot{z_1} \partial_y P_1 + \mu_1 \partial_y P_1 \right\|_{L^2} \simeq \frac{1}{z^{2\alpha+2}},
\end{align*}
it remains equivalent to have $\mu_1\partial_y P_1$ instead of 
$\dot{z}_1\partial_y P_1$.

Considering the two previous approximations, we define $\mathscr{T}$ by:
\begin{align}\label{defi:T_cal}
    \forall i \in \{ 1,2\}, \quad \mathscr{T}(i,\beta_0) :=  \frac{\beta_0}{z^{1+\alpha}} \partial_y \left( \mu_i \Lambda \tilde{R}_i \right) - \mu_i  \partial_y P_i,
\end{align}
the remaining term of the error by:
\begin{align}
    T
        & := \sum_{i=1}^2 b(z) \left( \Lambda \tilde{R}_i - \Lambda R_i \right) +\beta(\Gamma) \left( -\partial_y R_1 + \partial_y \tilde{R}_1 \right) -\delta(\Gamma) \left( -\partial_y R_2 + \partial_y \tilde{R}_2 \right) \notag \\
        & \quad + \mathscr{T}(1,\beta_0) - \mathscr{T}(2,\delta_0) +\frac{d}{dt}\left( -P_1 +P_2+bW\right), \label{defi:T}
\end{align}
and the terms to inverse by, for any $i \neq j \in \{1,2\}$:
\begin{align}\label{defi:S_cal}
    \mathscr{S}(i,j) := 3 \tilde{R}_i^2 \tilde{R}_j - 6\mu_i \tilde{R_i}\Lambda \tilde{R}_i P_i +6\mu_i \tilde{R_i} \Lambda \tilde{R}_i \tilde{R}_j + 3\mu_j \tilde{R}_i^2 \Lambda \tilde{R}_j +\mu_i P_i - \beta_0 \frac{\mu_i}{z^{1+\alpha}}\Lambda \tilde{R}_i.
\end{align}
Thus the energy of the error is given by:
\begin{align}
    \cE_V
        & = \sum_{i=1}^2 (-1)^i \left( \left(\dot{\mu}_i + (-1)^ib  \right) \Lambda R_i + \left( -\dot{z}_i + \mu_i \right) \partial_y R_i \right) + \beta \partial_yR_1 - \delta\partial_yR_2 \\
        & \quad + \partial_y \left( \left(- \vert D \vert^ \alpha-1 +3\tilde{R}_1^2 \right) ( -P_1 +bS_0(\cdot - z_1) ) +\mathscr{S}(1,2) \right) -b\Lambda \tilde{R}_1-\beta \partial_y \tilde{R}_1 \\
        & \quad +\partial_y \left( \left(- \vert D \vert^ \alpha-1 +3\tilde{R}_2^2 \right) \left( P_2 +b(l-S_0(\cdot - z_1))\right) - \mathscr{S}(2,1) \right) -b\Lambda \tilde{R}_2+\delta \partial_y \tilde{R}_2\\
        & \quad+ \partial_y (S_{V^3} + \tilde{S}) + T. 
\end{align}

To shorten the notations, we use the definition of $\overrightarrow{m}$ and $\overrightarrow{MV}$ in \eqref{defi:m_MV}, and the functions $S_1$ and $S_2$, equal to $0$ at $+\infty$ and satisfying:
\begin{align}\label{defi:S_1}
    \partial_y S_1 & := \partial_y \left( \left( - \vert D \vert^\alpha - 1 +3 \tilde{R}_1^2\right)\left(-P_1 + b(z)S_0(\cdot-z_1)\right) + \mathscr{S}(1,2) \right) - b(z) \Lambda \tilde{R}_1 - \beta(\Gamma) \partial_y \tilde{R}_1, \\
    \partial_y S_2 & := \partial_y \left( \left( - \vert D \vert^\alpha - 1 +3 \tilde{R}_2^2\right) \left(P_2 + b(z)(l-S_0(\cdot-z_2))\right) - \mathscr{S}(2,1) \right) - b(z) \Lambda \tilde{R}_2 + \delta(\Gamma) \partial_y \tilde{R}_2. \quad \quad  \label{defi:S_2}
\end{align}

By an adequate choice of $P_i$, $b$, $\beta$ and $\delta$, the functions $S_1$ and $S_2$ will not have a tail at $-\infty$, see \eqref{S_1_without_tail} and \eqref{S_2_without_tail}.

We conclude with the following decomposition:
\begin{align*}
        \cE_V = \overrightarrow{m}\cdot\overrightarrow{MV} + \partial_yS + T,
\end{align*}
where $T$ is defined in \eqref{defi:T}, and $S$ by
\begin{align}
    S=S_1+S_2+S_{V^3}+\tilde{S}.
\end{align}

Let us continue the construction in the next subsection by the choices of $P_1$ and $P_2$.

\subsection{Construction of the profiles}

\paragraph{}

This part is dedicated to the construction of the profiles $P_1$ and $P_2$. The goal is to minimise the quantities $S_1$ and $S_2$ by exploiting the intrinsic directions of the problems $\partial_y R_1$, $\partial_y R_2$, $\Lambda R_1$ and $\Lambda R_2$. In particular, the coefficient $b$ defined in \eqref{defi:b} is central in the study of the interaction. The profiles $P_i$ are established term by term in $\mathscr{S}$ in \eqref{defi:S_cal}, and using the expansion of the interaction terms, given by $Q_{app}$.

Due to the definition of $\mathscr{S}$ in \eqref{defi:S_cal}, for any $i\neq j \in \{1,2\}$, we define an approximate value of the function $\mathscr{S}$, where it is located, by:
\begin{align}
    \mathscr{F}(i,j,\beta_0,B_0) 
        & := 3\tilde{R}_i^2 Q_{app}((-1)^j(\cdot-z_i),z) -6\frac{\mu_i}{z^{1+\alpha}} \tilde{R_i} \Lambda \tilde{R}_i B_0(\cdot-z_i) + 6\frac{\mu_i}{z^{1+\alpha}} \tilde{R}_i \Lambda \tilde{R}_i a_1 \notag \\
            & \quad - 3\frac{\mu_j}{z^{1+\alpha}} \tilde{R}_i^2 \frac{a_1(\alpha+2)}{2(\alpha+1)} + \frac{\mu_i}{z^{1+\alpha}} B_0(\cdot-z_i)- \beta_0 \frac{\mu_i}{z^{1+\alpha}} \Lambda \tilde{R}_i.
\end{align}

The definitions of $b$, $S_0$ and $l$ are respectively given in \eqref{defi:b}, \eqref{definition:S_0} and \eqref{defi:l}. 

\begin{propo}\label{propo:est:profils} 
There exist two constants $\beta_0$ and $\delta_0$, two functions $\beta(\Gamma)$ and $\delta(\Gamma)$ in $\mathcal{C}^1(I)$ satisfying  \eqref{ineq:beta_delta}, two even functions $B_0,D_0 \in X^\infty(\R)$
and two profile functions $P_1(\Gamma,y)$ and $P_2(\Gamma,y)$ in $\mathcal{C}(I, X^{\infty}(\R))$ satisfying:
\begin{align}
   \left\vert P_1(\Gamma,y+z_1) - \frac{\beta_0}{z^{1+\alpha}}B_0(y) \right\vert + \left\vert P_2(\Gamma,y+z_2) - \frac{\delta_0}{z^{1+\alpha}}D_0(y) \right\vert & \leq \frac{C}{z^{2+\alpha}}\frac{1}{\langle y \rangle^{1+\alpha}}, \label{iden:prem_ordre_P_i} \\
    \partial_y \left( \left( -\vert D \vert^\alpha -1 + 3 \tilde{R}_1^2 \right) \left( -P_1 +b(z)S_0(\cdot -z_1) \right) + \mathscr{F}(1,2,\beta_0,B_0) \right) & = b(z) \Lambda \tilde{R}_1 + \beta(\Gamma) \partial_y \tilde{R}_1,\quad \quad\label{iden:P_1} \\
    \partial_y \left( \left( -\vert D \vert^\alpha -1 + 3 \tilde{R}_2^2 \right) \left( P_2 +b(z)(l-S_0(\cdot -z_2)) \right) - \mathscr{F}(2,1,\delta_0, D_0) \right) & = b(z) \Lambda \tilde{R}_2 - \delta(\Gamma) \partial_y \tilde{R}_2,\quad\quad \label{iden:P_2}
\end{align}
with the orthogonality conditions:
\begin{align}
    P_1 -b(z)S_0(\cdot- z_1) \perp \tilde{R}_1, \partial_y \tilde{R}_1, \quad \text{and} \quad P_2 +b(z)\left(l-S_0(\cdot- z_1)\right) \perp \tilde{R}_2, \partial_y \tilde{R}_2.
\end{align}

Moreover, the profiles $P_1$, $P_2$ verify:
\begin{align}\label{derivee:d_yP_i}
    \vert P_i(\Gamma,y) \vert + |\partial_yP_i(\Gamma,y)| & \leq \frac{C}{z^{1+\alpha}}\frac{1}{\langle y-z_i \rangle^{1+\alpha}}, \\
    \left\vert\frac{d}{dt} P_i(\Gamma,y) \right\vert & \leq  \frac{C}{z^{\frac{3+3\alpha}{2}}}\frac{1}{\langle y-z_i \rangle^{1+\alpha}}, \label{derivee:d_tP_i} \\
    \left\vert\frac{d}{dt} P_i(\Gamma,y) + \dot{z}_i \partial_y P_i(\Gamma,y) \right\vert & \leq  \frac{C}{z^{\frac{5+3\alpha}{2}}}\frac{1}{\langle y-z_i \rangle^{1+\alpha}}. \label{derivee:d_tP_i_asymptotic}
\end{align}
\end{propo}

The \textbf{\textit{profiles}} $P_1$ and $P_2$ are defined by:
\begin{align}\label{def:profile}
    P_1(\Gamma(t),y):= \overrightarrow{f}\left(\Gamma(t)\right) \cdot \overrightarrow{B}(y-z_1(t)), \quad  P_2(\Gamma(t),y):= \overrightarrow{f}\left(\Gamma(t)\right) \cdot\overrightarrow{D}(y-z_2(t)),
\end{align}
where the functions $\overrightarrow{f}$, $\overrightarrow{B}$ and $\overrightarrow{D}$ are established in the next proposition. The proof of Proposition \ref{propo:est:profils} is postponed after the proof of the next proposition.

\begin{propo}\label{propo:profils}
Let us define the vector functions:
\begin{align}
    \overrightarrow{f}(\Gamma):=\left( \frac{1}{z^{1+\alpha}}, \frac{1}{z^{2+\alpha}}, \frac{\mu_1}{z^{1+\alpha}}, \frac{\mu_2}{z^{1+\alpha}}, \frac{1}{z^{2\alpha+1}}, \frac{1}{z^{3+\alpha}}\right),
\end{align}
and for all $i \in \{1,2\}$:
\begin{align}
    \overrightarrow{F}(i,\beta_0,B_0) & := \left( 3Q^2 a_1, \quad 3Q^2a_1(-1)^i(\alpha+1)y, \quad -6Q\Lambda Q B_0+6Q\Lambda Q a_1 +B_0- \beta_0 \Lambda Q ,\right. \notag\\
        & \quad \quad \quad \quad \left. -3a_1 \frac{\alpha+2}{\alpha+1} Q^2, \quad 3Q^2a_2, \quad  3Q^2 \left(a_1 \frac{(\alpha+1)(\alpha+2)}{2}y^2 + a_3\right)
        \right).
\end{align}

There exist unique $\beta_0\in \mathbb{R}$, $\beta(\Gamma)\in C^1(I)$ satisfying \eqref{ineq:beta_delta}, $B_0, B_1, B_2, B_3, B_4$ and $B_5 \in X^\infty(\R)$, with $B_0$ an even function, and $B_0, B_1+b_1 S_0, B_2, B_3, B_4, B_5 \perp Q, Q'$ such that, :
\begin{align}\label{defi:B}
    \partial_y\left( L \left( \overrightarrow{f}(\Gamma) \cdot \overrightarrow{B} -b(z)S_0 \right)+ \overrightarrow{f}(\Gamma) \cdot\overrightarrow{F}(1,\beta_0,B_0) \right)= b(z)\Lambda Q + \beta(\Gamma)Q,
\end{align}
with $\overrightarrow{B}:=(B_0,B_1,B_2,B_3,B_4,B_5)$.

Similarly, there exist unique $\delta_0 \in \mathbb{R}$, $\delta(\Gamma)\in \mathcal{C}^1(I)$, satisfying \eqref{ineq:beta_delta}, $D_0, D_1, D_2, D_3, D_4$ and $D_5\in X^{\infty}(\R)$ with $D_0$ an even function, $D_0, D_1+b_1(l-S_0),D_2,D_3, D_4,D_5 \perp Q, Q'$ such that:
\begin{align}\label{defi:D}
    \partial_y\left( L \left( \overrightarrow{f}(\Gamma) \cdot \overrightarrow{D} +b(z) (l-S_0)\right) + \overrightarrow{f}(\Gamma) \cdot \overrightarrow{F}(2,\delta_0,D_0) \right)= -b(z)\Lambda Q + \delta(\Gamma) Q',
\end{align}
with $\overrightarrow{D}:=(D_0, D_1,D_2,D_3,D_4,D_5)$.
\end{propo}

Notice that in the previous decomposition, the tail of the profile of the first solitary wave, given by $B_1+b_1S_0$, has an influence on the profiles around the second solitary wave, on $D_1$. However, this tail does not change the coefficient $-b(z)\Lambda Q$, which is of great importance in the system of ODEs ruling the equations of $\mu$ and of $z$. 

Remark that the quantities $\mathscr{F}$ involved in Proposition \ref{propo:est:profils} correspond to translated versions of $\left(\overrightarrow{f}(\Gamma)\cdot \overrightarrow{F}(i,\beta_0,B_0)\right)(y-z_1)$.

To prove Proposition \ref{propo:profils}, we need Lemma \ref{lemm:ant1} and \ref{lemm:ant2} to find the adequate profiles.

\begin{proof}[Proof]

We define from Lemma \ref{lemm:ant1} the unique function $B_0\in X^{\infty}(\R)$ and the unique coefficient $\beta_0\in \mathbb{R}$ satisfying:
\begin{align}\label{defi:B_0}
    \begin{cases}
    L B_0(y) = -3a_1Q^2(y) + \beta_0 Q(y), \\
    B_0\perp Q,\quad B_0\perp Q'.
    \end{cases}
\end{align}

Notice that since $L$ keeps stable the parity of the functions, $B_0$ is an even function.

For the second term, we use Lemma \ref{lemm:ant2} by defining the function $B_1$, and the coefficients $\beta_1$ and $b_1$ as the unique solution of the following problem:

\begin{align}\label{defi:B_1}
    \begin{cases}
    \partial_yL (B_1(y)-b_1S_0(y)) =  \partial_y \left( 3(\alpha+1)a_1yQ^2(y) \right) + \beta_1 Q'(y) + b_1 \Lambda Q, \\
    B_1-b_1S_0\perp Q,\quad B_1-b_1S_0\perp Q'.
    \end{cases}
\end{align}

Notice in particular that $b_1$ is defined by the formula \eqref{defi:tilde_a}:
\begin{align}\label{defi:b_1}
   b_1= -\frac{2(\alpha+1)^2a_1\|Q\|_{L^3}^3}{(\alpha-1)\|Q\|_{L^2}^2}<0,
\end{align}
since the sign of $a_1>0$ is given in Lemma \ref{thm:asympQ}. This justifies the choice of definition of $b(z):=\displaystyle\frac{b_1}{z^{2+\alpha}}$, as stated in \eqref{defi:b}.

The third, fourth, fifth and sixth terms are defined as for $B_0$ and $\beta_0$. With Lemma \ref{lemm:ant1}, we define $B_2, B_3, B_4, B_5$ in $X^{\infty}(\R)$, and the coefficients $\beta_2, \beta_3$, $\beta_4$ and $\beta_5$ as the solutions of the following problems:
\begin{align}\label{defi:B_2_and_B_3}
    \begin{cases}
    L B_2(y) = -6Q\Lambda Q B_0-6Q\Lambda Q a_1 -B_0+ \beta_0 \Lambda Q + \beta_2 Q(y) \\
    B_2\perp Q,\quad B_2\perp Q'
    \end{cases}
    , 
    \begin{cases}
    L B_3(y) = \displaystyle3a_1 \frac{\alpha+2}{\alpha+1} Q^2+ \beta_3 Q(y) \\
    B_3\perp Q,\quad B_3\perp Q'
    \end{cases}
\end{align}
and
\begin{align}\label{defi:B_4}
    \begin{cases}
    L B_4(y) = \displaystyle -3a_2Q^2(y)+ \beta_4 Q(y) \\
    B_4\perp Q,\quad B_4\perp Q'.
    \end{cases}
     , 
    \begin{cases}
    L B_5(y) = \displaystyle -3\left(a_1 \frac{(\alpha+1)(\alpha+2)}{2}y^2 + a_3\right)Q^2(y)+ \beta_5 Q(y) \\
    B_5\perp Q,\quad B_5\perp Q'.
    \end{cases}
\end{align}
Therefore, we set: 
\begin{align}\label{defi:beta}
    \beta(\Gamma)=(\beta_0,\beta_1,\beta_2,\beta_3,\beta_4,\beta_5)\cdot\overrightarrow{f}(\Gamma).
\end{align}
Now, we continue with the construction of $\overrightarrow{D}$. Since the first, forth, fifth and sixth coordinates in $\overrightarrow{F}(2,\delta_0,D_0)$ are respectively equal to the first, forth, fifth and sixth terms in $\overrightarrow{F}(1,\beta_0,D_0)$, the functions $D_0$, $D_3$, $D_4$, $D_5$ will solve respectively the same problem as $B_0$, $B_3$, $B_4$ and $B_5$. Then, we take: 
\begin{align}
D_0=B_0,\quad  D_3=B_3, \quad  D_4=B_4, \quad D_5=B_5,
\end{align}
and 
\begin{align}
\delta_0=\beta_0, \quad \delta_3=\beta_3, \quad  \delta_4=\beta_4, \quad  \delta_5=\beta_5. \label{delta0=beta0}
\end{align}
The situation is similar for $D_2=B_2$ and for $\beta_2=\delta_2$. To construct $D_1$, as for the function $B_1$, we use Lemma \ref{lemm:ant2}. Since $z^{2+\alpha}b(z)=b_1$, there exist a unique function $D_1\in X^{\infty}(\R)$ and coefficients $\delta_1,d_1\in\R$ such that:  \begin{align}\label{defi:D_1}
    \begin{cases}
    \partial_yL \left(D_1(y)+d_1(l-S_0(y))\right) = \partial_y \left(-(\alpha+1)a_1 y 3Q^2(y) \right) + \delta_1 Q'(y) + d_1 \Lambda Q(y), \\
    D_1+d_1(l-S_0)\perp Q,\quad D_1+d_1(l-S_0)\perp Q'.
    \end{cases}
\end{align}
Moreover, $Q^2$ is orthogonal to $Q'$. Therefore by the formula \eqref{defi:tilde_a}, we obtain that:
\begin{align}
    d_1=-b_1.
\end{align}
Thus, we conclude the proof of Proposition \ref{propo:profils} by defining: 
\begin{align}
    \delta(\Gamma) :=(\delta_0,\delta_1,\delta_2,\delta_3,\delta_4,\delta_5)\cdot\overrightarrow{f}(\Gamma).
\end{align}
\end{proof}

\begin{proof}[Proof of Proposition \ref{propo:est:profils}]
The two identities \eqref{iden:P_1} and \eqref{iden:P_2} are deduced from the one of $\overrightarrow{B}$ and $\overrightarrow{D}$ in \eqref{defi:B} and \eqref{defi:D}, as well as the orthogonality conditions.

We continue with the estimate \eqref{derivee:d_yP_i} and \eqref{derivee:d_tP_i}. First, we deal with the term $\partial_y B_0$. From \eqref{defi:B_0}, we deduce that:
\begin{align*}
    \partial_y B_0 =\left(|D|^{\alpha}+1 \right)^{-1}\partial_y \left( 3Q^2B_0 + 3a_1Q^2 + \beta_0Q\right).
\end{align*}
Since $B_0\in X^{2+\alpha}(\R)$, we have that $\partial_yB_0 \in L^{\infty}(\R)\cap \mathcal{C}(\mathbb{R})$. Then, by Lemma \ref{lemma:est:kernel}, we obtain that $\partial_y B_0\in X^{2+\alpha}(\R)$. By a similar argument on $B_1$, $B_2$, $B_3$, $B_4$ and $B_5$ with \eqref{hypothese:II}, we conclude that: 
$$
|\partial_yP_1(\Gamma,y) |=|\overrightarrow{f}(\Gamma)\cdot (\partial_y\overrightarrow{B})(y-z_1)|\leq \frac{C}{z^{1+\alpha}}\frac{1}{\langle y-z_1 \rangle^{1+\alpha}}.
$$
The same estimate holds for $P_2$.

Now, we estimate $\partial_t P_i$ for $i\in\{1,2\}$. Note that the profiles $P_1(\Gamma)$ and $P_2(\Gamma)$ are $C^{1}(I)$, since $\Gamma\in C^{1}(I)$. By direct computation, we obtain that:
\begin{align}\label{derive:profil1}
    \frac{d}{dt} P_1(\Gamma(t),y)= \left(\frac{d}{dt}\Gamma(t) \cdot \nabla_{\Gamma} \right) \overrightarrow{f}(\Gamma(t))\cdot \overrightarrow{B}(y-z_1(t))-\dot{z_1}(t) \partial_y P_1(\Gamma(t),y).
\end{align}
By Proposition \ref{propo:profils}, we have that $B_j\in X^{\infty}$ for $j\in\{0,\cdots,5\}$. Therefore, we deduce with \eqref{hypothese:II} and \eqref{hypothese:III}, that: 
\begin{align} 
    \left|\left( \dot{\Gamma} \cdot \nabla_{\Gamma} \right) \overrightarrow{f}(\Gamma)\cdot \overrightarrow{B}(y-z_1)\right| \leq C\left( \frac{|\dot{z}|}{z^{2+\alpha}} + \frac{\vert \dot{\mu}_1 \vert+ \vert \dot{\mu_2}\vert }{z^{1+\alpha}}\right) \frac{1}{\langle y-z_1 \rangle^{1+\alpha}} \leq \frac{C}{z^{\frac{5+3\alpha}{2}}}\frac{1}{\langle y-z_1 \rangle^{1+\alpha}}. \label{eq:d_t_P_i_bis}
\end{align}
We conclude that: 
\begin{align*}
    \left|\frac{d}{dt} P_1(\Gamma)+\dot{z}_1 \partial_y P_1(\Gamma)\right| \leq \frac{C}{z^{\frac{5+3\alpha}{2}}\langle y-z_1 \rangle^{1+\alpha}}. 
\end{align*}
The same arguments hold to estimate the profile $P_2$. This finishes the proof of Proposition \ref{propo:est:profils}.  
\end{proof}

\subsection{Proof of Proposition \ref{propo:estimates_flow} and Theorem \ref{theo:construction} }

\paragraph{}

Once the construction of the profiles is finished we continue with the estimates of the different terms involved in the error. 

\begin{proof}[Proof of Proposition \ref{propo:estimates_flow}]

To obtain \eqref{eq:partial_y_Ri}, we have the decomposition on $\partial_y R_1$:
\begin{align*}
    \| \partial_y R_1 \phi \|_{L^{\infty}(\{y\leq \frac{z_1}{2}\})} & \leq \left\| \frac{C}{\langle y -z_1 \rangle ^{1+\alpha}} \frac{1}{\langle y \rangle ^{\alpha}} \right\|_{L^\infty} \leq \frac{C}{z^{\alpha}}, \\
    \| \partial_y R_1 \phi \|_{L^{\infty}(\{y\leq \frac{z_1}{2}\})} & \leq \| \partial_y R_1 \|_{L^{\infty}(\{y\leq \frac{z_1}{2}\})} \leq \frac{C}{z^{1+\alpha}}.
\end{align*}
The same estimate holds for $R_2$. Applying the same argument for the $H^1$-norm, we deduce \eqref{eq:partial_y_Ri}. We can replace $\partial_yR_i$ by $\Lambda R_i$ in the former estimates and we get \eqref{eq:Ri_Phi}.

The estimate \eqref{partialyPibWchi} and \eqref{partialtPibWchi} are direct consequences of Proposition \ref{propo:est:profils} and the definition of $b$.

	By Proposition \ref{propo:profils} the profiles $P_i$ and $\partial_yP_i$ for $i=1,2$ belong to $L^{\infty}(\R)$. Moreover, by definition, $W$ and  $\partial_y W$ are also in $L^{\infty}(\R)$. Then we deduce \eqref{V}.
	
	By Proposition \ref{propo:profils} for the profiles, and since $b(z)=\displaystyle\frac{b_1}{z^{\alpha+2}}$, we deduce that :
	\begin{align*}
	\|P_1\|_{L^{\infty}} + \|P_2\|_{L^{\infty}} +\|bW \|_{L^{\infty}}\leq \frac{C}{z^{1+\alpha}}.	
	\end{align*}
Furthermore, using $\Omega_i=\displaystyle\{y\in\R : y\leq \frac{z_i}{2}\}$, we get for $i=1,2$ that:
\begin{align*}
    \|R_i^k \Phi^2\|_{L^{\infty}}\leq \|R_i^k \Phi^2\|_{L^{\infty}(\Omega_i)} + \|R_i^k \Phi^2\|_{L^{\infty}(\Omega_i^{c})}\leq \frac{C}{z^{1+\alpha}}.
	\end{align*}
Gathering these estimates, we conclude the first part of \eqref{Vphi}. Concerning the second term:
\begin{align*}
\MoveEqLeft
    \|(V^2-R_1^2)\partial_y R_1\|_{L^2}
    \\
    & \leq \|2R_1(R_2 -P_1+P_2+bW)\partial_y R_1 \|_{L^2} + \|R_2^2 \partial_y R_1\|_{L^2} + C \| -P_1+P_2+bW \|_{L^\infty}^2 \leq \frac{C}{z^{1+\alpha}}.
\end{align*}
By differentiating $V$ and using Proposition \ref{propo:est:profils}, therefore we obtain the estimate \eqref{partialtV}. 
	\end{proof}

\begin{proof}[Proof of Theorem \ref{theo:construction}]

We continue with the inequalities \eqref{est:S}, \eqref{est:dtS} and \eqref{est:T}.

We first begin with the estimate on the $L^2$-norm of the term $S=S_{V^3}+S_1+S_2+ \tilde{S}$ with $S_{V^3},S_1,S_2,\tilde{S}$ are respectively defined in  \eqref{defi:S_V_3}, \eqref{defi:S_1}, \eqref{defi:S_2} and \eqref{defi:S_tilde}.

We begin with $S_{V^3}$, by decomposing the different terms. We have, using the decomposition of Proposition \ref{propo:DL_Q_mu}:
\begin{align}
    \left\| \left( R_1^2-\widetilde{R}_1^2\right) P_2 \right\|_{L^2} \leq C \vert \mu_1 \vert \left\| \frac{1}{\langle y-z_1\rangle^{1+\alpha}} P_2 \right\|_{L^2}.
\end{align}
Let $\Omega := \left\{ y \leq \frac{z_1+z_2}{2}\right\}$. By \eqref{derivee:d_yP_i} and \eqref{hypothese:II}, we obtain that:
\begin{align}
\MoveEqLeft
    \left\| \left( R_1^2 - \tilde{R}_1^2 \right) P_2 \right\|_{L^2} \notag\\
    & \leq C\frac{\vert \mu_1 \vert}{z^{1+\alpha}} \left( \left\| \frac{1}{\langle y-z_1\rangle^{1+\alpha}} \frac{1}{\langle y-z_2 \rangle^{1+\alpha}} \right\|_{L^2(\Omega)} + \left\| \frac{1}{\langle y-z_1\rangle^{1+\alpha}} \frac{1}{\langle y-z_2 \rangle^{1+\alpha}} \right\|_{L^2(\Omega^C)} \right) \leq \frac{C}{z^{\frac{5+3\alpha}{2}}}. \label{eq:intersection}
\end{align}
By similar computations, we have that:
\begin{align}
    \left\| R_1^2 P_2\right\|_{L^2} \leq \left\| \tilde{R}_1^2 P_2 \right\|_{L^2} + \left\|\left( R_1^2-\tilde{R}_1^2\right) P_2 \right\|_{L^2} \leq \frac{C}{z^{\frac{5+3\alpha}{2}}}.
\end{align}
Similarly:
\begin{align*}
    \left\| 3R_1^2 (-bS_0(\cdot-z_2)) + 3 R_2^2 \left( -P_1 + b(S_0(y-z_1)-l)\right)\right\|_{L^2} \leq \frac{C}{z^{\frac{5+3\alpha}{2}}}.
\end{align*}
For the third and forth terms of $S_{V^3}$, by \eqref{partialyPibWchi}, we have:
\begin{align*}
    \left\| R_1 R_2 (-P_1+P_2+bW)\right\|_{L^2} \leq \| R_1 R_2 \|_{L^2} \| -P_1+P_2+bW \|_{L^{\infty}} \leq \frac{C}{z^{\frac{5+3\alpha}{2}}}
\end{align*}
and 
\begin{align*}
    \left\|(-R_1+R_2)\left(-P_1+P_2+bW \right)^2 \right\|_{L^2} \leq \| -R_1+R_2\|_{L^2} \| -P_1+P_2+bW \|_{L^\infty}^2 \leq \frac{C}{z^{\frac{5+3\alpha}{2}}}.
\end{align*}
Finally, we compute the $L^2$-norm of the bump function $W$:
\begin{align*}
    \| b W \|_{L^2} \leq \frac{C}{z^{2+\alpha}} \sqrt{z} = \frac{C}{z^{\frac{3}{2}+\alpha}},
\end{align*}
and therefore:
\begin{align*}
    \left\| (-P_1+P_2+bW)^3 \right\|_{L^2} \leq \left\| (-P_1+P_2+bW) \right\|_{L^\infty}^2 \left\| (-P_1+P_2+bW) \right\|_{L^2} \leq \frac{C}{z^{\frac{5+3\alpha}{2}}}.
\end{align*}
With the previous computations, we conclude that :
\begin{align*}
    \|S_{V^3}\|_{L^2} \leq \frac{C}{z^{\frac{5+3\alpha}{2}}}.
\end{align*}

We continue with $S_1$. Notice that by definition of $P_i$, another formulation of $S_1$ and $S_2$ is available:
\begin{align}
    S_1 & = \mathscr{S}(1,2) - \mathscr{F}(1,2,\beta_0,B_0) \label{S_1_without_tail} \\
    S_2 & = -\mathscr{S}(2,1) + \mathscr{F}(2,1,\delta_0,D_0). \label{S_2_without_tail}
\end{align}
We focus on $S_1$, the computations are similar for $S_2$. We separate each term of $\mathscr{S}(1,2)- \mathscr{F}(1,2,\beta_0,B_0)$. First we look at $\displaystyle  \left\| 3 \tilde{R}_1^2 (\tilde{R}_2 - Q_{app}(\cdot -z_2,z)) \right\|_{L^2}$. The approximation of $Q(\cdot+z)$ by $Q_{app}(\cdot,z)$ in \eqref{eq:estimate_Qapp} holds on a certain region, thus we begin with $\left\{y\in \mathbb{R}; \vert y - z_1 \vert \leq \frac{z}{2} \right\}$. In this region, we have:
\begin{align*}
    \left\| 3 \tilde{R}_1^2 (\tilde{R}_2 - Q_{app}(\cdot -z_1,z)) \right\|_{L^2(\vert y - z_1 \vert \leq \frac{z}{2})} & \leq C \left\| 3 \tilde{R}_1^2 \left(\frac{1}{z^{3\alpha+1}} + \frac{\langle y-z_1 \rangle}{z^{2\alpha+2}} + \frac{\langle y-z_1 \rangle^ 3}{z^{\alpha+4}} \right) \right\|_{L^2} \\
    & \leq \frac{C}{z^{2\alpha+2}} \leq \frac{C}{z^{\frac{5+3\alpha}{2}}}.
\end{align*}
In the other part, we get:
\begin{align*}
    \left\| 3 \tilde{R}_1^2 (\tilde{R}_2 - Q_{app}(\cdot -z_1,z)) \right\|_{L^2(\vert y - z_1 \vert \geq\frac{z}{2})} \leq \left\| 3 \tilde{R}_1^2 \tilde{R}_2  \right\|_{L^2(\vert y - z_1 \vert \geq\frac{z}{2})} + \left\| 3 \tilde{R}_1^2 Q_{app}(\cdot -z_2,z) \right\|_{L^2(\vert y - z_1 \vert \geq\frac{z}{2})}.
\end{align*}
The first term on the right hand side of the former estimate is bounded by:
\begin{align*}
    \left\| 3 \tilde{R}_1^2\tilde{R}_2  \right\|_{L^2(\vert y - z_1 \vert \geq\frac{z}{2})}\leq C \| \widetilde{R}_1^2\|_{L^{\infty}(\vert y - z_1 \vert \geq\frac{z}{2})}\|\widetilde{R}_2\|_{L^2}\leq \frac{C}{z^{\frac{5+3\alpha}{2}}}.
\end{align*}
We estimate the second term on the right hand side by:
\begin{align*}
\MoveEqLeft
    \left\|  \tilde{R}_1^2 Q_{app}(\cdot -z_1,z) \right\|_{L^2(\vert y - z_1 \vert \geq\frac{z}{2})} \\
    &\leq \frac{C}{z^{\frac{5+3\alpha}{2}}} 
    +\frac{C}{z^{2+\alpha}}\left\|  \tilde{R}_1^2(y-z_1) \right\|_{L^2(\vert y - z_1 \vert \geq\frac{z}{2})}
    +\frac{C}{z^{3+\alpha}}\left\|  \tilde{R}_1^2(y-z_1)^3 \right\|_{L^2(\vert y - z_1 \vert \geq\frac{z}{2})} \leq \frac{C}{z^{\frac{5+3\alpha}{2}}}.
\end{align*}
Thus we conclude:
\begin{align*}
    \displaystyle  \left\| 3 \tilde{R}_1^2 (\tilde{R}_2 - Q_{app}(\cdot -z_1,z)) \right\|_{L^2}\leq \frac{C}{z^{\frac{5+3\alpha}{2}}}.
\end{align*}
The estimates on the other terms of $S_1$ are obtained by similar computations:
\begin{align*}
\MoveEqLeft
    \left\| 6 \mu_1 \tilde{R_1}\Lambda \tilde{R}_1 \left(P_1 -\frac{B_0(\cdot -z_1)}{z^{1+\alpha}} \right) \right\|_{L^2} + \left\| 6 \mu_1 \tilde{R_1} \Lambda \tilde{R}_1 \left(\tilde{R}_2 -\frac{a_0}{z^{1+\alpha}} \right) \right\|_{L^2} \\
    & + \left\| 3\mu_2 \tilde{R}_1^2 \left( \Lambda \tilde{R}_2 + \frac{1}{z^{1+\alpha}} \frac{a_0(\alpha+2)}{2(\alpha+1)}\right) \right\|_{L^2} +\left\| \mu_1 \left( P_1 -\frac{B_0(\cdot-z_1) }{z^{1+\alpha}} \right)\right\|_{L^2} \leq \frac{C}{z^{\frac{5+3\alpha}{2}}}
\end{align*}
then we conclude:
\begin{align*}
  \|S_1\|_{L^2}\leq \frac{C}{z^{\frac{5+3\alpha}{2}}}.  
\end{align*}
To finish the proof on S, we have to estimate $\widetilde{S}$. We focus on the first part of $\widetilde{S}$, which contains $\widetilde{\mathscr{S}}(1,2)$:
$$3\left(R_1^2 - \widetilde{R}_1^2\right)\left(-P_1+b S_0(y-z_1) \right) +6\mu_1 \tilde{R}_1 \Lambda \tilde{R}_1 P_1 + 3 R_1^2 R_2 - 3 \tilde{R}_1^2 \tilde{R}_2 - 6\mu_1 \tilde{R}_1 \Lambda \tilde{R}_1 \tilde{R}_2 - 3 \mu_2 \tilde{R}_1^2 \Lambda \tilde{R}_2 $$
 since the computations are similar for the other part. By using Proposition \ref{propo:DL_Q_mu}, \eqref{hypothese:II} and \eqref{defi:b} we deduce:
\begin{align*}
    \|3(R_1^2-\widetilde{R}_1^2)bS_0(y-z_1)\|_{L^2}+\|3(R_1^2-\widetilde{R}_1^2-2\mu_1\widetilde{R}_1\Lambda \widetilde{R}_1)P_1\|_{L^2} \leq \frac{C}{z^{\frac{5+3\alpha}{2}}}.
\end{align*}
To estimate the next terms in $\widetilde{S}$, we remark:
\begin{align*}
    R_1^2R_2-\widetilde{R}_1^2\widetilde{R}_2= R_1^2(R_2-\widetilde{R}_2) + \widetilde{R}_2(R_1^2-\widetilde{R}_1^2).
\end{align*}
From Proposition \ref{propo:DL_Q_mu} and \eqref{hypothese:II}, we obtain that:
\begin{align*}
    \|3R_1^2(R_2-\widetilde{R}_2) - 3\mu_2\Lambda \widetilde{R}_2\widetilde{R}_1^2\|_{L^2}\leq & \|3R_1^2(R_2-\widetilde{R}_2) - 3\mu_2\Lambda \widetilde{R}_2 R_1^2\|_{L^2} + \|3\mu_2\Lambda \widetilde{R}_2\left(R_1^2-\widetilde{R}_1^2\right) \|_{L^2})\\
    \leq &C\left(\frac{\mu_2^2}{z^{1+\alpha}} + \frac{|\mu_1||\mu_2|}{z^{1+\alpha}}  \right)\leq \frac{C}{z^{\frac{5+3\alpha}{2}}}.
\end{align*}
Arguing similarly, we obtain:
\begin{align*}
    \|3\widetilde{R}_2(R_1^2-\widetilde{R}_1^2) - 6\mu_1\widetilde{R}_1\Lambda\widetilde{R}_1\widetilde{R}_2\|_{L^{2}}\leq \frac{C}{z^{\frac{5+3\alpha}{2}}}.
\end{align*}
This concludes the estimate on $\widetilde{S}$.

We continue with the estimate on $T$. We decompose each term of its definition in \eqref{defi:T}. First, we have with \eqref{eq:lambda_Q_DL} and \eqref{hypothese:II} :
\begin{align*}
    \left\| b(z) \left( \Lambda R_1 -\Lambda \tilde{R}_1 \right)\right\|_{L^2} \leq C \frac{\vert \mu_1 \vert}{z^{2+\alpha}} \left\|\frac{1}{\langle x-z_1 \rangle^{\alpha+1}}\right\|_{L^2} \leq \frac{C}{z^{\frac{5+3\alpha}{2}}}.
\end{align*}
Second, we use the inequality \eqref{ineq:beta_delta} as used in Proposition \ref{propo:est:profils}, and from the asymptotic development of $\partial_y Q$, by \eqref{eq:lambda_Q_H_1} and \eqref{hypothese:III}:
\begin{align*}
\MoveEqLeft
    \left\| \beta(\Gamma)\left( - \partial_y R_1 + \partial_y \tilde{R}_1 \right) + \frac{\beta_0}{z^{1+\alpha}}\partial_y (\mu_i \Lambda \tilde{R}_i) \right\|_{L^2} \\
       & \leq \left\| \left( \beta(\Gamma) - \frac{\beta_0}{z^{1+\alpha}} \right) \partial_y (-R_1 + \tilde{R}_1) \right\|_{L^2} + \left\| \frac{\beta_0}{z^{1+\alpha}} \partial_y \left( -R_1 + \tilde{R}_1 +\mu_1 \Lambda \tilde{R_1} \right)  \right\|_{L^2} \\
       & \leq C \frac{\vert \mu_1 \vert}{z^{2+\alpha}} + C\frac{\mu_1^2}{z^{1+\alpha}} \leq \frac{C}{z^{\frac{5+3\alpha}{2}}}.
\end{align*}
Then, we consider the case of the time derivative on $-P_1$. We have, by \eqref{derivee:d_tP_i_asymptotic}:
\begin{align*}
    \left\| \frac{d}{dt}(-P_1) - \mu_1 \partial_y P_i  \right\|_{L^2} \leq \frac{C}{z^{\frac{5+3\alpha}{2}}} + \frac{C}{z^{1+\alpha}} \left\vert \dot{z}_1 -\mu_1\right\vert 
\end{align*}

We continue with the the term $\frac{d}{dt}W$:
\begin{align}\label{defi:d_t_W}
\frac{d}{dt} W(\Gamma(t)) = \left( \vert D \vert^\alpha+1 \right)^{-1} \left( \dot{z}_1(t) \Lambda \tilde{R}_1 - \dot{z}_2 \Lambda \tilde{R}_2\right),
\end{align}
which with \eqref{hypothese:II} and \eqref{defi:d_t_W} give:
\begin{align*}
    & \left\| \frac{d}{dt} \left(b(z(t))W(\Gamma(t)) \right) \right\|_{L^2} \leq C \frac{\vert \dot{z}_1 \vert + \vert \dot{z}_2 \vert}{z^{2+\alpha}} + C\frac{\vert \dot{z}\vert}{z^{3+\alpha}} \|W \|_{L^2} \\
        & \quad \quad \quad \quad \quad \quad \quad \quad \quad \quad \leq \frac{C}{z^{\frac{5+3\alpha}{2}}} +C\frac{\sqrt{z}}{z^{\frac{7+3\alpha}{2}}} \leq \frac{C}{z^{\frac{5+3\alpha}{2}}}.
\end{align*}
Those previous estimates conclude the bound \eqref{est:T} on $T$.

Since all the estimates have been established in $L^2$, we need to continue with the first derivative to establish the bound in $H^1$. We can notice that all the estimates are based on two main arguments:
\begin{itemize}
    \item An argument of localisation : if two functions are located at a distance $z$ large, and if the two functions have an explicit decay at infinity, then the product of the two functions can be quantified in terms of $z$. The spatial derivative either leaves unchanged the decay property in terms of $z$ of this product or improves it.
    \item An argument of smallness of the objects: the objects already have a quantified bound in terms of $z$, see for example the $L^\infty$-norm of $P_i$ in \eqref{derivee:d_yP_i}.
\end{itemize}

Therefore the computations made on the $L^2$-norm are similar to those on the $H^1$-norm. 

Concerning the time derivative of $S$ in \eqref{est:dtS}, let us deal with a generic example of a function $\frac{1}{z(t)^{1+\alpha}} g_{1+\mu(t)}(y-z(t))$, since all the involved functions, except $W$, are on this form. Either the time derivative applies to $\frac{1}{z(t)^{1+\alpha}}$, or to the scaling parameter $1+\mu(t)$ of the function $g$ or to the translation parameter $-z(t)$. However, we get in each case either $\dot{\mu}(t)$ or $\dot{z}(t)$, which by \eqref{hypothese:II} and \eqref{hypothese:III} are bounded by $z^{-\frac{1+\alpha}{2}}$. Notice also that the time derivative of the considered functions leaves unchanged or improves the space decay at infinity, and from the remark on the space derivative above, the bound in $z$ still holds.
The time derivative of $W$ has been developed in \eqref{defi:d_t_W}, and $\partial_t W$ fits in the previous discussion. As a result, the estimate on $\| \partial_t S \|_{L^2}$ is reduced to the product of two terms: one whose bound is the one of $\| S\|_{L^2}$, and one bounded by $z^{-\frac{1+\alpha}{2}}$.
\end{proof}

\section{Modulation}\label{sec:modulation}

The previous section was dedicated to the expected approximate solution. Here, we prove that if a solution is close to the approximation $V$, for two solitary waves far enough one to each other, then the solution stays close to this approximation on a certain time interval. Furthermore, we can impose some orthogonality conditions to the error between the solution and the approximation.

Let us define some conditions \eqref{condition:Cond_Z} on a vector $\Gamma=(z_1,z_2,\mu_1,\mu_2)\in \mathbb{R}^4$ dependent on a parameter $Z$:
\begin{align}\label{condition:Cond_Z}\tag{Cond$_Z$}
    z_1 > \frac{Z}{4}, \quad z_2 < - \frac{Z}{4}, \quad 0 <- \mu_1 <\frac{1}{Z}, \quad \text{and} \quad 0<\mu_2 < \frac{1}{Z},
\end{align}
and the tube:
\begin{equation}
    \mathcal{U}(Z,\nu) := \left\{ u \in H^{\frac{\alpha}{2}}(\R); \inf_{\Gamma \text{ satisfying \eqref{condition:Cond_Z}}} \| u - V(\Gamma) \|_{H^{\frac{\alpha}{2}}} \leq \nu\right\}.
\end{equation}
We also shorten the notations by:

\begin{align}\label{defi:R_i_modulation}
R_i(y)= R_i(\Gamma,y) := Q_{1+\mu_i}(y-z_i).
\end{align}

This proposition is time-dependent, and can be found, for example, in \cite{MMT02,eychenne2021asymptotic}.

\begin{propo}\label{prop:modulation}
    There exist $Z^*>0$, $\nu^*>0$ and a constant $K^*>0$ such that the following holds. Let $v$ be a solution of \eqref{mBO}  in $\mathcal{C}(\mathbb{R}, H^{\frac{\alpha}{2}})$. Let us define a time interval $I$. If for $Z>2Z^*$ and $\nu \in (0,\frac{\nu^*}{2})$, we have :
    \begin{align*}
        \sup_{t\in I} \left( \inf_{\Gamma \text{ satisfying \eqref{condition:Cond_Z}}} \left\| v(t,\cdot) - V(\Gamma,\cdot) \right\|_{H^\frac{\alpha}{2}} \right) < \nu,
    \end{align*}
    then there exists a unique $\mathcal{C}^1$-function $\Gamma: I \rightarrow \mathbb{R}^4$ such that:
    \begin{align*}
        \eps(t,\cdot):= v(t,\cdot) - V(\Gamma(t),\cdot)
    \end{align*}
    satisfies for any $i \in \{1,2\}$ and for any $t\in I$:
    \begin{align}\label{eq:ortho_modulation}
        \eps(t,\cdot) \perp R_i(t,\cdot) \quad \text{ and } \quad \eps(t,,\cdot) \perp \partial_y R_i(t,\cdot).
    \end{align}
Moreover, for any $t\in I$:
    \begin{align}
        & \| \eps(t,\cdot)  \|_{H^\frac{\alpha}{2}} + \vert \mu_1(t) \vert + \vert \mu_2(t) \vert \leq K^* \nu ,\label{eps:petit} \\
        &|\dot{z}_1(t)| + |\dot{z}_2(t)| +|\dot{\mu}_1(t)| +|\dot{\mu}_2(t)| \leq K^{*}, \label{derive:modulation:borne}\\
        & z_1(t) > \frac{Z}{8} , \quad z_2(t) \leq - \frac{Z}{8} \label{z_i:grand} .
    \end{align}
\end{propo}

\begin{proof}
We give here some insights of the proof. The proof is composed of two steps. The first part involves a qualitative version of the implicit function theorem, see section 2.2 in \cite{chow1982bifurcation}, to obtain the existence of the continuous function $\Gamma$. To this end, we study the functional  
\begin{equation}
    \begin{aligned}
        g : & \mathcal{U}(Z,\nu) \times \mathbb{R}_+^* \times \mathbb{R}_-^*\times \mathbb{R} \times \mathbb{R} & \longrightarrow & \hspace{3cm} \mathbb{R}^4 \\
            & \hspace{1.5cm} (w,z_1,z_2,\mu_1,\mu_2) & \longmapsto & \begin{pmatrix}
                \int (w-V(\Gamma)) R_1, \quad \int (w-V(\Gamma)) \partial_y R_1  \\ \int (w-V(\Gamma)) R_2 
                 ,\quad  \int (w-V(\Gamma)) \partial_y R_2
            \end{pmatrix},
    \end{aligned}
\end{equation}
at the point $(V(\tilde{\Gamma}),\tilde{\Gamma})$ with $V$ defined in \eqref{defi:V} and $\tilde{\Gamma}$ satisfying $\eqref{condition:Cond_Z}$. Note that the estimates obtained on $g$ and $d_{\Gamma}g$ used to verify the implicit function theorem, are uniform in $\tilde{\Gamma}$ satisfying \eqref{condition:Cond_Z}, for $Z>2Z^*$ with $Z^*$ large enough, and $\nu< \frac{\nu^*}{2}$ with $\nu^*$ small enough. In other words, for all $\tilde{\Gamma}$, the function $\Gamma$ associated with $\tilde{\Gamma}$ given by the implicit function theorem is  defined on a ball $B(V(\tilde{\Gamma}),\nu)$, with $\nu$ independent of the point $V(\tilde{\Gamma})$. Since $\nu$ is chosen independently of $\tilde{\Gamma}$ satisfying \eqref{condition:Cond_Z}, we can extend by uniqueness the parameters to the whole tube $\mathcal{U}(Z,\nu)$. Therefore, we get $\Gamma\in C^{1}(\mathcal{U}(Z,\nu))$.    

However, the solution $u$ of \eqref{mBO} is only continuous, then we obtain that the function $\Gamma(t):=\Gamma(v(t,\cdot))$ is only continuous. To get more regularity, we use the Cauchy-Lipischtz theorem. By differentiating the orthogonality condition, we have that the parameters verify an ODE system. By using the Cauchy-Lipischtz theorem, we obtain the regularity of the parameters even though $u$ is only continuous.

\end{proof}

\begin{rema}\label{global_well_modul}
The parameters $z_1$, $z_2$, $\mu_1$, $\mu_2$ defined in Proposition \ref{prop:modulation}, verify an ODE system which is globally Lipschitz. In other words, the function $\Gamma$ is well-defined and $C^1(\R)$. However, the conclusion of the Proposition \ref{prop:modulation} are only verified for $t\in I$.   
\end{rema}

\section{Proof of the Theorem \ref{main_theo}}\label{sec:proof_main_theorem}

\subsection{Bootstrap setting}\label{subsec:bootstrap_setting}
Let $(S_n)_{n=0}^{+\infty}$ be a increasing sequence of times going to infinity, with $S_n>T_0$, for $T_0>1$ large enough to be  chosen later. Recall that $V$ is defined in \eqref{defi:V}. For all $n\in\N$, we define $u_n$ as being the solution of \eqref{mBO} verifying 
\begin{align}\label{CIu_n}
    v_n(S_n,\cdot)=V(\Gamma^{in}_n,\cdot),
\end{align}
with 
\begin{align}
    \Gamma_n^{in} & :=(z_{1,n}^{in},z_{2,n}^{in},\mu_{1,n}^{in},\mu_{2,n}^{in}), \\
    z_{1,n}^{in}& =-z_{2,n}^{in} :=\frac{z_n^{in}}{2},\quad \mu^{in}_{1,n}=-\mu_{2,n}^{in}:=\frac{\mu_n^{in}}{2}, \quad 
    \mu_n^{in}:=\sqrt{\frac{-4b_1}{\alpha+1}}\left(z^{in}_n\right)^{-\frac{\alpha+1}{2}},\label{eq:mu_n_in}\\
    \left(z^{in}_n\right)^{\frac{\alpha+3}{2}} & \in [a^{\frac{\alpha+3}{2}}S_n - S_n^{\frac{1}{2}+r}, a^{\frac{\alpha+3}{2}}S_n + S_n^{\frac{1}{2}+r}], \label{CI:z_n_in}
\end{align}
with $b_1$ defined in \eqref{defi:b_1}, $a=\left( \frac{\alpha+3}{2}\sqrt{\frac{-4b_1}{\alpha+1}}\right)^{\frac{2}{\alpha+3}}$ and $r=\frac{\alpha-1}{4(\alpha+3)}$. The constant $z^{in}_n$ will be fixed later.

By choosing $T_0$ large enough and $C_0= 2 \sqrt{\frac{-4b_1}{\alpha+1}}$, we can suppose that \eqref{hypothese:I}-\eqref{hypothese:III} and  \eqref{condition:Cond_Z} are satisfied by $\Gamma_n^{in}$ for any $n\in \mathbb{N}$. By \eqref{CIu_n}, $v_n(S_n)\in \mathcal{U}(Z,\nu)$ and $V(\Gamma_n^{in})$ satisfies the assumption of theorem \ref{theo:construction}. By continuity of $v_n$ (see Corollary \ref{solution:global_well_posed}), on an open time interval $I_n\ni S_n$, $\left\{v_n(t); t \in I_n \right\}$ is in $\mathcal{U}(Z,\nu)$. By applying Proposition \ref{prop:modulation}, we define a unique function $\Gamma_n=(z_{1,n},z_{2,n},\mu_{1,n},\mu_{2,n},)$ on $I_n$ such that the conditions \eqref{eq:ortho_modulation}, \eqref{eps:petit} and \eqref{z_i:grand} are satisfied and $\Gamma_n(S_n)=\Gamma_n^{in}$ by construction. $\Gamma_n$ also satisfies \eqref{hypothese:I}-\eqref{hypothese:III}, which justifies the setting of Theorem \ref{theo:construction}.

By sake of clarity, we drop the index $n$, and denote $v$, $\Gamma$, $z_{1}$, $z_{2}$, $\mu_{1}$, $\mu_{2}$ instead of $v_n$, $\Gamma_n$, $z_{1,n}$, $z_{2,n}$, $\mu_{1,n}$, $\mu_{2,n}$ for the subsections \ref{sec:system_ODE} and \ref{monotonicity}. 

As in Section \ref{sec:approximation}, we denote:
\begin{align}\label{defi:mu_bar}
    z:= z_1-z_2, \quad \mu:= \mu_1 -\mu_2, \quad \bar{z}:=z_1+z_2, \quad \bar{\mu}:=\mu_1+\mu_2 \quad \text{and} \quad \eps := v-V(\Gamma).
\end{align}

We introduce the bootstrap estimates 
\begin{numcases}{}
\|\eps(t)\|^2_{H^{\frac{\alpha}{2}}} \leq t^{-\frac{3\alpha+5}{\alpha+3}},\label{boot:eps}\\
\lvert z^{\frac{\alpha+3}{2}}(t) - a^{\frac{\alpha+3}{2}}t \rvert \leq t^{\frac{1}{2}+r}, \label{boot:z}\\
\left\lvert \mu(t) - \sqrt{\frac{-4b_1}{\alpha+1}} \frac{t^{-\frac{\alpha+1}{\alpha+3}}}{a^{\frac{\alpha+1}{2}}} \right\rvert \leq C^{*}t^{-\frac{5\alpha+11}{4(\alpha+3)}},\label{boot:mu}\\
\lvert \bar{z}(t) \rvert \leq C^{*}t^{-\frac{\alpha-1}{2(\alpha+3)}},\label{boot:z_bar}\\
\lvert \bar{\mu}(t) \rvert \leq C^{*}t^{-2\frac{\alpha+1}{\alpha+3}},\label{boot:mu_bar}
\end{numcases}

with $C^{*}>1$ to be fixed later. Note that the condition \eqref{boot:z} implies 
\begin{align}\label{boot:z2}
    \lvert z(t) -a t^{\frac{2}{\alpha+3}} \rvert \leq C t^{-r}.
\end{align}
We define 
\begin{align}
    t^{*}(z^{in}_n)=\inf\{t\in [T_0,S_n]: \forall \tilde{t}\in[t,S_n], \eqref{boot:eps}-\eqref{boot:mu_bar} \text{ is true }\}.
\end{align}

We want to prove that for an adequate choice of $z_n^{in}$ in \eqref{CI:z_n_in},  $t^*(z_n^{in})=T_0$. 

By the previous choice of $C_0$, the assumptions \eqref{hypothese:I}-\eqref{hypothese:III} on the approximation and the condition \eqref{condition:Cond_Z} on the modulation are satisfied on $(t^*(z_n^{in}),S_n]$, increasing $T_0$ if necessary.

\vspace{0.5cm}

The section \ref{sec:system_ODE} provides the tools to get a bound of $z_1$, $z_2$, $\mu_1$ and $\mu_2$, and the section \ref{monotonicity} the bound on $\|\eps \|_{H^{\frac{\alpha}{2}}}$. Next, in the section \ref{sec:topological_argument}, we prove that we can choose $z_n^{in}$ to close the bootstrap. We finish the proof of Theorem \ref{main_theo} in the section \ref{sec:conclusion}.

\begin{rema}
Notice that different parameters are involved along this section. We clarify the order in which they are fixed. First, we fix the parameter $A$, introduced in subsection \ref{monotonicity}; then the parameter $C^*$ involved in the bootstrap dependently of $A$, and finally, the initial time $T_0$ dependently of $A$ and $C^*$. 
\end{rema}

\subsection{System of ODE}\label{sec:system_ODE}

We now continue with the system of ODEs ruling the parameters $z_1$, $z_2$, $\mu_1$ and $\mu_2$. To do so, we compute the time derivative of the orthogonality conditions.

\begin{propo}\label{propo:ODEs}
The functions $z_1$, $z_2$, $\mu_1$ and $\mu_2$ satisfy that for all $i \in \{1,2\}$ :
\begin{align}
      \sum_{i=1}^2\vert \dot{\mu}_i(t) + (-1)^i b(z(t)) \vert \leq& C\left( \frac{1}{z^{\frac{3\alpha+5}{2}}(t)} + \frac{1}{z^{\alpha+5}(t)}\|\eps(t)\|_{H^{\frac{\alpha}{2}}} + \|\eps(t)\|_{H^{\frac{\alpha}{2}}}^2 \right), \label{eq:mupt1_2}
\end{align}
      and 
\begin{align} \label{eq:zpt1_2}
     \vert \dot{z}_1(t) - \mu_1(t) +\beta(\Gamma(t)) \vert +\vert \dot{z}_2(t) - \mu_2(t) +\delta(\Gamma(t)) \vert
         \leq&  C\left(\frac{1}{z^{\frac{3\alpha+5}{2}}(t)}+ \|\eps(t)\|_{H^{\frac{\alpha}{2}}} \right).
\end{align}
\end{propo}

\begin{proof}
We begin with the first orthogonality condition $\int \eps R_1$. Since $\eps=v-V$ and $v$ solves \eqref{mBO}, we deduce that: 
\begin{align*}
    \partial_t\eps +\partial_y\left(-|D|^{\alpha}\eps - \eps + \left(\eps+V\right)^3 -V^3 \right)=-\cE_V.
\end{align*}
By differentiating in time the equality $0=\int \eps R_1$ and using the fact $\int \eps \partial_yR_1=0$, we obtain that:
\begin{align*}
    0=\frac{d}{dt}\int \eps R_1 =& \int \left(-|D|^{\alpha}\eps - \eps + 3R_1^2\eps \right) \partial_yR_1 + \int \left( (V+\eps)^3 - V^3 - 3R_1^2\eps\right)\partial_yR_1 \\
    &-\int \overrightarrow{m}\cdot\overrightarrow{MV}R_1 -\int\partial_yS R_1 - \int T R_1 + \dot{\mu}_1\int\eps \Lambda R_1.
\end{align*}
By using the equation of $R_1$ and the condition $\eps\perp \partial_yR_1$, we deduce that:
\begin{align*}
    \int \left(-|D|^{\alpha}\eps - \eps + 3R_1^2\eps \right)\partial_yR_1=
    \int \left(-|D|^{\alpha}\eps - (1+\mu_1)\eps + 3R_1^2\eps \right)\partial_yR_1=0.
\end{align*}
Now, we continue with $\int \left( (V+\eps)^3 - V^3 - 3R_1^2\eps\right)\partial_yR_1$. First, note that: 
\begin{align*}
    (V+\eps)^3 - V^3 - 3R_1^2\eps= 3V\eps^2+\eps^3 +3\eps\left( -2R_1\left(R_2 - P_1+P_2+bW \right) + \left( R_2 -P_1 + P_2 + bW \right)^2 \right).
\end{align*}
We recall $\|V\|_{L^{\infty}}+\| \partial_yR_1 \|_{L^\infty}\leq C$. Therefore, using the Sobolev embedding $H^{\frac{1}{6}}(\R)\xhookrightarrow{}L^3(\R)$, we have that: 
\begin{align*}
    \bigg|\int \left(3\eps^2 V + \eps^3\right)\partial_yR_1 \bigg| \leq C\left(\|\eps\|_{L^2}^2+\|\eps\|^3_{H^{\frac{\alpha}{2}}} \right) .
\end{align*}
Furthermore, $| R_2\partial_yR_1|\leq \frac{C}{z^{1+\alpha}}$ and $|P_1|+|P_2|+|bW|\leq \frac{C}{z^{1+\alpha}}$, we conclude that: 
\begin{align*}
    \bigg|\int \left( (V+\eps)^3 - V^3 - 3R_1^2\eps\right)\partial_yR_1\bigg| \leq C\left(\frac{\|\eps\|_{L^2}}{z^{\alpha+1}} + \|\eps\|_{L^2}^2 + \|\eps\|_{H^{\frac{\alpha}{2}}}^3 \right).
\end{align*}
Let us estimate $\int \overrightarrow{m}\cdot\overrightarrow{MV}R_1$. Using the set $\{y\in\R:  y \leq \frac{z_1+z_2}{2}\}$, we get that: 
\begin{align*}
    \bigg|\int\Lambda R_2 R_1\bigg| + \bigg|\int\partial_y R_2 R_1\bigg|\leq \frac{C}{z^{\alpha+1}}.
\end{align*}
Moreover, with $R_1\perp \partial_yR_1$, we obtain that: 
\begin{align*}
   \bigg| \int \overrightarrow{m}\cdot\overrightarrow{MV}R_1 -\left(-\dot{\mu_1}+b(z) \right)\int\Lambda R_1 R_1 \bigg| \leq \frac{C}{z^{\alpha+1}}\left(|\dot{\mu}_2+b(z)|+|\dot{z}_2-\mu_2+\delta(\Gamma)| \right).
\end{align*}
Finally, using Cauchy-Schwarz inequality, \eqref{est:S} and \eqref{est:T} we get that: 
\begin{align*}
    \bigg|\int \partial_yS R_1 \bigg|+ \bigg|\int TR_1 \bigg| +|\dot{\mu}_1|\bigg|\int \eps \Lambda R_1 \bigg|\leq C \left( \frac{1}{z^{\frac{3\alpha+5}{2}}} + |\dot{\mu}_1|\|\eps\|_{L^2} \right).
\end{align*}
Gathering these estimates, and thanks to the facts  $\|\eps\|_{H^{\frac{\alpha}{2}}}\leq C\kappa$ and $|\dot{\mu}_i|+|\dot{z}_i|\leq C$ from \eqref{eps:petit} and \eqref{derive:modulation:borne}, we obtain that: 
\begin{align}\label{eq:mupt1}
\frac{\alpha-1}{2(\alpha+1)}\|Q\|_{L^2}^2|\dot{\mu_1}-b(z)| \leq& \frac{C}{z^{\alpha+1}}\left(|\dot{\mu}_2+b(z)| + |\dot{z}_2-\mu_2+\delta(\Gamma)|+\|\eps\|_{L^2} \right) + \frac{C}{z^{\frac{3\alpha+5}{2}}}\\
&+ C\left(|\dot{\mu}_1|\|\eps\|_{H^{\frac{\alpha}{2}}} + \|\eps\|_{H^{\frac{\alpha}{2}}}^2  \right)  \notag.  
\end{align}   
By similar computations, we also deduce that: 
\begin{align}\label{eq:mupt2}
\frac{\alpha-1}{2(\alpha+1)}\|Q\|_{L^2}^2|\dot{\mu_2}+b(z)| \leq& \frac{C}{z^{\alpha+1}}\left(|\dot{\mu}_1-b(z)| + |\dot{z}_1-\mu_1+\beta(\Gamma)|+\|\eps\|_{L^2} \right) + \frac{C}{z^{\frac{3\alpha+5}{2}}}\\
&+ C\left(|\dot{\mu}_1|\|\eps\|_{H^{\frac{\alpha}{2}}} + \|\eps\|_{H^{\frac{\alpha}{2}}}^2 \right)  \notag. 
\end{align}
Therefore, by adding \eqref{eq:mupt1} and \eqref{eq:mupt2} we obtain:
\begin{align}
      \sum_{i=1}^2\vert \dot{\mu}_i + (-1)^i b(z) \vert \leq& \frac{C}{z^{\alpha+1}}\left( \lvert\dot{z}_1-\mu_1+\beta(\Gamma)\rvert + \lvert\dot{z}_2-\mu_2+\delta(\Gamma)\rvert \right) \notag \\ +& C\left( \frac{1}{z^{\frac{3\alpha+5}{2}}} + \left(|\dot{\mu}_1|+|\dot{\mu}_2|+\frac{1}{z^{1+\alpha}}\right)\|\eps\|_{H^{\frac{\alpha}{2}}} + \|\eps\|_{H^{\frac{\alpha}{2}}}^2 \right) . \label{eq:mupt1_2_bis}
\end{align}

 Let us continue with the second orthogonality condition:
\begin{align*}
    0 = \frac{d}{dt} \int \eps \partial_y R_1 =& \int \left(-|D|^{\alpha}\eps - \eps + (V+\eps)^3 - V^3 \right) \partial_y^2R_1  -\int \overrightarrow{m}\cdot\overrightarrow{MV}\partial_yR_1 \\ & +\int S \partial_y^2R_1 - \int T \partial_yR_1 + \dot{\mu}_1\int\eps \partial_y\Lambda R_1 - \dot{z}_1\int\eps \partial_y^2R_1.
\end{align*}
Since $\lvert V \rvert+\left\lvert(\partial_y^2 +|D|^{\alpha}\partial_y^2)R_1\right\rvert\leq C$ and using the Sobolev embedding , $H^{\frac{1}{6}}(\R)\xhookrightarrow{}L^3(\R)$,  we deduce that: 
\begin{align*}
    \left\lvert \int \left(-|D|^{\alpha}\eps - \eps + (V+\eps)^3 - V^3 \right) \partial_y^2R_1 \right\rvert \leq C\left(\|\eps\|_{H^{\frac{\alpha}{2}}} + \|\eps\|_{H^{\frac{\alpha}{2}}}^2 + \|\eps\|_{H^{\frac{\alpha}{2}}}^3 \right).
\end{align*}
By developing $\overrightarrow{m}\cdot\overrightarrow{MV}$ and using the facts $\left\lvert\int \partial_yR_1\left(\partial_yR_2 +\Lambda R_2 \right) \right\rvert \leq \frac{C}{z^{\alpha+1}}$ and $\int\partial_yR_1\Lambda R_1=0$ since $\partial_yR_1$ is odd, we get that:
\begin{align*}
    \left\lvert \int \overrightarrow{m}\cdot\overrightarrow{MV} \partial_yR_1 - \left(\dot{z}_1- \mu_1 +\beta(\Gamma)\right)\int \left(\partial_yR_1 \right)^2  \right\rvert\leq \frac{C}{z^{\alpha+1}}\left(\left\lvert \dot{\mu}_2+b(z) \right\rvert +\left\lvert \dot{z}_2-\mu_2+\delta(\Gamma) \right\rvert\right).
\end{align*}
We estimate the last terms by applying Cauchy-Schwarz inequality, \eqref{est:S} and \eqref{est:T}. We have that: 
\begin{align*}
 \left\lvert\int S \partial_y^2R_1\right\rvert + \left\lvert\int T \partial_yR_1 \right\rvert + \left\lvert \dot{\mu}_1\int\eps \partial_y \Lambda R_1 \right\rvert + \left\lvert \dot{z}_1\int\eps \partial_y^2R_1  \right\rvert \leq C\left(\frac{1}{z^{\frac{3\alpha+5}{2}}} + \left(\left\lvert \dot{\mu}_1 \right\rvert +\left\lvert \dot{z}_1 \right\rvert\right) \|\eps\|_{L^2} \right). 
\end{align*}
Gathering these estimates and using $\|\eps\|_{H^{\frac{\alpha}{2}}}\leq C\kappa$ and the fact $|\dot{\mu}_i|+|\dot{z}_i|\leq C$ \eqref{derive:modulation:borne}, we conclude that: 
\begin{align}\label{eq:zpt1}
    \left\lvert \dot{z}_1 - \mu_1 +\beta(\Gamma) \right\rvert \int\left(\partial_yQ \right)^2 \leq C \left(\frac{1}{z^{\alpha+1}} \left\lvert \dot{\mu}_2+b(z) \right\rvert + \left\lvert \dot{z}_2 -\mu_2 +\delta(\Gamma) \right\rvert + \frac{1}{z^{\frac{3\alpha+5}{2}}}+ \|\eps\|_{H^{\frac{\alpha}{2}}}  \right) . 
\end{align}
By similar arguments, we deduce that: 
\begin{align}\label{eq:zpt2}
    \left\lvert \dot{z}_2 - \mu_2 +\delta(\Gamma) \right\rvert \int\left(\partial_yQ \right)^2 \leq C \left( \frac{1}{z^{\alpha+1}} \left\lvert \dot{\mu}_1-b(z) \right\rvert + \left\lvert \dot{z}_1 -\mu_1 +\beta(\Gamma) \right\rvert + \frac{1}{z^{\frac{3\alpha+5}{2}}}+ \|\eps\|_{H^{\frac{\alpha}{2}}}  \right). 
\end{align}
Then, by adding \eqref{eq:zpt1} and \eqref{eq:zpt2}, we obtain:
\begin{align} \label{eq:zpt1_2_bis}
     \vert \dot{z}_1 - \mu_1 +\beta(\Gamma) \vert +\vert \dot{z}_2 - \mu_2 +\delta(\Gamma) \vert
         \leq&  C\left(\frac{1}{z^{\alpha+1}} \left( \left\lvert \dot{\mu}_1-b(z) \right\rvert + \left\lvert \dot{\mu}_2+b(z) \right\rvert \right) + \frac{1}{z^{\frac{3\alpha+5}{2}}}+ \|\eps\|_{H^{\frac{\alpha}{2}}} \right).
\end{align}
Gathering \eqref{eq:mupt1_2_bis} and \eqref{eq:zpt1_2_bis}, we obtain \eqref{eq:zpt1_2}, and 
\begin{align*}
      \sum_{i=1}^2\vert \dot{\mu}_i + (-1)^i b(z) \vert \leq& C\left( \frac{1}{z^{\frac{3\alpha+5}{2}}} + \left(|\dot{\mu}_1|+|\dot{\mu}_2|+\frac{1}{z^{\alpha+1}}\right)\|\eps\|_{H^{\frac{\alpha}{2}}} + \|\eps\|_{H^{\frac{\alpha}{2}}}^2 \right).
\end{align*}
Since $|\dot{\mu_i}|\leq |\dot{\mu}_i + (-1)^ib(z)| + b(z)$, by applying the former inequality and \eqref{defi:b}, we conclude \eqref{eq:mupt1_2}. 
\end{proof}

\subsection{Monotonicity} \label{monotonicity}

We define:
\begin{align}\label{defi:phi}
    \phi(y)=\left(\int_{-\infty}^{+\infty} \frac{ds}{\langle s \rangle^{1+\alpha}}\right)^{-1}&\int_{y}^{+\infty} \frac{ds}{\langle s \rangle^{1+\alpha}},
\end{align}
and
\begin{align}
    \phi_{1}(t,y) := \frac{1-\phi(y)}{(1+\mu_1(t))^2} + \frac{ \phi(y)}{(1+\mu_2(t))^2} \quad \text{and} \quad \phi_{2}(t,y) :=& \frac{\mu_1(t)}{(1+\mu_1(t))^2} (1-\phi(y)) + \frac{\mu_2(t)}{(1+\mu_2(t))^2} \phi(y).
\end{align}
Let $A>0$, we define the rescaled functions:
\begin{align}
\phi_A(y)=\phi\left(\frac{y}{A}\right),\quad 
       \phi_{1,A}(t,y):=\phi_1\left(t,\frac{y}{A}\right), \quad \phi_{2,A}(t,y):=\phi_2\left(t,\frac{y}{A}\right),
\end{align}
the derivatives by:
\begin{align}\label{Phi'}
\Phi(y)=\sqrt{|\phi'(y)|}, \quad \Phi_i(t,y)=\sqrt{|\phi'_i(t,y)|}, \quad \Phi_{i,A}(t,y)=\Phi_{i}\left(t,\frac{y}{A} \right).
\end{align}

By direct computation, we have:
\begin{align}\label{eq:equiv_Phi}
    \Phi_1(y) = \frac{c}{\langle y \rangle^{\frac{1+\alpha}{2}}} \frac{\mu^{\frac{1}{2}}(2+\bar{\mu})^{\frac{1}{2}}}{(1+\mu_1)(1+\mu_2)} \quad \text{and} \quad \Phi_2(y) = \frac{c}{\langle y \rangle^{\frac{1+\alpha}{2}}} \frac{\mu^{\frac{1}{2}}(1-\mu_1\mu_2)^{\frac{1}{2}}}{(1+\mu_1)(1+\mu_2)}.
\end{align}

We also define the functional:
\begin{align}\label{defi:F}
    F(t)=\displaystyle\int\left(\frac{\eps|D|^{\alpha}\eps}{2} + \frac{\eps^2}{2} - \frac{\left(V+\eps\right)^4}{4} + \frac{V^4}{4} + V^3\eps- S\eps\right)\phi_{1,A} + \frac{\eps^2}{2}\phi_{2,A}.
\end{align}

We claim the following theorem that will help us to get the estimate \eqref{boot:eps} on the error $\eps$.

\begin{theo}\label{theo:functional}
The following bound on the functional holds:
\begin{align}
    F(t)\leq Ct^{-\frac{7\alpha+9}{2(\alpha+3)}}.
\end{align}
\end{theo}

\subsubsection{Preliminary results}
To get the monotinicity properties of the modified energy, we need to recall a result from Lemma 6 and Lemma 7 from \cite{kenig2011local} and Lemma 3.2 from \cite{eychenne2021asymptotic}.
\begin{lemm}\label{commG}
	Let  $\alpha\in]0,2[$. In the symmetric case, there exists $C>0$ such that:
	\begin{eqnarray}
		\left|\displaystyle\int \left(|D|^{\alpha} u\right) u\Phi^2_{j,A}-\int\left(|D|^{\frac{\alpha}{2}}\left(u\Phi_{j,A} \right) \right)^2 \right|\leq \frac{C}{A^{\alpha}}\int u^2\Phi^2_{j,A}\label{est:comm1},
	\end{eqnarray}
	and
	\begin{eqnarray}
		\left|\displaystyle\int \left(|D|^{\alpha} u\right) \partial_x u \phi_{j,A}  +(-1)^{j+1}\frac{\alpha-1}{2}\int\left(|D|^{\frac{\alpha}{2}}\left(u\Phi_{j,A} \right) \right)^2  \right|\leq \frac{C}{A^{\alpha}}\int u^2\Phi^2_{j,A}\label{est:comm2},
	\end{eqnarray}
	for any $u\in \mathcal{S}(\mathbb{R})$,  $A>1$ and $j \in \{1,\cdots,N\}$. 
	
	In the non-symmetric case, there exists $C>0$ such that:
	\begin{align}  \label{est:ncomm1}
		\left| \displaystyle\int \left( \left( |D|^{\alpha} u \right)v - \left( |D|^{\alpha} v \right)u \right) \Phi^2_{j,A}  \right| \leq& 
		\begin{cases} \displaystyle\frac{C}{A^{\alpha}}\int \left(u^2 + v^2 \right) \Phi^2_{j,A} , & \text{if} \  \alpha\in]0,1], \\
	\displaystyle \frac{C}{A^{\frac{\alpha}{2}}} \int \left( u^2+ v^2 + \left(|D|^{\frac{\alpha}{2}}u\right)^2\right) \Phi^2_{j,A}, & \text{if} \ \alpha\in]1,2[ ,
	\end{cases}
		\end{align}
		and
	\begin{align} \label{est:ncomm2}
		\bigg| \displaystyle\int \left( \left( |D|^{\alpha} u \right)\partial_xv + \left( |D|^{\alpha} v \right)\partial_xu \right)& \phi_{j,A}  +(-1)^{j+1} (\alpha-1)\int |D|^{\frac{\alpha}{2}}\left(u\Phi_{j,A} \right) |D|^{\frac{\alpha}{2}}\left(v\Phi_{j,A} \right)   \bigg| \notag\\ 
		&\leq \begin{cases} \displaystyle\frac{C}{A^{\alpha}}\int \left(u^2 + v^2 \right) \Phi^2_{j,A} , & \text{if} \  \alpha\in]0,1], \\
	\displaystyle \frac{C}{A^{\frac{\alpha}{2}}} \int \left( u^2+ v^2 + \left(|D|^{\frac{\alpha}{2}}u\right)^2\right) \Phi^2_{j,A}, & \text{if} \ \alpha\in]1,2[ ,
	\end{cases},
	\end{align}
	for any $u,v\in \mathcal{S}(\mathbb{R})$,  $A>1$ and $j \in \{1,\cdots,N\}$. 
\end{lemm}
The estimates \eqref{est:comm1}-\eqref{est:comm2} are proved in Lemmas 6 and 7 in \cite{kenig2011local} for $\alpha\in[1,2]$. Observe however
that their proofs extend easily to the case $\alpha\in]0,2[$. Note also that while only one side of the
inequalities in \eqref{est:comm1}-\eqref{est:comm2} is stated in Lemmas 6 and 7 in \cite{kenig2011local} , both sides are actually proved. 

\begin{lemm}[\cite{eychenne2021asymptotic}, Lemma 3.3]
Let $0\leq\alpha\leq2$. For all $u\in\mathcal{S}(\R)$, we have that:
\begin{align}\label{eq:interchange_Dalpha}
    \left\vert\int \left(|D|^{\frac{\alpha}{2}}(u\Phi_{1,A})\right)^2 - \left(|D|^{\frac{\alpha}{2}} u\right)^2\Phi^2_{1,A}\right\vert \leq \frac{C}{A^{\frac{\alpha}{2}}}\int \left( u^2 + (|D|^{\frac{\alpha}{2}}u)^2 \right)\Phi_{1,A}^2.
\end{align}
\end{lemm}

The following estimates are proved in Appendix \ref{proof:preliminary}.

\begin{lemm}\label{est:commsimpl}
For $\alpha\in ]0,2[$, then for all $u\in\mathcal{S}(\R)$ we have that:
\begin{align}
    \bigg\|\left[ |D|^{\alpha} , \Phi_{j,A} \right]u \bigg\|_{L^2}^2\leq \begin{cases}\displaystyle\frac{C}{A^{2\alpha}} \int u^2\Phi^2_{j,A}, \quad \text{ if } \quad  \alpha \in ]0,1] \\
    \displaystyle\frac{C}{A^{\alpha}} \int \left(u^2+ \left(|D|^{\frac{\alpha}{2}}u\right)^2\right)\Phi^2_{j,A},\quad \text{ if }\quad  \alpha \in ]1,2] \end{cases}
\end{align}
\end{lemm}

\begin{lemm}\label{esttc}
	Let $\alpha\in]0,2[$, then for all $u\in\mathcal{S}(\R)$ there exists  $C>0$ such that: 
	\begin{align} \label{est:esttc}
		\bigg|\int |D|^{\alpha}\left(u\Phi_{j,A} \right)& \left((|D|^{\alpha}u)\Phi_{j,A} \right)   -\int \left(|D|^{\alpha} u\right)^2\Phi^2_{j,A} \bigg|\notag\\ \leq&  \frac{C}{A^{\frac{\alpha}{2}}} \int \left( u^2 + \left(|D|^{\frac{\alpha}{2}}u\right)^2 + \left(|D|^{\alpha} u\right)^2 \right) \Phi^2_{j,A}, 
	\end{align}
	for all $u\in \mathcal{S}(\R)$, $A>1$ and $j \in \{1,\cdots,N\}$. 
\end{lemm}

\begin{lemm}\label{lemm:commH1}
Let $1 \leq\alpha\leq2$. For all $u\in \mathcal{S}(\R)$, we have that: 
\begin{align*}
    \left\|[|D|^{\alpha},\phi_{1,A}]u \right\|_{L^2}\leq C \bigg|\frac{1}{(1+\mu_1)^2} - \frac{1}{(1+\mu_2)^2}\bigg|^{\frac{1}{2}}  \|u\Phi_{1,A} \|_{H^{1}}. 
\end{align*}
\end{lemm}

\begin{rema}
Notice that the scaling in $A$ is not coherent with the previous inequality. In the proof in the appendix, we establish this inequality in $H^{\frac{\alpha}{2}}(\R)$ and use at the very end the embedding $H^{\frac{\alpha}{2}}(\R) \subset H^1(\R)$.
\end{rema}

\begin{lemm}\label{lemm:commL2}
Let $0\leq \alpha\leq 2$. For all $u\in \mathcal{S}(\R)$, we have that: 
\begin{equation*}
\|[|D|^{\alpha},\sqrt{\phi_A}]u \|_{L^2} +  \|[|D|^{\alpha},\sqrt{1-\phi_A}]u \|_{L^2} \leq
    \left\{ \begin{aligned} & \frac{C}{A^{\alpha}} \|u\|_{L^2}, & \alpha\in(0,1]\\
    & \frac{C}{A^{\frac{\alpha}{2}}}\|u\|_{H^{\frac{\alpha}{2}}}, & \alpha\in(1,2]
    \end{aligned} \right. .
\end{equation*}
\end{lemm}

\subsubsection{Proof of the Theorem \ref{theo:functional}}

In this part, we study the functional $F$ defined in \eqref{defi:F}, dependent on the two functions $\phi_{1,A}$ and $\phi_{2,A}$. For sake of clearness, we drop the indices $A$ in this part only and denote those functions by $\phi_1$ and $\phi_2$. The parameter $A$ will appear explicitly when needed.

We recall the equation satisfied by $\eps$: 
\begin{align*}
    \partial_t\eps +\partial_y\left(-|D|^{\alpha}\eps - \eps + \left(\eps+V\right)^3 -V^3 \right)=-\cE_V.
\end{align*}
We differentiate in time the functional $F$ defined in \eqref{defi:F}, by using \eqref{definition:energy} we deduce that:
\begin{align*}
    \frac{d}{dt}F(t)=&\int (\partial_t\eps) \left( |D|^{\alpha}\eps + \eps - \left(\eps+V\right)^3 +V^3 -S\right)\phi_{1} +\frac{1}{2}\int\left( \eps|D|^{\alpha}\partial_t\eps - (\partial_t\eps) |D|^{\alpha}\eps \right)\phi_{1}\\
    & +\int \left(-(\partial_tV)\left( \left(V+\eps \right)^3-V^3-3V^2\eps \right)\phi_{1}+ (\partial_t\eps)\eps\phi_{2}\right) -\int (\partial_tS)\eps \phi_{1} \\
    +& \int\left(\frac{\eps|D|^{\alpha}\eps}{2} + \frac{\eps^2}{2} - \frac{\left(V+\eps\right)^4}{4} + \frac{V^4}{4} + V^3\eps- S\eps\right)\partial_t\phi_{1} + \int \frac{\eps^2}{2}\partial_t\phi_{2}\\
    =&I_1+\cdots + I_6.
\end{align*}

\noindent\textit{\textbf{Estimate on $I_1$:}}
Using integration by parts and the definition of $\cE_V$, we deduce that:
\begin{align*}
    I_1&=\frac{1}{2}\int\left(|D|^{\alpha}\eps+\eps-(V+\eps)^3+V^3 \right)^2\Phi_1^2 - \int \partial_y\left( |D|^{\alpha}\eps+\eps-(V+\eps)^3+V^3\right)S\phi_1 \\
    & \quad- \int \cE_V\left(|D|^{\alpha}\eps+\eps-(V+\eps)+V^3 -S \right)\phi_1\\
    &=\frac{1}{2}\int\left(|D|^{\alpha}\eps+\eps-(V+\eps)^3+V^3 \right)^2\Phi_1^2 - \int\left( \left(\overrightarrow{m} \cdot \overrightarrow{MV} + T \right) \left(|D|^{\alpha}\eps+\eps-(V+\eps)^3+V^3 \right)-TS \right)\phi_1 \\
    & \quad- \int S\left(|D|^{\alpha}\eps + \eps - \left( V+\eps \right)^3 +V^3  \right)\Phi_1^2 + \int \overrightarrow{m}\cdot\overrightarrow{MV}S\phi_1 + \frac{1}{2} \int S^2\Phi_1^2\\
    &=I_{1,1}+\cdots+I_{1,5}.
\end{align*}
We start with $I_{1,1}$. By direct computations, we get that: 
\begin{align*}
\frac{1}{2}&\int(|D|^{\alpha}\eps+\eps-(V+\eps)^3+V^3 )^2\Phi_1^2 - \frac{1}{2}\int\left( \left(|D|^{\alpha}\eps\right)^2 +\eps^2+ \left(-(V+\eps)^3+V^3\right)^2 \right)\Phi_1^2 \\
&= \int \left(|D|^{\alpha}\eps\right)\eps\Phi_1^2 + \int \left(|D|^{\alpha}\eps+\eps\right)\left(V^3-\left(V+\eps \right)^3 \right)\Phi_1^2=I_{1,1,1}+I_{1,1,2}.
\end{align*}

By using the estimate \eqref{est:comm1}, we obtain that: 
\begin{align*}
    \left\lvert I_{1,1,1}- \int\left(|D|^{\frac{\alpha}{2}}\left(\eps\Phi_1 \right) \right)^2\right\rvert \leq\frac{C}{A^{\alpha}}\int\eps^2\Phi_1^2.
\end{align*}
Since $V^3-\left(V+\eps \right)^3=-3V^2\eps - 3V\eps^2-\eps^3$, by applying Young's inequality, the bound on $V$ \eqref{V} and Cauchy-Schwarz' inequality, we have that: 
\begin{align*}
    \left\lvert I_{1,1,2} \right\rvert \leq \frac{C}{A^{\alpha}}\int \left( |D|^{\alpha}\eps\right)^2\Phi_1^2 +CA^{\alpha}\int\left(V^4\eps^2+\eps^4+\eps^6 \right)\Phi_1^2 + C\int \left(V^4\eps^2+\eps^3+\eps^4  \right)\Phi_1^2. 
\end{align*}
We recall that $2\mu_1=\mu+\bar{\mu}$, $2\mu_2=\bar{\mu}-\mu$ and $\alpha>1$. Therefore, using the bootstrap estimates \eqref{boot:eps}, \eqref{boot:mu} and \eqref{boot:mu_bar} and \eqref{Vphi}, we conclude for $I_{1,1}$ that: 
\begin{align*}
     I_{1,1} - \frac{1}{2}\int(\eps^2+(|D|^{\alpha}\eps)^2)\Phi_1^2 - \int\left(|D|^{\frac{\alpha}{2}}\left(\eps\Phi_1 \right) \right)^2  \geq  -\frac{C}{A^{\alpha}}\int \left(\eps^2+(|D|^{\alpha}\eps)^{2} \right)\Phi_1^2 - CA^{\alpha}t^{-\frac{3(3\alpha+5)}{2(\alpha+3)}}.
\end{align*}

Let us estimate $I_{1,2}$. By using the definition of $\overrightarrow{m}\cdot\overrightarrow{MV}$ in \eqref{defi:m_MV}, we obtain that:
\begin{align*}
\MoveEqLeft
  \int\overrightarrow{m}\cdot\overrightarrow{MV}  \left(|D|^{\alpha}\eps+\eps-(V+\eps)^3+V^3\right)\phi_1 \\
  & =  \sum_{i=1}^2\int \left( (-1)^i \dot{\mu_i}-b(z) \right)\Lambda R_i \left(|D|^{\alpha}\eps+\eps-(V+\eps)^3+V^3 \right)\phi_1 \\
  & \quad + \int \left( (\dot{z_1}-\mu_1+\beta(\Gamma)) \partial_y R_1 - (\dot{z_2}-\mu_2+\delta(\Gamma)) \partial_y R_2 \right) \left(|D|^{\alpha}\eps+\eps-(V+\eps)^3+V^3 \right)\phi_1\\
  &=J_1+J_2.
\end{align*}
Since $\frac{1}{1+\mu_i}\leq C$, we deduce that: 
\begin{align*}
    |J_1|\leq&\sum_{i=1}^2 |(-1)^i \dot{\mu_i}-b(z)|\bigg| \int \Lambda R_i \left(|D|^{\alpha}\eps+\eps-(V+\eps)^3+V^3 \right)\left( \phi_1 + (-1)^i\frac{1}{(1+\mu_i)^2} \right)\bigg|\\
    +& C \sum_{i=1}^2 |(-1)^i \dot{\mu_i}-b(z)|\bigg| \int \Lambda R_i \left(|D|^{\alpha}\eps+(1+\mu_i)\eps-R_i^2\eps \right) \bigg| \\
    +& C\sum_{i=1}^{2}|(-1)^i \dot{\mu_i}-b(z)|\bigg|\int \Lambda R_i\left( -\mu_i\eps-(V+\eps)^3 + V^3 -3R_i^2\eps\right)\bigg|=J_{1,1}+J_{1,2}+J_{1,3}.
\end{align*}
Thanks to the identity $(V+\eps)- V^3=\eps^3+3\eps^2V+3\eps V^2$, the fact $\alpha<2$, and by Cauchy-Schwarz' inequality, we get that: 
\begin{align*}
    J_{1,1}\leq C \sum_{i=1}^2|(-1)^i\dot{\mu_i}-b(z)|\bigg|\frac{1}{(1+\mu_1)^2}-\frac{1}{(1+\mu_2)^2}\bigg|\left\|\Lambda R_i (\phi-\delta_{2i})\right\|_{H^1}\|\eps\|_{H^{\frac{\alpha}{2}}}.
\end{align*}
Moreover, we recall $L\Lambda Q=-Q$ and since $\eps\perp R_i$, we deduce that:
\begin{align*}
    J_{1,2}= C\sum_{i=1}^2 |(-1)^i \dot{\mu_i}-b(z)|\bigg| \int R_i \eps \bigg|=0.
\end{align*}
Applying Cauchy-Schwarz' inequality, and Sobolev embedding $H^{\frac{1}{3}}(\R)\xhookrightarrow{}L^6(\R)$ , we have that: 
\begin{align*}
    J_{1,3}\leq C\sum_{i=1}^2|(-1)^i \dot{\mu_i}-b(z)|\left(\left(\vert \mu_i \vert +\left\|(V^2-R_i^2)\Lambda R_i\right\|_{L^2}\right)\|\eps\|_{L^2} + \|\eps\|_{L^2}^2 + \| \eps\|_{H^{\frac{\alpha}{2}}}^3 \right).
\end{align*}
Now, let us estimate $J_2$. We focus on the first term of $J_2$ with $\partial_y R_1$, the second is similar. We decompose this term into: 
\begin{align*}
\MoveEqLeft
    (\dot{z_1}-\mu_1+\beta(\Gamma)) \int \partial_y R_1 \left(|D|^{\alpha}\eps+\eps-(V+\eps)^3+V^3 \right)\left(\phi_1-\frac{1}{(1+\mu_1)^2}\right) \\
    & \quad +\frac{\dot{z_1}-\mu_1+\beta(\Gamma)}{(1+\mu_1)^2} \left(\int \partial_y R_1 \left(|D|^{\alpha}\eps+\eps-3R_1^2 \eps\right) 
    + \int \partial_yR_1 \left(V^3-(V+\eps)^3+3R_1^2\eps\right)\right) \\
     & \quad \quad =J_{2,1}+J_{2,2}+J_{2,3}.
\end{align*}
By applying the Cauchy-Schwarz' inequality and Sobolev embedding $H^{\frac{1}{3}}(\R)\xhookrightarrow{}L^6(\R)$, we obtain that: 
\begin{align*}
    |J_{2,1}|\leq C|\dot{z_1}-\mu_1+ \beta(\Gamma)|\bigg| \frac{1}{(1+\mu_1)^2}-\frac{1}{(1+\mu_2)^2}\bigg|\left(\left\|\partial_y R_1 \phi \right\|_{H^1} \|\eps\|_{H^{\frac{\alpha}{2}}} + \|\eps\|_{L^2}^2 +  \|\eps \|^3_{H^{\frac{\alpha}{2}}} \right).
\end{align*}
Since $\eps\perp \partial_yR_i$ and $LQ'=0$, we deduce that:
\begin{align*}
    J_{2,2}=0.
\end{align*}
Moreover, by Cauchy-Schwarz inequality and Sobolev embedding , $H^{\frac{1}{6}}(\R)\xhookrightarrow{}L^3(\R)$ we have that:
\begin{align*}
   |J_{2,3}|\leq C |\dot{z_1}-\mu_1+\beta(\Gamma)|\left(\|\eps\|_{L^2} \|(V^2-R_1^2) \partial_yR_1\|_{L^2} + \|\eps\|_{L^2}^2 + \| \eps\|_{H^{\frac{\alpha}{2}}}^3 \right). 
\end{align*}

By Cauchy-Schwarz' inequality, we have that: 
\begin{align*}
   \left|\int T \left(|D|^{\alpha}\eps+\eps-(V+\eps)^3+V^3 -S \right)\phi_1 \right|\leq C\|T\phi_1\|_{H^{\frac{\alpha}{2}}}\left(\|\eps\|_{H^{\frac{\alpha}{2}}}+ \|\eps\|_{H^{\frac{\alpha}{2}}}^2 + \|\eps\|_{H^{\frac{\alpha}{2}}}^3 +\|S\|_{L^2} \right).
\end{align*}

Gathering those identities, and using the estimate on $T$ \eqref{est:T}, the estimates on the solitary waves \eqref{eq:partial_y_Ri}, \eqref{eq:Ri_Phi}, \eqref{Vphi}, the bootstrap estimates \eqref{boot:eps}, \eqref{boot:mu}, \eqref{boot:mu_bar} and the equation on $\dot{\mu_i}$ \eqref{eq:mupt1_2} and $\dot{z_i}$ \eqref{eq:zpt1_2}, we get that: 
\begin{align*}
    |I_{1,2}|\leq Ct^{-\frac{3(3\alpha+5)}{2(\alpha+3)}}.
\end{align*}
Let us estimate $I_{1,3}$.
By Cauchy-Schwarz inequality, the estimate on $S$ \eqref{est:S} and the bootstrap estimate on $\eps$ \eqref{boot:eps}, we obtain that: 
\begin{align*}
    |I_{1,3}|\leq C\left\|S\Phi_1^2\right\|_{H^{\frac{\alpha}{2}}}\left(\|\eps\|_{H^{\frac{\alpha}{2}}} + \left\| \left(V+\eps\right)^3 - V^3 \right\|_{L^2} \right)\leq Ct^{-\frac{3(3\alpha+5)}{2(\alpha+3)}}.
\end{align*}
Using the definition of $\overrightarrow{m}\cdot\overrightarrow{MV}$,  the estimate on $\dot{\mu_i}$ \eqref{eq:mupt1_2}, $\dot{z_i}$ \eqref{eq:zpt1_2} and the estimate on $S$ \eqref{est:S}, we deduce that: 
\begin{align*}
    |I_{1,4}|\leq
     C\|S\|_{L^2} \left(\sum_{i=1}^2 |(-1)^i \dot{\mu_i}-b(z)|+ | \dot{z_1}-\mu_1-\beta(\Gamma)|+ | \dot{z_2}-\mu_2-\delta(\Gamma)| \right)\leq Ct^{-\frac{3(3\alpha+5)}{2(\alpha+3)}}.
\end{align*}
Finally, by the estimate on $S$ \eqref{est:S}:
\begin{align*}
    \vert I_{1,5} \vert \leq C t^{-\frac{3(3\alpha+5)}{2(\alpha+3)}}.
\end{align*}

\noindent\textit{Conclusion:}

\begin{align*}
     I_1 - \frac{1}{2}\int(\eps^2+(|D|^{\alpha}\eps)^2)\Phi_1^2 - \int\left(|D|^{\frac{\alpha}{2}}\left(\eps\Phi_1 \right) \right)^2 \geq -\frac{C}{A^{\alpha}}\int \left(\eps^2+(|D|^{\alpha}\eps)^{2} \right)\Phi_1^2 - CA^{\alpha}t^{-\frac{3(3\alpha+5)}{2(\alpha+3)}}.
\end{align*}

\noindent\textit{\textbf{Estimate on $I_2$:}}
From the equation of $\eps$, since $\phi_1$ is decreasing and integration by parts, we deduce that:
\begin{align*}
\MoveEqLeft
   2 I_2= -\frac{1}{2}\int\left(|D|^{\alpha}\eps\right)^2\Phi^2_{1} -\int \partial_y\eps |D|^{\alpha}\left(|D|^{\alpha}\eps +2\eps \right)\phi_{1} + \int \eps |D|^{\alpha}\left(|D|^{\alpha}\eps +\eps \right)\Phi^2_{1} \\ 
   & + \int \left( \cE_V|D|^{\alpha}\eps - \eps|D|^{\alpha}\cE_V \right)\phi_1
   + \int \left(\eps |D|^{\alpha}\partial_y\left( - (V+\eps)^3 + V^{3} \right)-\partial_y\left(-(V+\eps)^3 + V^3 \right) (|D|^{\alpha}\eps)\right) \phi_1 \\
   & =I_{2,1}+\cdots I_{2,5}.
\end{align*}
Let us estimate $I_{2,2}$ and $I_{2,3}$. Using the commutator estimates in the non-symmetric case \eqref{est:ncomm1}, \eqref{est:ncomm2} with $v=|D|^{\alpha}\eps$, the commutator estimates in the symmetric case \eqref{est:comm1}, \eqref{est:comm2}  , and Lemma \ref{esttc} we get that: 
\begin{align*}
    \left\vert I_{2,2}+I_{2,3} -\alpha\int\left(|D|^{\frac{\alpha}{2}}\left(\eps\Phi_{1}\right)\right)^2  - \left(\alpha + \frac{1}{2}\right)\int \left(|D|^{\alpha}\eps\right)^2\Phi^2_{1} \right\vert &\leq \frac{C}{A^{\frac{\alpha}{2}}}\int \left(\eps^2 + \left(|D|^{\frac{\alpha}{2}}\eps\right)^2 +\left(|D|^{\alpha}\eps\right)^2\right) \Phi^2_{1}  
\end{align*}
From Cauchy-Schwarz inequality and Lemma \ref{lemm:commH1}, we get that:   
\begin{align*}
    |I_{2,4}|
    & = \left\lvert\int \eps\left(|D|^{\alpha}\left(\cE_V \phi_1\right) -|D|^{\alpha}\left(\cE_V \right)\phi_1 \right)\right\rvert \leq \|\eps\|_{L^2}\|[|D|^{\alpha},\phi_1]\cE_V \|_{L^2}\\
    &\leq C \bigg|\frac{1}{(1+\mu_1)^2} - \frac{1}{(1+\mu_2)^2}\bigg| \|\eps\|_{L^2}\left\|\cE_V\sqrt{\phi'} \right\|_{H^1}\\
    &\leq C \bigg|\frac{1}{(1+\mu_1)^2} - \frac{1}{(1+\mu_2)^2}\bigg| \|\eps\|_{L^2}\left(\left\|\overrightarrow{m}.\overrightarrow{MV} \sqrt{\phi'} \right\|_{H^1} + \left\|S \right\|_{H^1} +\left\|T \right\|_{H^1} \right).
\end{align*}
Therefore, by using the estimates on $\dot{\mu_i}$ \eqref{eq:mupt1_2}, on $\dot{z_i}$ \eqref{eq:zpt1_2}, the estimates on $S$ \eqref{est:S}, $T$ \eqref{est:T}, the interaction between $\partial_y R_i$ or $\Lambda R_i$ and $\Phi$ \eqref{eq:Ri_Phi} and the bootstrap estimates \eqref{boot:eps}-\eqref{boot:mu_bar}, we have that: 
\begin{align*}
    |I_{2,4}|\leq t^{-\frac{3(3\alpha+5)}{2(\alpha+3)}}.
\end{align*}
Now, we estimate $I_{2,5}$. Note that: 
$$
(V+\eps)^3 - V^3 = 3V^2\eps + 3V\eps^2 + \eps^3.
$$
Then, we decompose $I_{2,5}$ as: 
\begin{align*}
I_{2,5} &=\int \left(\partial_y\left(3V^2\eps \right) (|D|^{\alpha}\eps)-\eps |D|^{\alpha}\partial_y\left( 3V^2\eps \right)\right) \phi_1 +\int \left(\partial_y\left(3V\eps^2 \right) (|D|^{\alpha}\eps)-\eps |D|^{\alpha}\partial_y\left( 3V\eps^2 \right)\right) \phi_1 \\
&+ \int \left(\partial_y\left(\eps^3 \right) (|D|^{\alpha}\eps)-\eps |D|^{\alpha}\partial_y\left( \eps^3 \right)\right) \phi_1.
\end{align*}
Let $v\in\{3\eps V^2,3\eps^2V, \eps^3\}$  . Using integration by parts, the commutator estimates in the non-symmetric case \eqref{est:ncomm1} and \eqref{est:ncomm2}, we get that: 
\begin{align*}
\left\lvert\int \left(\partial_yv (|D|^{\alpha}\eps)-\eps |D|^{\alpha}\partial_y v\right) \phi_1\right\rvert&= \left\lvert \int \partial_yv (|D|^{\alpha}\eps)\phi_1 + \int \partial_y\eps (|D|^{\alpha}v) \phi_1  -  \int \eps (|D|^{\alpha}v) \Phi_1^2\right\rvert \\
&\leq (\alpha-1) \left\lvert\int|D|^{\frac{\alpha}{2}}\left(v\Phi_1 \right)|D|^{\frac{\alpha}{2}}\left(\eps\Phi_1 \right)\right\rvert + \left\lvert\int(|D|^{\alpha}\eps)v\Phi_1^2\right\rvert\\
&+\frac{C}{A^{\frac{\alpha}{2}}}\int \left(\eps^2 + v^2 + (|D|^{\frac{\alpha}{2}}\eps)^2 \right)\Phi_1^2. 
\end{align*}
Moreover, from Young's inequality, we obtain that:
\begin{align*}
\left|\int(|D|^{\alpha}\eps)v\Phi_1^2 \right|\leq \frac{C}{A^{\alpha}}\int (|D|^{\alpha}\eps)^{2}\Phi_1^2 + CA^{\alpha} \int v^2\Phi_1^2.    
\end{align*}
By using Young's inequality and \eqref{eq:interchange_Dalpha}, we deduce that: 
\begin{align*}
    \left|\int|D|^{\frac{\alpha}{2}}\left(v\Phi_1 \right)|D|^{\frac{\alpha}{2}}\left(\eps\Phi_1 \right) \right|\leq \frac{C}{A^{\alpha}}\int (\eps^2 + (|D|^{\alpha}\eps)^2)\Phi_1^2 + CA^{\alpha}\int v^2\Phi_1^2
\end{align*}
By the Sobolev's embeddings, $H^{\frac{1}{4}}(\R)\xhookrightarrow{}L^4(\R)$, $H^{\frac{1}{3}}(\R)\xhookrightarrow{}L^6(\R)$, we obtain that: \begin{align*}
|I_{2,5}|\leq CA^{\alpha}\left(\| V^4\Phi_1^2 \|_{L^\infty}\|\eps \|_{L^2}^{2} +\| V^2\Phi_1^2 \|_{L^\infty}\|\eps \|^4_{H^{\frac{\alpha}{2}}} + \|\eps \|^6_{H^{\frac{\alpha}{2}}}\right) + \frac{C}{A^{\frac{\alpha}{2}}}\int \left(\eps^2 + (|D|^{\frac{\alpha}{2}}\eps)^2 + (|D|^{\alpha}\eps)^2 \right)\Phi_1^2.
\end{align*}  
Moreover, applying the estimate \eqref{Vphi}, \eqref{boot:eps} and \eqref{boot:z}, we get that: 
\begin{align*}
\left\vert I_{2,5}\right\vert\leq CA^{\alpha}t^{-\frac{3(3\alpha+5)}{2(\alpha+3)}} + \frac{C}{A^{\frac{\alpha}{2}}}\int \left(\eps^2 + \left(|D|^{\frac{\alpha}{2}}\eps \right)^2 +\left(|D|^{\alpha}\eps \right)^2 \right)\Phi_1^2.
\end{align*}
\noindent\textit{Conclusion:}
\begin{align}\label{est:I2}
    \left\vert I_2 - \frac{\alpha}{2}\int \left(|D|^{\alpha}\eps\right)^2\Phi_1^2 - \frac{\alpha}{2}\int\left(|D|^{\frac{\alpha}{2}}\left(\eps \Phi_1\right) \right)^2 \right\vert \leq&   CA^{\alpha}t^{-\frac{3(3\alpha+5)}{2(\alpha+3)}}
    + \frac{C}{A^{\frac{\alpha}{2}}}\int \left(\eps^2 + (|D|^{\frac{\alpha}{2}}\eps)^2 + (|D|^{\alpha}\eps)^2 \right)\Phi_1^2.
\end{align}

\noindent\textit{\textbf{Estimate on $I_3$:}}
We decompose $I_3$ as: 
\begin{align*}
    I_3&= -3\int \partial_y(V^2\eps)\eps\phi_2+\partial_tV  V\eps^2\phi_1 -\int  \partial_y\left( 3V\eps^2 + \eps^3\right)\eps \phi_2 + \partial_tV\eps^3\phi_1 \\
    &+ \int \partial_y\left(|D|^{\alpha}\eps +\eps \right)\eps\phi_2 - \int \cE_V\eps\phi_2= I_{3,1} + I_{3,2}+ I_{3,3} +I_{3,4}.
\end{align*}
By adding $0$ and integrating by part, we deduce that: 
\begin{align*}
    I_{3,1}&=3\int \partial_y R_1(\phi_2-\dot{z_1}\phi_1)V \eps^2 + 3 \int \partial_y R_2(\dot{z_2}\phi_1 - \phi_2 )V \eps^2 - 3\int\left( \partial_yV + \partial_yR_1 - \partial_y R_2 \right) V\eps^2\phi_2 \\
    &-  3\int\left( \partial_tV -\dot{z_1} \partial_yR_1 + \dot{z_2} \partial_y R_2 \right) V\eps^2\phi_1 +\frac{3}{2}\int V^2\eps^2\Phi_2^2=I_{3,1,1}+\cdots +I_{3,1,5}.
\end{align*}
Using the definition of $\phi_1$ and $\phi_2$, we obtain that: 
\begin{align*}
    |I_{3,1,1}|&= \bigg|3\frac{\mu_1-\dot{z_1}}{1+\mu_1}\int\partial_y R_1 V\eps^2(1-\phi) + 3\frac{\dot{z_1}- \mu_2}{1+\mu_2}\int\partial_y R_1 V\eps^2\phi  \bigg|\\
    &\leq C \|\eps\|_{L^2}^2 \left( |\dot{z_1}-\mu_1| + \left(|\dot{z_1}|+|\mu_2| \right) \|\partial_y R_1 V\phi\|_{L^\infty}\right).
\end{align*}
Using the same argument, we deduce that: 
\begin{align*}
    |I_{3,1,2}|\leq  C \|\eps\|_{L^2}^2 \left( |\dot{z_2}-\mu_2| + \left(|\dot{z_2}|+|\mu_1| \right) \|\partial_y R_2 V(1-\phi)\|_{L^\infty}\right).
\end{align*}
Using the definition of $V$ and $\phi_2$:
\begin{align*}
    |I_{3,1,3}|\leq C\|\eps\|_{L^2}^2 \left(|\mu_1|+|\mu_2|\right)\|\partial_y\left(P_2 -P_1 +bW\chi \right) \|_{L^{\infty}}
\end{align*}
and 
\begin{align*}
    |I_{3,1,4}|\leq C \|\eps\|_{L^2}^2 \left(\left(|\dot{\mu_1}|+|\dot{\mu_2}|\right) +\left(|\mu_1|+|\mu_2|\right)\|\partial_t\left(P_2 -P_1 +bW\chi \right) \|_{L^{\infty}}\right).
\end{align*}

Gathering these identities, and using the bootstrap hypothesis, the time estimate of the different terms and \eqref{partialyPibWchi} and \eqref{partialtPibWchi}, we conclude that: 
\begin{align*}
   \vert I_{3,1} \vert \leq Ct^{-\frac{3(3\alpha+5)}{2(\alpha+3)}}.
\end{align*}
For $I_{3,2}$, using integration by parts and Sobolev embedding, and the bootstrap hypothesis, we deduce that: 
\begin{align*}
    |I_{3,2}| & \leq C\left(\|\eps\|_{H^{\frac{\alpha}{2}}}^3\|\partial_t V\|_{L^{\infty}} + (|\mu_1| + |\mu_2|)(\| V \|_{L^{\infty}}+\|\partial_y V \|_{L^{\infty}})\|\eps\|_{H^{\frac{\alpha}{2}}}^3+(|\mu_1| + |\mu_2|)\|\eps\|_{H^{\frac{\alpha}{2}}}^4 \right) \\
    & \leq Ct^{-\frac{3(3\alpha+5)}{2(\alpha+3)}}. 
\end{align*}
Using integration by part, the commutator estimates in the symmetric case \eqref{est:comm1} and \eqref{est:comm2}, and since $\partial_y\phi_2<0$, we obtain that: 
\begin{align*}
    I_{3,3}\geq -\frac{\alpha+1}{2}\int\left(|D|^{\frac{\alpha}{2}}\left(\eps\Phi_2 \right)\right)^2 -\left(\frac{1}{2}+\frac{C}{A^{\alpha}} \right)\int \eps^2\Phi_2^2.
\end{align*}
Moreover with \eqref{eq:equiv_Phi}, we have: 
\begin{align*}
    \Phi_2^2=\left\lvert\frac{\mu_1\mu_2-1}{2+\mu_1+\mu_2}\right\rvert\Phi_1^2.
\end{align*}
Then, we get that:
\begin{align*}
    I_{3,3}\geq -\frac{\alpha+1}{2}\frac{1-\mu_1\mu_2}{2+\mu_1+\mu_2}\int\left(|D|^{\frac{\alpha}{2}}\left(\eps\Phi_1 \right)\right)^2 -\left(\frac{1}{2}+\frac{C}{A^{\alpha}} \right)\frac{1-\mu_1\mu_2}{2+\mu_1+\mu_2}\int \eps^2\Phi_1^2.
\end{align*}
Since $\frac{1-\mu_1\mu_2}{2+\mu_1+\mu_2}\leq \frac{3}{4}$ by \eqref{boot:mu} and \eqref{boot:mu_bar}, we deduce that: 
\begin{align*}
   I_{3,3} \geq - \frac{3(\alpha+1)}{8}\int\left(|D|^{\frac{\alpha}{2}}\left(\eps\Phi_1 \right)\right)^2 - \left(\frac{3}{8}+\frac{C}{A^{\alpha}} \right)\int \eps^2\Phi_1^2.
\end{align*}
Let us estimate the last term of $I_3$. Using the definition of $\cE_V$ and Cauchy-Schwarz inequality, we have that:
\begin{align*}
    |I_{3,4}|\leq C (|\mu_1|+\vert \mu_2 \vert )\|\eps\|_{L^2}\left(\|\partial_yS\|_{L^2} + \|T\|_{L^2} \right) + \left\vert \int \overrightarrow{m}\cdot\overrightarrow{MV}\eps \phi_2\right\vert.
\end{align*}
Using the definition of $\overrightarrow{m}\cdot\overrightarrow{MV}$ and the orthogonality condition $\eps\perp \partial_y R_i$, we deduce that: 
\begin{align*}
    \bigg|\int \overrightarrow{m}\cdot\overrightarrow{MV}\eps \phi_2\bigg| \leq C \|\eps\|_{L^2} (|\mu_1|+|\mu_2|) \left( \sum_{i=1}^2|(-1)^{i}\dot{\mu_i}-b(z)| + |\dot{z_i}-\mu_i|\left\|\partial_yR_i(\phi - \delta_{2,i}) \right\|_{L^2} \right).
\end{align*}
Therefore with \eqref{eq:partial_y_Ri}, we get that: 
$$|I_{3,4}|\leq Ct^{-\frac{3(3\alpha+5)}{2(\alpha+3)}}.$$
\noindent\textit{Conclusion:}
\begin{align*}
     I_3 \geq  -\frac{3(\alpha+1)}{8}\int\left(|D|^{\frac{\alpha}{2}}\left(\eps\Phi_1 \right)\right)^2 - \left(\frac{3}{8} +\frac{C}{A^{\alpha}}\right)\int \eps^2\Phi_1^2 -Ct^{-\frac{3(3\alpha+5)}{2(\alpha+3)}}.
\end{align*}

\noindent\textit{\textbf{Estimate on $I_4$:}}
Applying Cauchy-Schwarz inequality and the estimate on the time derivative of $S$ \eqref{est:dtS}, we obtain that: 
\begin{align*}
    |I_4|\leq C\|\partial_tS\|_{L^2}\|\eps\|_{L^2}\leq Ct^{-\frac{3(3\alpha+5)}{2(\alpha+3)}}.
\end{align*}
\noindent\textit{\textbf{Estimate on $I_5$:}}
First, note by direct computation, we have: 
\begin{align*}
    |\partial_t\phi_1|=\bigg|\frac{2\dot{\mu_1}}{(1+\mu_1)^3}(1-\phi) +\frac{2\dot{\mu_2}}{(1+\mu_2)^3} \phi \bigg|\leq C\left(|\dot{\mu_1}| + |\dot{\mu_2}| \right).
\end{align*}
Then, by the Sobolev embedding $H^{\frac{1}{3}}(\R)\xhookrightarrow{}L^6(\R)$ and $H^{\frac{1}{4}}(\R)\xhookrightarrow{}L^4(\R)$, we deduce that: \begin{align*}
    \bigg| \int \left(\frac{\eps^2}{2} - \frac{\left(V+\eps\right)^4}{4}+\frac{V^4}{4} + V^3\eps \right)\partial_t\phi_1\bigg|\leq C\left(|\dot{\mu_1}| + |\dot{\mu_2}| \right) \left(\|\eps \|_{H^{\frac{\alpha}{2}}}^2 + \|\eps \|_{H^{\frac{\alpha}{2}}}^3 +\|\eps\|_{H^{\frac{\alpha}{2}}}^4 \right).
\end{align*}
Moreover, by Cauchy-Schwarz inequality, we get: 
\begin{align*}
    \bigg| \int S\eps\partial_t\phi_1 \bigg|\leq C \left(|\dot{\mu_1}| + |\dot{\mu_2}| \right)\|S\|_{L^2}\|\eps\|_{L^2}.
\end{align*}
Now, let us estimate the first term in $I_5$. By direct computations, we have that:
\begin{align*}
    \int \eps|D|^{\alpha}\eps\partial_t\phi_1&= -\frac{2\dot{\mu_1}}{(1+\mu_1)^3} \left(\int D^{\frac{\alpha}{2}}\eps [|D|^{\frac{\alpha}{2}},(1-\phi)]\eps  + \int \left(|D|^{\frac{\alpha}{2}}\eps\right)^2(1-\phi)\right)\\
    &-\frac{2\dot{\mu_2}}{(1+\mu_2)^3}  \left(\int D^{\frac{\alpha}{2}}\eps [|D|^{\frac{\alpha}{2}},\phi]\eps  + \int \left(|D|^{\frac{\alpha}{2}}\eps\right)^2\phi\right).
\end{align*}
Using Lemma \ref{lemm:commL2}, we deduce that: 
\begin{align*}
    \bigg|\int \eps|D|^{\alpha}\eps\partial_t\phi_1  \bigg| \leq C(|\dot{\mu_1}| + |\dot{\mu_2}|) \|\eps\|_{H^{\frac{\alpha}{2}}}^2.
\end{align*}
\noindent\textit{Conclusion:} 
\begin{align*}
    |I_5|\leq C(|\dot{\mu_1}| + |\dot{\mu_2}|)\left( \| S \|_{L^2}\|\eps\|_{H^{\frac{\alpha}{2}}}+ \|\eps\|_{H^{\frac{\alpha}{2}}}^2  +\|\eps\|_{H^{\frac{\alpha}{2}}}^3 +\|\eps\|_{H^{\frac{\alpha}{2}}}^4\right)\leq Ct^{-\frac{3(3\alpha+5)}{2(\alpha+3)}}.
\end{align*}
\noindent\textit{\textbf{Estimate on $I_6$:}}
By definition of $\phi_2$, we obtain that: 
\begin{align*}
    \left\lvert\partial_t\phi_2\right\rvert\leq C \left( \left\lvert \dot{\mu_1}\right\rvert + \left\lvert \dot{\mu_2}\right\rvert\right).
\end{align*}
then, by using the estimate on $\dot{\mu_i}$ \eqref{eq:mupt1_2}, the bootstrap estimates \eqref{boot:eps}, \eqref{boot:z}, we have that:
\begin{align*}
    |I_6|\leq Ct^{-\frac{3(3\alpha+5)}{2(\alpha+3)}}.
\end{align*}
Gathering the estimates on $I_1,...,I_6$, we obtain that: 
\begin{align*}
    \frac{d}{dt}F(t)\geq& \frac{\alpha+1}{2}\int\left(|D|^{\alpha}\eps \right)^2\Phi^2_1 + \left(1+\frac{\alpha}{2} -\frac{3(\alpha+1)}{8} \right) \int(|D|^{\frac{\alpha}{2}}(\eps\Phi_1))^2 
    +\left( \frac{1}{2} -\frac{3}{8}  \right) \int\eps^2\Phi^2_1\\
    &-\frac{C}{A^{\frac{\alpha}{2}}}\int (\eps^2 + (|D|^{\frac{\alpha}{2}}\eps)^2 + (|D|^{\alpha}\eps)^2)\Phi_1^2 - CA^{\alpha}t^{-\frac{3(3\alpha+5)}{2(\alpha+3)}}.
\end{align*}
To compare the quantities $\int\left(|D|^{\alpha}\eps \right)^2\Phi^2_1$ and $\int\left(|D|^{\alpha}\eps\Phi_1 \right)^2$ we use \eqref{eq:interchange_Dalpha}, thus we have:
\begin{align*}
    \frac{d}{dt}F(t)\geq& \frac{\alpha+1}{2}\int\left(|D|^{\alpha}\eps \right)^2\Phi^2_1 + \left(1+\frac{\alpha}{2} -\frac{3(\alpha+1)}{8} \right) \int(|D|^{\frac{\alpha}{2}}\eps)^2\Phi_1^2
    +\left( \frac{1}{2} -\frac{3}{8}  \right) \int\eps^2\Phi^2_1\\
    &-\frac{C}{A^{\frac{\alpha}{2}}}\int (\eps^2 + (|D|^{\frac{\alpha}{2}}\eps)^2 + (|D|^{\alpha}\eps)^2)\Phi_1^2 - CA^{\alpha}t^{-\frac{3(3\alpha+5)}{2(\alpha+3)}}.
\end{align*}
By taking $A>A_1$ large enough, $T_0$ large enough, we deduce that:
$$
\frac{d}{dt}F(t)\geq -CA^{\alpha}t^{-\frac{3(3\alpha+5)}{2(\alpha+3)}}
$$
However, the choice of $A$ is independent of parameters. We set $A>\max(A_1,A_2)$, with $A_2$ defined in Claim \ref{claim:coercivite} for the coercivity of the localized  linearized operator. For now, $A$ is a constant. Then, integrating in time from $t$ to $S_n$ we conclude that: 
\begin{align*}
    F(t)\leq Ct^{-\frac{7\alpha+9}{2(\alpha+3)}},
\end{align*}
with the constant $C$ independent of the different parameters.
\subsection{Topological argument}\label{sec:topological_argument}

We argue by contradiction. Let suppose for all $z^{in}_n$ in \eqref{CI:z_n_in}, we have $t^*(z^{in}_n)>T_0$.

 Suppose first that one of the bootstrap estimates \eqref{boot:eps}, \eqref{boot:mu}, \eqref{boot:z_bar} or \eqref{boot:mu_bar} is saturated, in the sense that the equality is achieved. 
 
 \textit{1) Closing bootstrap for $\eps$.} First we start to show we can improve \eqref{boot:eps}. We recall that the notations $\phi$, $\phi_1$ and $\phi_2$ holds respectively for $\phi_A$, $\phi_{1,A}$ and $\phi_{2,A}$. Using the Cauchy-Schwarz inequality, \eqref{est:S}, \eqref{boot:eps} and the definition of $\phi_1$ , we get that:
\begin{align}
    F(t)\geq& -C t^{\frac{3(3\alpha+5)}{2(\alpha+3)}} + \frac{1}{2}\int\left( \eps|D|^{\alpha}\eps + \eps^2 - 3\widetilde{R}^2_1\eps^2\right)\frac{1-\phi}{(1+\mu_1)^2} +\frac{1}{2}\int\left( \eps|D|^{\alpha}\eps + \eps^2 - 3\widetilde{R}^2_2\eps^2 \right)\frac{\phi}{(1+\mu_2)^2}\notag\\
    &+ \int \frac{\eps^2}{2} \phi_2 +\int \left( \frac{V^4}{4} + V^{3}\eps - \frac{\left(V+\eps \right)^4}{4} \right)\phi_1 + \frac{3}{2}\widetilde{R}^2_1\eps^2\frac{1-\phi}{(1+\mu_1)^2} + \frac{3}{2}\widetilde{R}^2_2\eps^2\frac{\phi}{(1+\mu_2)^2} .  \label{coercoercivite:F}
\end{align}
First of all, we estimate the last term on the right hand side. We get that, by straight forward computations:
\begin{align*}
     \frac{V^4}{4} + V^{3}\eps - \frac{\left(V+\eps \right)^4}{4}= -\frac{3}{2}V^{2}\eps^2 - \eps^3V - \frac{\eps^4}{4}.
\end{align*}
Using the Sobolev embedding and  the bootstrap estimates on $\eps$ \eqref{boot:eps}, we deduce that: 
\begin{align*}
    \left\vert\int (\eps^3V + \frac{1}{4}\eps^4)\phi_1\right\vert\leq C t^{-\frac{3(3\alpha+5)}{2(\alpha+3)}}.
\end{align*}
Moreover, we have that: 
\begin{align*}
\MoveEqLeft
      \widetilde{R}^2_1\frac{1-\phi}{(1+\mu_1)^2} + \widetilde{R}^2_2\frac{\phi}{(1+\mu_2)^2} -V^{2}\phi_1= \left( \widetilde{R}^2_1 - R^2_1 \right) \phi_1 +\left( \widetilde{R}^2_2 - R^2_2 \right) \phi_1 - \widetilde{R}^2_1\frac{\phi}{(1+\mu_2)^2}- \widetilde{R}^2_2\frac{1-\phi}{(1+\mu_1)^2} \\
      &+2R_1R_2\phi_1-2(-R_1+R_2)(-P_1+P_2+bW)\phi_1 - (-P_1+P_2+bW)^2\phi_1
\end{align*}
Therefore, by applying the bootstrap estimate on $\eps$ \eqref{boot:eps}, the estimate on the profile $P_i$ \eqref{derivee:d_yP_i}, the estimate on the solitary waves \eqref{eq:partial_y_Ri}, the estimate on $\Lambda Q$ \eqref{eq:lambda_Q_DL} and finally the bootstrap estimate on $z$ \eqref{boot:z}, we get that:
\begin{align*}
\left\vert \int \left( \frac{V^4}{4} + V^{3}\eps - \frac{\left(V+\eps \right)^4}{4} \right)\phi_1 + \frac{3}{2}\widetilde{R}^2_1\eps^2\frac{1-\phi}{(1+\mu_1)^2} + \frac{3}{2}\widetilde{R}^2_2\eps^2\frac{\phi}{(1+\mu_2)^2} \right\vert \leq C t^{-\frac{4\alpha+6}{\alpha+3}} .
\end{align*}

Moreover, from the bootstrap estimates on $\mu$ \eqref{boot:mu} and $\bar{\mu}$ \eqref{boot:mu_bar} we have that: 
\begin{align*}
\left\vert\int\eps^2\phi_2\right\vert\leq Ct^{-\frac{\alpha+1}{\alpha+3}}\|\eps\|_{L^2}^2.
\end{align*}

Now, we estimate the two first integrals in \eqref{coercoercivite:F}. We claim the following:
\begin{align}
   \int\left( \eps|D|^{\alpha}\eps + \eps^2 - 3\widetilde{R}^2_i\eps^2\right)\frac{1-\phi}{(1+\mu_1)^2} +\left( \eps|D|^{\alpha}\eps + \eps^2 - 3\widetilde{R}^2_i\eps^2\right)\frac{\phi}{(1+\mu_2)^2} \geq \kappa \|\eps\|^2_{H^{\frac{\alpha}{2}}}, \quad i=1,2.
\end{align}

The proof of this inequality is given in Claim \ref{claim:coercivite} in the Appendix \ref{proof:claim_coer}. The proof is based on the coercivity of the linearized operator $L$. By combining the former inequalities and using Theorem \ref{monotonicity}, we deduce that: 
\begin{align*}
   \kappa\|\eps\|_{H^{\frac{\alpha}{2}}}^2 - Ct^{-\frac{4\alpha+6}{\alpha+3}} -Ct^{-\frac{\alpha+1}{\alpha+3}}\|\eps\|_{L^2}^2 \leq  F(t)\leq C t^{-\frac{7\alpha+9}{2(\alpha+3)}} .
\end{align*}
Therefore for $T_0$ large enough, we conclude that: 
\begin{align*}
    \|\eps\|_{H^{\frac{\alpha}{2}}}^2\leq C t^{-\frac{7\alpha+9}{2(\alpha+3)}}.
\end{align*}
Therefore, we strictly improved the bound \eqref{boot:eps} on $\eps$. This concludes the proof for $\eps$. 

\textit{2) Closing bootstrap for $\mu$,$\bar{\mu}$ and $\bar{z}$.}
Now, we improve the bound on $\mu$ \eqref{boot:mu}. We recall $\mu=\mu_1-\mu_2$ and $z=z_1-z_2$. Combining the bootstrap estimate on $\eps$ \eqref{boot:eps} and  $z$ \eqref{boot:z} on the right hand side of the estimate of $\dot{\mu_i}$ in \eqref{eq:mupt1_2} we deduce that:
\begin{align*}
    \left\vert\dot{\mu} - \frac{2b_1}{z^{\alpha+2}}\right\vert \leq Ct^{-\frac{3\alpha+5}{\alpha+3}}.
\end{align*}
Because $b_1<0$ and by the equivalent of $z$ in \eqref{boot:z2}, we have $\dot{\mu}<0$. By the initial condition $\mu(S_n)>0$, see \eqref{eq:mu_n_in}, $\mu$ is positive on $(t^*, S_n]$.

Then, multiplying by $\mu$, using the estimate on $\dot{z_i}$ \eqref{eq:zpt1_2} and the bootstrap on $z$ \eqref{boot:z} and $\mu$ \eqref{boot:mu}, we obtain that: 
\begin{align*}
    \left\vert\frac{\dot{\overbrace{\mu^2}}}{2} + \frac{2b_1}{\alpha+1}\dot{\overbrace{\frac{1}{z^{\alpha+1}}}}\right\vert \leq Ct^{-\frac{4\alpha+6}{\alpha+3}}.
\end{align*}
By the choice of the initial data, we have that:
\begin{align*}
\mu^2(S_n)=-\frac{4b_1}{\alpha+1}\frac{1}{z^{\alpha+1}(S_n)}.    
\end{align*}
Therefore, by integrating from $t$ to $S_n$, we get that:
\begin{align}
    \left\vert\frac{\mu^2}{2} + \frac{2b_1}{\alpha+1}\frac{1}{z^{\alpha+1}}\right\vert \leq Ct^{-\frac{3(\alpha+1)}{\alpha+3}}.\label{eq:mu^2}
\end{align}
 With the bootstrap hypothesis on $z$ \eqref{boot:z}, we deduce that: 
\begin{align*}
    \left\vert \mu -\sqrt{\frac{-4b_1}{\alpha+1}}\frac{t^{-\frac{\alpha+1}{\alpha+3}}}{a^{\frac{\alpha+1}{2}}} \right\vert\leq C_1 t^{-\frac{5\alpha+11}{4(\alpha+3)}},
\end{align*}
with the constant $C_1>0$.

Let us compute the bound on $\bar{\mu}$. From the estimate on $\dot{\mu_i}$ \eqref{eq:mupt1_2} and the bootstrap estimate on $\eps$ \eqref{boot:eps} and $z$ \eqref{boot:z}, we obtain that:
\begin{align*}
    |\dot{\bar{\mu}}|\leq C t^{-\frac{3\alpha+5}{\alpha+3}}.
\end{align*}
By the choice of the initial data, we have that $\mu_1(S_n)=-\mu_2(S_n)$. Thus, by integrating we deduce that: 
\begin{align}
    |\bar{\mu}|\leq C_2 t^{-\frac{2(\alpha+1)}{\alpha+3}} \label{opti:barmu}
\end{align}
with the constant $C_2>0$.

Let us get a bound on $\bar{z}$. Using the fact that $|\beta(\Gamma)|+|\delta(\Gamma)|\leq\frac{2(\beta_0+\delta_0)}{z^{\alpha+1}}$, the bound obtain for $\bar{\mu}$ \eqref{opti:barmu} and the estimate on $\dot{z_i}$ \eqref{eq:zpt1_2}, we deduce that: 
\begin{align*}
    |\dot{\bar{z}}|\leq& |\dot{\bar{z}}-\bar{\mu}+\beta(\Gamma) + \delta(\Gamma)| + |\bar{\mu}| + |\beta(\Gamma) + \delta(\Gamma)| 
    & leq C_3 t^{-\frac{3\alpha+5}{2(\alpha+3)}} + (2(\beta_0+\delta_0) + C_2)t^{-\frac{2(\alpha+1)}{\alpha+3}} \\
    & \leq 2C_3t^{-\frac{3\alpha+5}{2(\alpha+3)}}.
\end{align*}
Therefore by integrating, we conclude that: 
\begin{align*}
    |\bar{z}|\leq \frac{2C_3(2(\alpha+3))}{\alpha-1} t^{-\frac{\alpha-1}{2(\alpha+3)}}
\end{align*}

Hence, by taking the constant $C^{*}>\displaystyle\max\left(C_1,C_2,\frac{2C_3(2(\alpha+3))}{\alpha-1}\right)$, we can close the bootstrap estimate on $\mu$, $\bar{\mu}$ and $\bar{z}$. Then, none of the previous inequalities on $\dot{\mu}$, $\dot{\bar{\mu}}$ and $\dot{\bar{z}}$ can saturate independently of the initial condition $z_n^{in}$.

\textit{3) Closing bootstrap for $z$.}
Subsequently, the inequality \eqref{boot:z} saturates for any $z_n^{in}$. We now prove that this equality is the source of a contradiction on $t^*(z_n^{in})$.

First, we remark $z_n^{in}=\left(a^{\frac{\alpha+3}{2}}S_n + \lambda_n S_n^{\frac{1}{2} +r}\right)^{\frac{2}{\alpha+3}}$, for some $\lambda_n\in [-1,1]$. Therefore, we can write $t^*(z^{in}_n)=t^{*}(\lambda_n)$. We set: 
\begin{align}
    \Phi: [-1,1]&\longrightarrow\{-1,1\}\\
            \lambda &\longmapsto \left(z^{\frac{\alpha+3}{2}}(t^{*}(\lambda))-a^{\frac{\alpha+3}{2}}t^{*}(\lambda) \right)\left(t^{*}(\lambda)\right)^{-\frac{1}{2}-r},\notag
\end{align}
and 
\begin{align}
    f:\R&\longrightarrow\R^{+}\\
        s&\longmapsto \left(z^{\frac{\alpha+3}{2}}(s)-a^{\frac{\alpha+3}{2}}s \right)^2s^{-1-2r}.\notag
\end{align}
By assumption, we have for any $\lambda\in [-1,1]$, $t^{*}(\lambda)>T_0$ and thus: 
\begin{align}
    |z^{\frac{\alpha+3}{2}}(t^{*}(\lambda))-a^{\frac{\alpha+3}{2}}t^{*}(\lambda)| = (t^{*}(\lambda))^{\frac{1}{2}+r}.\label{sature}
\end{align}
We claim:
\begin{claim}
\begin{enumerate}
    \item Transversality condition: Let $s_0>T_0$ such that \eqref{sature} is verified at $s_0$, then: 
    \begin{align}\label{transversality}
        f \text{ is decreasing on a neighbourhood of } s_0.
    \end{align}
    \item Continuity: $\Phi\in C^{0}([-1,1]:\{-1,1\})$.
\end{enumerate}
\end{claim}
Let us assume the claim and finish the proof. The transversality condition \eqref{transversality} implies that $t^{*}(\pm 1)=S_n$. Moreover, $\Phi(\pm1)=\pm1$. This contradicts (2) of the former claim. Now, we prove the claim. First, we prove the transversality condition \eqref{transversality}. By direct computations, we have that: 
\begin{align*}
    f'(s) = 2\left(\dot{\wideparen {z^{\frac{\alpha+3}{2}}}}(s)-a^{\frac{\alpha+3}{2}} \right) \left(z^{\frac{\alpha+3}{2}}(s)-a^{\frac{\alpha+3}{2}}s \right)  s^{-1-2r} -\left(1+2r \right)\left(z^{\frac{\alpha+3}{2}}(s)-a^{\frac{\alpha+3}{2}}s \right)^2s^{-2-2r}.
\end{align*}
From the estimate obtain on $\mu^2$ \eqref{eq:mu^2} and the estimate on $\dot{z_i}$ \eqref{eq:zpt1_2}, we obtain that: 
\begin{align}\label{equation_modulation}
    \left\vert\dot{\wideparen {z^{\frac{\alpha+3}{2}}}}(t)-\frac{\alpha+3}{2}\sqrt{\frac{-4b_1}{\alpha+1}}\right\vert\leq Ct^{-\frac{\alpha+1}{\alpha+3}}.
\end{align}
Therefore, by using \eqref{sature} and \eqref{equation_modulation}, and since $a^{\frac{\alpha+3}{2}}=\frac{\alpha+3}{2}\sqrt{\frac{-4b_1}{\alpha+1}}$,  we deduce that: 
\begin{align*}
    f'(s_0)< Cs_0^{-1-3r} -\left(1+2r \right)s_0^{-1}.
\end{align*}
Since $r>0$ and for $T_0$ large enough, we conclude that: 
\begin{align*}
    f'(s_0)<0.
\end{align*}

To prove the second part of the former claim, it is enough to show that $\lambda\mapsto t^{*}(\lambda)$ is continuous. Let us fix $\lambda\in [-1,1]$. From the transversality condition, there exists $\eps_\lambda>0$ such that $\forall\eps\in (0,\eps_\lambda)$, $\exists \delta>0$ and the two following conditions are verified: $f(t^{*}(\lambda)-\eps)>1+\delta$, and for all $t\in [t^{*}(\lambda)+\eps,S_n]$ (possibly empty), $f(t)< 1-\delta$.

Note that the function is well defined, since the function $z$ is globally well defined, see Remark \eqref{global_well_modul}. 
Then by the continuity of the flow, there exists $\eta>0$ such that for all $|\lambda-\bar{\lambda}|<\eta$, with $\bar{\lambda}\in[-1,1]$, the corresponding $\bar{f}$ verifies $|\bar{f}(s)-f(s)|<\frac{\delta}{2}$ for $s\in [t^{*}(\lambda)-\eps,S_n]$. Therefore, we obtain that for all $s\in[t^{*}(\lambda)+\eps,S_n]$:
\begin{align*}
    \bar{f}(s)< \lvert\bar{f}(s) - f(s) \rvert + f(s)< 1-\frac{\delta}{2}.
\end{align*}
Thus, $t^{*}(\bar{\lambda})<t^{*}(\lambda)+\eps$. Furthermore, 
\begin{align*}
    \bar{f}(t^{*}(\lambda)-\eps)> f(t^{*}(\lambda)-\eps) - \lvert\bar{f}(t^{*}(\lambda)-\eps) - f(t^{*}(\lambda)-\eps) \rvert > 1+\frac{\delta}{2}.
\end{align*}
In other words, $t^{*}(\lambda)-\eps<t^{*}(\bar{\lambda})$, and then $\Phi$ is continuous.

This contradicts the fact $t^{*}(\lambda)>T_0$ and implies the existence of  $z_{n}^{in}$ such that \eqref{boot:eps}-\eqref{boot:mu_bar} are true for all $t\in[T_0,S_n]$.

\subsection{Conclusion}\label{sec:conclusion}

In this section we have proved that there exists $(z_n^{in})^{\frac{\alpha+3}{2}}\in [a^{\frac{\alpha+3}{2}}S_n^{\frac{1}{2}+r} - S_n,a^{\frac{\alpha+3}{2}}S_n + S_n^{\frac{1}{2}+r}]$ such that the bootstrap estimates \eqref{boot:eps}-\eqref{boot:mu_bar} are true for all $t\in[T_0,S_n]$. Let us show this implies Theorem \ref{main_theo}. From \eqref{boot:eps}, we obtain that: 
\begin{align*}
    \|v_n(T_0,\cdot)\|_{H^{\frac{\alpha}{2}}}\leq \|\eps_n(T_0,\cdot)\|_{H^{\frac{\alpha}{2}}} + \|V(\Gamma_n(T_0),\cdot)\|_{H^{\frac{\alpha}{2}}}\leq C.
\end{align*}
Therefore, by Banach-Alaoglu, there exists $w_0\in H^{\frac{\alpha}{2}}(\R)$ and a sub-sequence also denoted by $(v_n)_n$ such that: 
\begin{align*}
    v_n(T_0)\rightharpoonup w_0.
\end{align*}
Thus, we denote by $w$ the solution of \eqref{mBO} such that $w(T_0)=w_0$.
Let $t>T_0$. From the weak continuity of the flow of Theorem \ref{weak_continuity}, we have that: 
\begin{align*}
    \MoveEqLeft
    \left\|w(t,\cdot) + Q\left( \cdot - \frac{a}{2}t^{\frac{2}{\alpha+3}} \right) - Q\left( \cdot + \frac{a}{2}t^{\frac{2}{\alpha+3}} \right)\right\|_{H^{\frac{\alpha}{2}}}\\
    & \leq \liminf_{n\to\infty}  \|\eps_n(t,\cdot) \|_{H^{\frac{\alpha}{2}}} +\liminf_{n\to\infty} \left\|V(\Gamma_n(t),\cdot) + Q\left( \cdot - \frac{a}{2}t^{\frac{2}{\alpha+3}} \right) - Q\left( \cdot + \frac{a}{2}t^{\frac{2}{\alpha+3}} \right)\right\|_{H^{\frac{\alpha}{2}}}.
\end{align*}
Then, by using \eqref{boot:eps}-\eqref{boot:z_bar}, we conclude that: 
\begin{align*}
    \left\|w(t,\cdot) + Q\left( \cdot - \frac{a}{2}t^{\frac{2}{\alpha+3}} \right) - Q\left( \cdot + \frac{a}{2}t^{\frac{2}{\alpha+3}} \right)\right\|_{H^{\frac{\alpha}{2}}}\leq C t^{-\frac{\alpha-1}{4(\alpha+3)}}.
\end{align*}


\textbf{Acknowledgements} The authors thank Razvan Mosincat for interesting discussions, and Didier Pilod for suggesting this problem and constant support. The authors were supported by a Trond Mohn foundation grant.


\appendix

\section{Local well-posedness}\label{appendix:LWP}

We recall the results of well-posedness of \eqref{mBO}.

\begin{theo}[\cite{Guo12}, Theorem 1.5]
    Let $\alpha \in (1,2)$, and $u_0\in H^{s}(\mathbb{R})$, with $s \geq \frac{1}{2}-\frac{\alpha}{4}$. There exists a time $T=T(\|u_0\|_{H^{\frac{1}{2}-\frac{\alpha}{4}}(\mathbb{R})})>0$, and a unique solution $u\in \mathcal{C}([-T,T],H^s(\mathbb{R}))$ of \eqref{mBO}. Furthermore, the flow $u_0 \mapsto u$ is locally Lipschitz continuous from $H^s(\mathbb{R})$ to $\mathcal{C}([-T,T], H^s(\mathbb{R}))$. 
\end{theo}

Because the equation is subcritical, we obtain as a corollary the global well-posedness.

\begin{coro}[\cite{Guo12}, Corollary 1.6]\label{solution:global_well_posed}
For any initial condition $u_0\in H^\frac{\alpha}{2}(\mathbb{R})$, there exists a unique global solution of \ref{mBO} in $\mathcal{C}(\mathbb{R},H^\frac{\alpha}{2}(\mathbb{R}))$.
\end{coro}

We continue with another property of the flow, which is the weak-continuity in $H^{\frac{\alpha}{2}}(\mathbb{R})$.

\begin{theo}[Weak continuity of the flow]\label{weak_continuity}
Let $\alpha\in(1,2)$. Suppose that $u_{0,n}\rightharpoonup u_0 \in H^{\frac{\alpha}{2}}(\R)$. We consider $u_n$ solutions of \eqref{mBO} corresponding to the initial data $u_n(0)=u_{n,0}$ and satisfying $u_n\in C([0,T]:H^{\frac{\alpha}{2}}(\R))$ for any $T>0$. Then, $u_n(t)\rightharpoonup u(t)$ in $H^{\frac{\alpha}{2}}(\R)$, for all $t\geq 0$. 
\end{theo}

The proof of the weak continuity of the flow relies on the well-posedness result given in the Corollary \ref{solution:global_well_posed}. We refer to \cite{eychenne2021asymptotic} Appendix A, \cite{goubet2009weak} for a proof of this result.

\section{Justification of the definition of $S_0$}\label{proof:lemma:S0}

First, we recall some well-known results on pseudo-differential operators (see \cite{alinhac2012operateurs}, or \cite{hormander2007analysis} chapter $18$).
Let $D=-i\partial_x$. We define the symbolic class $\mathcal{S}^{m,q}$ by 
\begin{align*}
	\mathcal{S}^{m,q} := \left\{ 
		a\in C^{\infty}(\R_x \times \mathbb{R}_\xi) ; \quad 
		\forall k,\beta\in \N,\exists C_{k,\beta}>0 \text{ such that } |\partial_x^k\partial_{\xi}^{\beta}a(x,\xi)|\leq C_{k,\beta}\langle x \rangle^{q-k} \langle \xi\rangle^{m-\beta} 
	\right\}.
\end{align*}
For all $u$ in the Schwartz space $\mathcal{S}(\R)$, we set the operator associated to the symbol $a(x,\xi)\in \mathcal{S}^{m,q}$ by 
$$
a(x,D)u:=\displaystyle\frac{1}{2\pi}\int e^{ix\xi} a(x,\xi) \mathcal{F}(u)(\xi) d\xi.
$$ 
We state the three following results 
\begin{enumerate}
	\item  Let $a\in \mathcal{S}^{m,q}$, there exists $C>0$, such that for all $u\in \mathcal{S}(\R)$ 
	\begin{align}
		\| a(x,D)u \|_{L^2}\leq C\|\langle x \rangle^q\langle D \rangle^m u \|_{L^2}. \label{pseudo:est}
	\end{align}
	\item Let $a\in \mathcal{S}^{m,q}$ and $b\in \mathcal{S}^{m',q'}$, then there exists $c\in \mathcal{S}^{m+m',q+q'}$ such that 
	\begin{align}
		a(x,D)b(x,D)=c(x,D).\label{pseudo:prod}
	\end{align}
	\item If $a\in \mathcal{S}^{m,q}$ and $b\in \mathcal{S}^{m',q'}$ are two operators, we define the commutator by $[a(xD),b(x,D)]:=a(x,D)b(x,D)-b(x,D)a(x,D)$.
	Moreover there exists $c\in \mathcal{S}^{m+m'-1,q+q'-1}$ such that 
	\begin{align}
		[a(x,D),b(x,D)]=c(x,D).\label{pseudo:comp}
	\end{align}
	\item Let $a\in \mathcal{S}^{m,q}$, we have the following development for the adjoint $a^*$ of $a$. Let $k\in\N$, then 
	\begin{align}\label{dev:adjoint}
	    a^*(x,\xi)=\sum_{\beta\leq k} \frac{1}{\beta!}\partial_{\xi}^{\beta}D_{x}^{\beta}\bar{a}(x,\xi) + R_k(x,\xi)
	\end{align}
	with $\partial_{\xi}^{\beta}D_{x}^{\beta}\bar{a}\in \mathcal{S}^{m-\beta,q-\beta}$ and $R_k\in \mathcal{S}^{m-\beta-1,q-\beta-1}$. Moreover the rest $R_k$ is given by 
	\begin{align}\label{rest:pseudodiff}
	    R_k(x,\xi)=\frac{1}{2\pi}\int_0^1(1-t)^{2k+1}dt\int e^{-iy\eta}\sum_{\beta+\gamma =2k+2}\frac{2k+2}{\beta!\gamma!}\partial_y^{\beta}\partial_{\eta}^{\gamma}\bar{a}(x-ty,\xi-t\eta)y^{\beta}\eta^{\beta}dyd\eta.
	\end{align}
\end{enumerate}
As a consequence of \eqref{pseudo:prod},  $\langle D \rangle^m\langle x \rangle^q\langle D \rangle^{-m}\in \mathcal{S}^{0,q}$. Therefore, by \eqref{pseudo:est}, we have 
\begin{align*}
	\|\langle D \rangle^m\langle x \rangle^qu \|_{L^2}&=\|\langle D \rangle^m\langle x \rangle^q\langle D \rangle^{-m}\langle D \rangle^mu \|_{L^2}\\
	&\leq C_2\|\langle x \rangle^q\langle D \rangle^m u \|_{L^2},
\end{align*}
for $C_2>0$.
By the same computations with $\langle x \rangle^q$ instead of $\langle D \rangle^{m}$, there exists $C_1>0$ such that
\begin{align*}
	C_1\|\langle x \rangle^q\langle D \rangle^m u \|_{L^2}&\leq  \|\langle D \rangle^m\langle x \rangle^qu \|_{L^2}.
\end{align*}
Gathering these two estimates, we conclude that
\begin{align}
	C_1\|\langle x \rangle^q\langle D \rangle^m u \|_{L^2}&\leq  \|\langle D \rangle^m\langle x \rangle^qu \|_{L^2}\leq C_2\|\langle x \rangle^q\langle D \rangle^m u \|_{L^2}.  \label{pseudo:equiv}
\end{align}
We recall also the Schur's test.
\begin{theo}[Schur's test \cite{halmos2012bounded}, Theorem 5.2]\label{theo:schur}
	Let $p, q$ be two non-negative measurable functions. If there exists $\alpha,\beta>0$ such that 
	\begin{enumerate}
		\item $\displaystyle\int_{\R}|K(x,y)| q(y) dy \leq \alpha p(x) \text{ a.e. } x\in \R $.
		\item $\displaystyle\int_{\R}|K(x,y)| p(x) dx \leq \beta q(y) \text{ a.e. } y\in \R$.
	\end{enumerate}
	Then $Tf:=\displaystyle\int_{\R}K(x,y)f(y) dy$ is a bounded operator on $L^2(\R)$.
\end{theo}

We recall two other lemmas useful for the rest of the appendix. The definition of $\phi$ is given in \eqref{defi:phi}.
\begin{lemm}[\cite{kenig2011local} Claim 5]\label{lemma:poid}
	There exists $C>0$ such that
	\begin{align*}
	|\phi(x) - \phi(y)| &\leq C\frac{|x-y|}{\left( \langle x \rangle  \langle y \rangle \right)^{\frac{\alpha+1}{2}}} + C\frac{|x-y|^2}{\left( \langle x \rangle + \langle y \rangle \right)^{\alpha+2}} \quad  \text{ if } \quad  |x-y|\leq \frac{1}{2}\left( \langle x \rangle + \langle y \rangle \right),\\
	|\phi(x) - \phi(y)| &\leq  C \quad  \text{ if } \quad  |x-y|\geq \frac{1}{2}\left( \langle x \rangle + \langle y \rangle \right).
	\end{align*}
\end{lemm}

\begin{lemm}[\cite{kenig2011local}, Lemma A.2]\label{lemm:noyaubassefreq}
	Let $p$ be a homogeneous function of degree $\beta>-1$. Let $\chi\in C^{\infty}_0(\R)$ such that $0\leq \chi\leq1$, $\chi(\xi)=1$ if $|\xi|<1$ and $\chi(\xi)=0$ if $|\xi|>2$. Let 
	$$
	k(x)= \frac{1}{2\pi}\int e^{ix\xi}p(\xi)\chi(\xi) d\xi.
	$$
	Then for all $q\in\mathbb{N}$, there exists $C_q>0$ such that, for all $x\in\R$,
	\begin{align}
		|\partial_x^q k(x)|\leq \frac{C_q}{\langle x \rangle^{\beta+q+1}}.
	\end{align}
\end{lemm}

Now, we can start the proof of the justification of the definition of $S_0$.
\begin{proof}
We recall the definition of $\Lambda Q$, and estimate on $Q$ from \cite{FLS16}:
\begin{align*}
\Lambda Q = \frac{\alpha}{2(\alpha+1)} Q + \frac{1}{\alpha +1} x \partial_x Q, \quad \vert Q \vert + \vert x \partial_x Q \vert \leq \frac{1}{1+ \vert x \vert^{1+\alpha}}.
\end{align*}

Since $\Lambda Q \in L^2(\R)$, we can define by the Fourier transform $(1+ \vert D \vert^{\alpha})^{-1} \Lambda Q \in H^{\alpha}(\R)$:
\begin{align*}
\left\| (1+ \vert D \vert^\alpha)^{-1} \Lambda Q \right\|_{H^\alpha}^2 = \left\| \frac{(1+ \vert \xi \vert^2)^{\frac{\alpha}{2}}}{1+\vert \xi \vert^\alpha} \widehat{\Lambda Q} \right\|_{L^2}^2 \lesssim \| \Lambda Q \|_{L^2} <\infty.
\end{align*}
The integral of $\left( 1+ \vert D \vert^\alpha \right)^{-1} \Lambda Q$ on a finite interval is well-defined since it is in $L^2(\R)$. However, it is not clear that the integral over an infinite interval is finite. We use the pseudo-differential theory to prove that the limit is finite. Let us define $\chi$, a cut-off function equal to $1$ in a neighbourhood of $0$, with compact support. Let $I$ be a compact interval. By the Cauchy-Schwarz inequality :
\begin{align*}
\MoveEqLeft
\int_I \left\vert \left( 1+ \vert D \vert^\alpha \right)^{-1} \Lambda Q \right\vert \leq \int_I \left\vert \left( 1-\chi(D) \right) \left( 1+ \vert D \vert^\alpha \right)^{-1} \Lambda Q \right\vert  + \int_I \left\vert \chi(D) \left( 1+ \vert D \vert^\alpha \right)^{-1} \Lambda Q \right\vert \\
	& \leq C \left\| \langle x \rangle^{\frac{3}{4}
} \left( 1-\chi(D) \right) \left( 1+ \vert D \vert^\alpha \right)^{-1} \Lambda Q \right\|_{L^2(I)} + C \left\| \langle x \rangle^{\frac{3}{4} } \chi(D) \left( 1+ \vert D \vert^\alpha \right)^{-1} \Lambda Q \right\|_{L^2(I)} \\ & = \mathcal{I}_1+ \mathcal{I}_2.
\end{align*}
Note that the previous constant can be chosen independently of $I$. We have from \eqref{pseudo:comp} that the symbol $\langle x \rangle^{\frac{3}{4}} (1-\chi(\xi))(1+\vert \xi \vert^\alpha)^{-1}$ belongs to $\mathcal{S}^{ -\alpha,\frac{3}{4}}\subset \mathcal{S}^{0,\frac{3}{4}} $. Thus, since $\langle x\rangle^{\frac{3}{4}} \Lambda Q\in L^2(\R)$:
\begin{align*}
\mathcal{I}_1 \lesssim \| \langle x \rangle ^{\frac{3}{4}} \Lambda Q \|_{L^2(\mathbb{R})} <\infty.
\end{align*}

We can not deal with the integral $\mathcal{I}_2$ with symbols only, because $\chi(\xi) (1+\vert \xi \vert^\alpha)^{-1}$ is not smooth around $0$. We use the commutator to bring the decay in $x$ close to $\Lambda Q$ (notice the integral is over $\mathbb{R}$):
\begin{align*}
\mathcal{I}_2^2 \lesssim \int_\mathbb{R} \left( \left[ \langle x \rangle^{\frac{3}{4}}, \chi(D)(1+\vert D \vert^\alpha )^{-1} \right] \Lambda Q \right)^2 + \int_\mathbb{R} \left( \chi(D)(1+\vert D \vert^\alpha )^{-1}\langle x \rangle^{\frac{3}{4}} \Lambda Q \right)^2.
\end{align*}
By the Plancherel formula, the second term can be bounded by $\left\| \langle x \rangle^{ \frac{3}{4}} \Lambda Q \right\|_{L^2}^2 <\infty$. The first term needs to develop the commutator. First, let us define the kernel $k$ satisfying:
\begin{align*}
\chi(D) (1+ \vert D \vert^\alpha)^{-1} u (x)=\frac{1}{2\pi} \int e^{i\xi x} \frac{\chi(\xi)}{1+ \vert \xi \vert^\alpha} \hat{u}(\xi) d\xi = k\star u(x), \quad \text{ so } \quad \hat{k}(\xi)= \frac{\chi(\xi)}{1+ \vert \xi \vert^\alpha}.
\end{align*}
The kernel $k$ is well-defined as the inverse Fourier transform of a function in $L^2$. We thus get:
\begin{align*}
\left[ \langle x \rangle^{\frac{3}{4}}, \chi(D) (1+ \vert D \vert^\alpha)^{-1} \right]u 
	& = \langle x \rangle^{\frac{3}{4}} k \star u (x) - k \star \left( \langle x \rangle^{\frac{3}{4}} u \right) (x) \\
	& = \int k(x-y) \left( \langle x \rangle^{\frac{3}{4}} - \langle y \rangle^{\frac{3}{4}} \right) u(y) dy.
\end{align*}
By Lemma \ref{theo:schur} and the symmetry of $k$, it is enough to prove that $y\mapsto k(x-y)\left(\langle x \rangle^{\frac{3}{4}} -\langle y \rangle^{\frac{3}{4}} \right)\in L^{1}(\R)$. First, we have to estimate $k$. By  integrating by parts twice, we deduce that :
\begin{align}\label{estim:k}
\frac{1}{1+x^2}(1-\partial_\xi^2) e^{ix \xi}=e^{ix\xi} \quad \text{ and } \quad \left\vert k(x) \right\vert = \left\vert \frac{1}{2\pi} \int e^{ix\xi} \frac{\chi(\xi)}{1+\vert \xi \vert^\alpha} d\xi \right\vert \leq \frac{C_\alpha}{\langle x \rangle^2}.
\end{align}
Let $A_1:=\{y\in \R:|x-y|\leq \frac{1}{2}\left(\langle x \rangle + \langle y \rangle \right) \}$,  and $A_2:=\{y\in \R : |x-y|> \frac{1}{2}\left(\langle x \rangle + \langle y \rangle \right) \}$. Notice the following equivalences :
\begin{align}\label{estim:sim1}
    |x-y|\leq \frac{1}{2}\left(\langle x \rangle + \langle y \rangle \right) \Rightarrow \langle x \rangle\sim \langle y \rangle,
\end{align}
and
\begin{align}\label{estim:sim2}
   |x-y|> \frac{1}{2}\left(\langle x \rangle + \langle y \rangle \right) \Rightarrow \langle x-y \rangle\sim |x-y| \sim  \langle x \rangle + \langle y \rangle.
\end{align}
Then, from \eqref{estim:k} and \eqref{estim:sim2}, we deduce that
\begin{align}\label{estim:A1}
   \left| \int_{A_2} k(x-y)\left( \langle x \rangle^{\frac{3}{4}} - \langle y \rangle^{\frac{3}{4}} \right)  dy \right|\leq \int_{A_2}\frac{\langle x \rangle^{\frac{3}{4}} + \langle y \rangle^{\frac{3}{4}}}{\left(\langle x \rangle + \langle y \rangle\right)^{\frac{3}{4}}} \frac{1}{\langle x-y \rangle^{\frac{5}{4}}} dy\leq C.
\end{align}
Moreover by $\eqref{estim:sim1}$, we obtain that $A_1\subset[-c_2|x|,-c_1| x|]\cup[c_1| x|,c_2| x|] $, for some $0<c_1<c_2<+\infty$ independent of $x$. Moreover,
by the mean value theorem and $\ln(x+\langle x\rangle)'=\frac{1}{\langle x \rangle}$, we get that 
\begin{align}\label{estim:A2}
    \left| \int_{A_1} k(x-y)\left( \langle x \rangle^{\frac{3}{4}} - \langle y \rangle^{\frac{3}{4}} \right)  dy \right| \leq C\langle x \rangle^{-\frac{1}{4}}\int_{A_1} \frac{1}{\langle x-y \rangle} dy
    \leq C\langle x \rangle^{-\frac{1}{4}}\ln(C(|x| + \langle x\rangle))
    \leq C.
\end{align}
Gathering \eqref{estim:A1} and \eqref{estim:A2}, we conclude that $k$ defines a bounded operator on $L^{2}(\R)$. It implies that $\mathcal{I}_2$ is bounded, and thus $\int_I \vert (1+\vert D \vert^{\alpha}) \Lambda Q \vert$ is bounded independently of $I$. This achieves the proof of the well-posedness of $S_0$, and that $S_0$ has a finite limit at $-\infty$.
\end{proof}

\section{Proof of the preliminary results}\label{proof:preliminary}

\begin{proof}[Proof of Lemma \ref{est:commsimpl}]
Let $\chi$ be a smooth cut-off function supported around $0$. To estimate this commutator we split the norm in low and high frequency. For the low frequency we use the Schur's Lemma (Lemma \ref{theo:schur}), and the pseudo-differential calculus for the high frequency. To get an explicit dependence in $A$ we prove the estimate
\begin{align*}
    \bigg\|\left[ |D|^{\alpha} , \Phi \right]u \bigg\|_{L^2}^2\leq \begin{cases}C \int u^2\Phi^2, \quad \text{ if } \quad  \alpha \in ]0,1] \\
    C \int \left(u^2+ \left(|D|^{\frac{\alpha}{2}}u\right)^2\right)\Phi^2,\quad \text{ if }\quad  \alpha \in ]1,2] \end{cases}
\end{align*}
Then, we conclude Lemma \ref{est:commsimpl} by changing the variable $x=\frac{x'}{A}$ and multiplying by $\left\vert \frac{1}{(1+\mu_1)^2} - \frac{1}{(1+\mu_2)^2}  \right\vert$.

Let us start the proof. By the Schur's lemma (Lemma \ref{theo:schur}), we deduce that 
\begin{align*}
    \bigg\|\left[\chi(D))|D|^{\alpha} , \Phi \right]u \bigg\|_{L^2}^2\leq C\int u^2\Phi^2
\end{align*}
From pseudo-differential calculus, and $\langle x\rangle^{\frac{\alpha}{2}}\sim 1+x^{\frac{\alpha}{2}} $, we get that 
\begin{align*}
    \bigg\|\left[ (1-\chi(D))|D|^{\alpha} , \Phi \right]u \bigg\|_{L^2}^2 \leq\begin{cases}
    C\displaystyle\int u^2 \Phi^2, \quad \text{ if } \quad \alpha\in ]0,1]     \\
     C\displaystyle\int u^2\Phi^2 +C \int \left( |D|^{\frac{\alpha}{2}} (u \Phi)\right)^2,\quad  \text{ if } \quad  \alpha\in ]1,2]  \end{cases}.
\end{align*}
Again, by applying the pseudo-differential calculus, we deduce that 
\begin{align*}
\int \left( |D|^{\frac{\alpha}{2}} (u \Phi)\right)^2 &\leq C\left( \int \left( \chi(D)|D|^{\frac{\alpha}{2}} (u \Phi)\right)^2 +  \int \left( (1-\chi(D))|D|^{\frac{\alpha}{2}} (u \Phi)\right)^2\right)\\ 
&\leq C\left(\int u^2\Phi^2 + \int (|D|^{\frac{\alpha}{2}}u)^2\Phi^2\right).
\end{align*}
Then, by changing the variable $x=\frac{x'}{A}$ and multiplying by $\bigg|\frac{1}{(1+\mu_1)^2} - \frac{1}{(1+\mu_2)^2}\bigg|$, we conclude the proof of Lemma \ref{est:commsimpl}.
\end{proof}

\begin{proof}[Proof of Lemma \ref{esttc}]
By direct computations and Young's inequality, we have that  
\begin{align}
\MoveEqLeft
\left|\int |D|^{\alpha}\left(u\Phi_{j,A} \right) \left((|D|^{\alpha}u)\Phi_{j,A} \right)-\int \left(|D|^{\alpha} u\right)^2\Phi^2_{j,A}\right|  \notag \\
    & =\left\vert  \int |D|^{\alpha}u\Phi_{1,A}[|D|^{\alpha},\Phi_{1,A}]u \right\vert \leq \frac{C}{A^{\frac{\alpha}{2}}} \int \left(|D|^{\alpha}u\right)^{2}\Phi_{1,A}^2 +CA^{\frac{\alpha}{2}}\left\| [|D|^{\alpha},\Phi_{1,A}] u\right\|_{L^2}^2.  \label{eq:sans_nom}
\end{align}
and by the change of variable $x'=\frac{x}{A}$ and $v(x')=u(x)$:
\begin{align*}
    \left\| [|D|^{\alpha},\Phi_{1,A}] u\right\|_{L^2}^2 = \frac{1}{ A^{2\alpha-1}} \left\| [|D|^{\alpha},\Phi_{1}] v\right\|_{L^2}^2.
\end{align*}

We write
\begin{align*}
    \left\| [|D|^{\alpha},\Phi_{1}] v \right\|_{L^2}^2 \leq C\left( \left\| [|D|^{\alpha}\chi(D),\Phi_{1}]v \right\|_{L^2}^2 + \left\| [|D|^{\alpha}(1-\chi(D)),\Phi_{1}] v\right\|_{L^2}^2 \right).
\end{align*}
Using Theorem \ref{theo:schur}, we deduce that 
\begin{align*}
\left\| [|D|^{\alpha}\chi(D),\Phi_{1}] v \right\|_{L^2}^2 \leq C\int v^2 \Phi_{1}^2 .
\end{align*}
Moreover, using pseudo-differential calculus, we deduce that 
\begin{align*}
  \left\|[|D|^{\alpha}(1-\chi(D)),\Phi_{1}] v \right\|_{L^2}^2\leq C \int \left(v^2 + \left(|D|^{\frac{\alpha}{2}}v\right)^2 \right) \Phi_{1}^2.
\end{align*}
Gathering those estimates and coming back to the initial data, we get:
\begin{align*}
    \left\| [\vert D \vert^\alpha (1-\chi(D)), \Phi_{1,A}] u\right\|_{L^2}^2 \leq \frac{C}{A^{\alpha}} \int \left( u^2 + \left( \vert D \vert^{\frac{\alpha}{2}} u \right) \right)^2 \Phi_{1,A}.
\end{align*}
Using this last inequality in \eqref{eq:sans_nom}, we conclude the lemma.
\end{proof}

\begin{proof}[Proof of Lemma \ref{lemm:commH1}]
We recall that if $A,B$ are two pseudo-differential operators then the commutator $[A,B]$ is also a pseudo-differential $C$. Moreover the principal symbol of $C$ is given by 
\begin{align}\label{eq:Poisson_bracket_symbol}
\{a,b\}=\partial_{\xi} a\partial_y b - \partial_ya\partial_{\xi}b,
\end{align}
with $a,b$ respectively symbol of $A$ and $B$. Therefore,
$[(1-\chi(D))|D|^{\alpha},\phi_{1}]\in \mathcal{S}^{\alpha-1,-\alpha-1} \subset \mathcal{S}^{\frac{\alpha}{2},-\alpha-1}$. Then, by applying the pseudo-differential calculus and the fact $\partial_y \phi_{1}=\left(\frac{1}{(1+\mu_2)^2} - \frac{1}{(1+\mu_1)^2}\right)\partial_y\phi$, we have that  
\begin{align*}
    \left\|[(1-\chi(D))|D|^{\alpha},\phi_{1}]u \right\|_{L^2}\leq C \bigg|\frac{1}{(1+\mu_1)^2} - \frac{1}{(1+\mu_2)^2}\bigg|^{\frac{1}{2}}  \|u\Phi_{1} \|_{H^{\frac{\alpha}{2}}}.
\end{align*}
Now, we estimate the low frequency. Let $k$ be the operator defined by $\mathcal{F}(k(u))(\xi)=\chi(\xi)|\xi|^{\alpha}\mathcal{F}(u)(\xi)$. Then, we have that 
\begin{align*}
    [\chi(D)|D|^{\alpha},\phi_{1}]u=\left( \frac{1}{(1+\mu_2)^2} - \frac{1}{(1+\mu_1)^2}\right)\int k(x-y)(\phi(y)-\phi(x))u(y)dy.
\end{align*}
To prove that $[\chi(D)|D|^{\alpha},\phi_{1}]$ defines an operator bounded on $L^2(\R)$, we use the Schur's lemma (Lemma \ref{theo:schur}) on $\displaystyle x \mapsto \int k(x-y)(\phi(y)-\phi(x)) u(y)dy$ and by using Lemma \ref{lemma:poid} and \ref{lemm:noyaubassefreq}. Notice that this process gives us an explicit constant in term of $\mu_1$ and $\mu_2$.
By changing the variable $\displaystyle x=\frac{x'}{A}$, we deduce that:
\begin{align*}
    \left\|[|D|^{\alpha},\phi_{1,A}]u \right\|_{L^2}\leq \frac{C}{A^{\frac{\alpha-1}{2}}} \bigg|\frac{1}{(1+\mu_1)^2} - \frac{1}{(1+\mu_2)^2}\bigg|^{\frac{1}{2}}  \|u\Phi_{1,A} \|_{H^{\frac{\alpha}{2}}}.
\end{align*}
We obtain by definition of the Sobolev space:
\begin{align*}
    \left\|[|D|^{\alpha},\phi_{1,A}]u \right\|_{L^2} \leq C \bigg|\frac{1}{(1+\mu_1)^2} - \frac{1}{(1+\mu_2)^2}\bigg|^{\frac{1}{2}}  \|u\Phi_{1,A} \|_{H^1}.
\end{align*}
This concludes the proof of Lemma \ref{lemm:commH1}.
\end{proof}

\begin{proof}[Proof of Lemma \ref{lemm:commL2}]
The proof is based on the same arguments as the former lemmas. For the high frequency we use the pseudo-differential calculus, except that we use the function $\sqrt{\phi}$ instead of $\phi$. Using the Poisson bracket in \eqref{eq:Poisson_bracket_symbol}, we deduce that the commutator satisfies $[(1-\chi(D))\vert D \vert^\alpha,\sqrt{\phi}] \in \mathcal{S}^{\alpha-1, -1-\frac{\alpha}{2}}\subset \mathcal{S}^{\frac{\alpha}{2},0}$, and we can use the same arguments as above. For the low frequency we use the Schur's lemma (Lemma \ref{theo:schur}).
\end{proof}

\section{Proof of the coercivity property}\label{proof:claim_coer}

We prove the following result of coercivity which is time-independent, with $R_1$, $R_2$, $\tilde{R}_1$ and $\tilde{R}_2$ defined in \eqref{defi:R_i_modulation} and dependent on $\Gamma$ satisfying the condition \eqref{condition:Cond_Z}:
\begin{claim}\label{claim:coercivite}
Let $\eps\in H^{\frac{\alpha}{2}}(\R)$  satisfying the four orthogonality conditions:
\begin{align*}
    0 = \int \eps {R}_1 = \int \eps \partial_y {R}_1 =\int \eps {R}_2 = \int \eps \partial_y {R}_2, 
\end{align*}
and $\Gamma=(z_1,z_2,\mu_1,\mu_2)$ satisfying \eqref{condition:Cond_Z}. Then, there exists $A_2$, $Z^{*}_1$, $\kappa>0$ such that for all $A>A_2$ and $\Gamma$ satisfying (Cond$_{Z^*_1}$):
\begin{align*}
   \sum_{i=1}^2 \int\left( \eps|D|^{\alpha}\eps + \eps^2 - 3\widetilde{R}^2_i\eps^2\right)\psi_{i,A} \geq \kappa \|\eps\|^2_{H^{\frac{\alpha}{2}}}, \quad i=1,2,
\end{align*}
with $\displaystyle\psi_{1,A}:=\frac{1-\phi_A}{(1+\mu_1)^2}$ or $\displaystyle\psi_{2,A}:=\frac{\phi_A}{(1+\mu_2)^2}$.
\end{claim}

\begin{proof}
Since $\psi_{i,A}>0$, and $L$ is coercive, see \eqref{eq:coercivite}, we deduce that: 
\begin{align*}
\MoveEqLeft
    \int\left( \eps|D|^{\alpha}\eps + \eps^2 - 3\widetilde{R}^2_i\eps^2\right)\psi_{i,A} \\
    &=\int \left(|D|^{\frac{\alpha}{2}}\left(\eps\sqrt{\psi_{i,A}}\right)\right)^2 +\left(\eps\sqrt{\psi_{i,A}}\right)^2 -3 \widetilde{R}_i^2\left(\eps\sqrt{\psi_{i,A}}\right)^2 + \int \eps\sqrt{\psi_{i,A}}\left[|D|^{\alpha},\sqrt{\psi_{i,A}}\right]\eps\\
    &\geq \kappa_1\left\|\eps\sqrt{\psi_{i,A}}\right\|_{H^{\frac{\alpha}{2}}}^2 + \int \eps\sqrt{\psi_{i,A}}\left[|D|^{\alpha},\sqrt{\psi_{i,A}}\right]\eps - \frac{1}{\kappa_1} \left(\int\eps\sqrt{\psi_{i,A}}\widetilde{R}_i \right)^2 - \frac{1}{\kappa_1} \left(\int \eps\sqrt{\psi_{i,A}}\partial_y\widetilde{R}_i \right)^2.
\end{align*}
Since $\langle \xi\rangle^{\frac{\alpha}{2}}\geq \kappa_2( 1+|\xi|^{\frac{\alpha}{2}})$, we obtain that: 
\begin{align*}
    \left\|\eps\sqrt{\psi_{i,A}}\right\|_{H^{\frac{\alpha}{2}}}^2\geq \kappa_2\int \left(\eps^2 + (|D|^{\frac{\alpha}{2}}\eps)^2\right)\psi_{i,A} +\kappa_2\int (|D|^{\frac{\alpha}{2}}(\eps\sqrt{\psi_{i,A}}))^2 - (|D|^{\frac{\alpha}{2}}\eps)^2\psi_{i,A}.
\end{align*}
Notice that: 
\begin{align*}
    \int (|D|^{\frac{\alpha}{2}}(\eps\sqrt{\psi_{i,A}}))^2 - (|D|^{\frac{\alpha}{2}}\eps)^2\psi_{i,A}= 2 \int\left(|D|^{\frac{\alpha}{2}}(\eps\sqrt{\psi_{i,A}}) \right)[|D|^{\frac{\alpha}{2}},\sqrt{\psi_{i,A}}]\eps - \int \left([|D|^{\frac{\alpha}{2}},\sqrt{\psi_{i,A}}]\eps \right)^2 .
\end{align*}
 Using Lemma \ref{lemm:commL2} and Young's inequality, we obtain that:
\begin{align*}
     \kappa_1\left\|\eps\sqrt{\psi_{i,A}}\right\|_{H^{\frac{\alpha}{2}}}^2 +\int\eps\sqrt{\psi_{i,A}} \left[|D|^{\alpha},\sqrt{\psi_{i,A}}\right]\eps \geq& \kappa_1\kappa_2\int \left(\eps^2 + (|D|^{\frac{\alpha}{2}}\eps)^2\right)\psi_{i,A} 
     -\frac{C}{A^{\frac{\alpha}{2}}}\int \eps^2 + (|D|^{\frac{\alpha}{2}}\eps)^2. 
\end{align*}
Note that since $\eps\perp R_i$, we have that: 
\begin{align*}
    \int \eps\sqrt{\psi_{i,A}}\widetilde{R}_i = \int \eps\left(\sqrt{\psi_{i,A}}-1\right) R_i + \int \eps \sqrt{\psi_{i,A}} \left(\widetilde{R}_i-R_i\right)   
\end{align*}
Then, by using the Cauchy-Schwarz' inequality, \eqref{eq:partial_y_Ri}, we get that:
\begin{align*}
    \left(\int \eps\sqrt{\psi_{i,A}}\widetilde{R}_i \right)^2 + \left(\int \eps\sqrt{\psi_{i,A}}\partial_y\widetilde{R}_i \right)^2 \leq C \|\eps\|_{L^2}^2 \left(\frac{1}{z^{\alpha}}+ \|R_i- \tilde{R}_i \|_{H^1}^2 \right).
\end{align*}
Moreover, we have that $\psi_{1,A}+\psi_{2,A}\geq \kappa_3>0$. Therefore, we can conclude, with \eqref{eq:lambda_Q_DL}:
\begin{align*}
    \| R_i - \tilde{R}_i \|_{H^1} \leq C \mu_i^2,
\end{align*}
by taking $Z$ and $A>A_2$ large enough, that there exists $\kappa>0$ such that:
\begin{align*}
   \sum_{i=1}^{2} \int\left( \eps|D|^{\alpha}\eps + \eps^2 - 3\widetilde{R}^2_i\eps^2\right)\psi_{i,A}\geq \kappa \|\eps\|_{H^{\frac{\alpha}{2}}}^2.
\end{align*}
\end{proof}


\bibliographystyle{plain}
\bibliography{biblio}

\begin{thebibliography}{10}

\bibitem{albert1997model}
John~P. Albert, Jerry~L. Bona, and Jean-Claude Saut.
\newblock Model equations for waves in stratified fluids.
\newblock {\em Proceedings of the Royal Society of London. Series A:
  Mathematical, Physical and Engineering Sciences}, 453(1961):1233--1260, 1997.

\bibitem{alinhac2012operateurs}
Serge Alinhac and Patrick G{\'e}rard.
\newblock {\em Op{\'e}rateurs pseudo-diff{\'e}rentiels et th{\'e}oreme de
  {Nash-Moser}}.
\newblock EDP Sciences, 2012.

\bibitem{Ang18}
Jaime Angulo~Pava.
\newblock Stability properties of solitary waves for fractional {K}d{V} and
  {BBM} equations.
\newblock {\em Nonlinearity}, 31(3):920--956, 2018.

\bibitem{Ary22}
Shrey Aryan.
\newblock Existence of two-solitary waves with logarithmic distance for the
  nonlinear {K}lein-{G}ordon equation.
\newblock {\em Commun. Contemp. Math.}, 24(1):Paper No. 2050091, 25, 2022.

\bibitem{BK04}
J.~L. Bona and H.~Kalisch.
\newblock Singularity formation in the generalized {B}enjamin-{O}no equation.
\newblock {\em Discrete Contin. Dyn. Syst.}, 11(1):27--45, 2004.

\bibitem{CS07}
Luis Caffarelli and Luis Silvestre.
\newblock An extension problem related to the fractional {L}aplacian.
\newblock {\em Comm. Partial Differential Equations}, 32(7-9):1245--1260, 2007.

\bibitem{chow1982bifurcation}
Shui~Nee Chow and Jack~K. Hale.
\newblock {\em Methods of bifurcation theory}, volume 251 of {\em Grundlehren
  der Mathematischen Wissenschaften [Fundamental Principles of Mathematical
  Sciences]}.
\newblock Springer-Verlag, New York-Berlin, 1982.

\bibitem{CM18}
Vianney Combet and Yvan Martel.
\newblock Construction of multibubble solutions for the critical {GKDV}
  equation.
\newblock {\em SIAM J. Math. Anal.}, 50(4):3715--3790, 2018.

\bibitem{CMM11}
Rapha\"{e}l C\^{o}te, Yvan Martel, and Frank Merle.
\newblock Construction of multi-soliton solutions for the {$L^2$}-supercritical
  g{K}d{V} and {NLS} equations.
\newblock {\em Rev. Mat. Iberoam.}, 27(1):273--302, 2011.

\bibitem{CMYZ21}
Rapha\"{e}l C\^{o}te, Yvan Martel, Xu~Yuan, and Lifeng Zhao.
\newblock Description and classification of 2-solitary waves for nonlinear
  damped {K}lein-{G}ordon equations.
\newblock {\em Comm. Math. Phys.}, 388(3):1557--1601, 2021.

\bibitem{DD07}
Fran\c{c}oise Demengel and Gilbert Demengel.
\newblock {\em Functional spaces for the theory of elliptic partial
  differential equations}.
\newblock Universitext. Springer, London; EDP Sciences, Les Ulis, 2012.
\newblock Translated from the 2007 French original by Reinie Ern\'{e}.

\bibitem{eychenne2021asymptotic}
Arnaud Eychenne.
\newblock Asymptotic ${N}$-soliton-like solutions of the fractional
  {K}orteweg-de {V}ries equation.
\newblock {\em arXiv preprint arXiv:2112.11278}, 2021.

\bibitem{EV22}
Arnaud Eychenne and Frederic Valet.
\newblock Asymptotic of non-linear ground states for fractional {L}aplacian.
\newblock {\em to appear}, 2022.

\bibitem{FL13}
Rupert~L. Frank and Enno Lenzmann.
\newblock Uniqueness of non-linear ground states for fractional {L}aplacians in
  {$\mathbb{R}$}.
\newblock {\em Acta Math.}, 210(2):261--318, 2013.

\bibitem{FLS16}
Rupert~L. Frank, Enno Lenzmann, and Luis Silvestre.
\newblock Uniqueness of radial solutions for the fractional {L}aplacian.
\newblock {\em Comm. Pure Appl. Math.}, 69(9):1671--1726, 2016.

\bibitem{GLPR18}
Patrick G\'{e}rard, Enno Lenzmann, Oana Pocovnicu, and Pierre Rapha\"{e}l.
\newblock A two-soliton with transient turbulent regime for the cubic half-wave
  equation on the real line.
\newblock {\em Ann. PDE}, 4(1):Paper No. 7, 166, 2018.

\bibitem{goubet2009weak}
O.~Goubet and L.~Molinet.
\newblock Global attractor for weakly damped nonlinear {S}chr\"{o}dinger
  equations in {$L^2(\mathbb{ R})$}.
\newblock {\em Nonlinear Anal.}, 71(1-2):317--320, 2009.

\bibitem{Guo12}
Zihua Guo.
\newblock Local well-posedness for dispersion generalized {B}enjamin-{O}no
  equations in {S}obolev spaces.
\newblock {\em J. Differential Equations}, 252(3):2053--2084, 2012.

\bibitem{GH22}
Zihua Guo and Chunyan Huang.
\newblock Well-posedness of the modified dispersion-generalized
  {B}enjamin-{O}no equations in modulation spaces.
\newblock {\em J. Math. Anal. Appl.}, 509(1):Paper No. 125933, 20, 2022.

\bibitem{halmos2012bounded}
Paul~R. Halmos and Viakalathur~S. Sunder.
\newblock {\em Bounded integral operators on {$L^2$} spaces}, volume~96.
\newblock Springer Science \& Business Media, 2012.

\bibitem{hormander2007analysis}
Lars H{\"o}rmander.
\newblock {\em The analysis of linear partial differential operators III:
  Pseudo-differential operators}.
\newblock Springer Science \& Business Media, 2007.

\bibitem{JL22}
Jacek Jendrej and Andrew Lawrie.
\newblock An asymptotic expansion of two-bubble wave maps in high equivariance
  classes.
\newblock {\em Anal. PDE}, 15(2):327--403, 2022.

\bibitem{KMV17}
Henrik Kalisch, Daulet Moldabayev, and Olivier Verdier.
\newblock A numerical study of nonlinear dispersive wave models with
  {S}pec{T}ra{VV}ave.
\newblock {\em Electron. J. Differential Equations}, pages Paper No. 62, 23,
  2017.

\bibitem{KM09}
Carlos~E. Kenig and Yvan Martel.
\newblock Asymptotic stability of solitons for the {B}enjamin-{O}no equation.
\newblock {\em Rev. Mat. Iberoam.}, 25(3):909--970, 2009.

\bibitem{kenig2011local}
Carlos~E. Kenig, Yvan Martel, and Luc Robbiano.
\newblock Local well-posedness and blow-up in the energy space for a class of
  ${L}^2$ critical dispersion generalized {B}enjamin--{O}no equations.
\newblock {\em Annales de l'Institut Henri Poincar\'{e}. Analyse Non
  Lin\'{e}aire}, 28(6):853--887, 2011.

\bibitem{KT06}
Carlos~E. Kenig and Hideo Takaoka.
\newblock Global wellposedness of the modified {B}enjamin-{O}no equation with
  initial data in {$H^{1/2}$}.
\newblock {\em Int. Math. Res. Not.}, pages Art. ID 95702, 44, 2006.

\bibitem{KS21}
Kihyun Kim and Robert Schippa.
\newblock Low regularity well-posedness for generalized {B}enjamin-{O}no
  equations on the circle.
\newblock {\em J. Hyperbolic Differ. Equ.}, 18(4):931--984, 2021.

\bibitem{KLPS18}
Christian Klein, Felipe Linares, Didier Pilod, and Jean-Claude Saut.
\newblock On {W}hitham and related equations.
\newblock {\em Stud. Appl. Math.}, 140(2):133--177, 2018.

\bibitem{KSW22}
Christian Klein, Jean-Claude Saut, and Yuexun Wang.
\newblock On the modified fractional {K}orteweg--de {V}ries and related
  equations.
\newblock {\em Nonlinearity}, 35(3):1170--1212, 2022.

\bibitem{KT20}
Herbert Koch and Daniel Tataru.
\newblock Multisolitons for the cubic nls in 1-d and their stability, 2020.

\bibitem{KMR09}
Joachim Krieger, Yvan Martel, and Pierre Rapha\"{e}l.
\newblock Two-soliton solutions to the three-dimensional gravitational
  {H}artree equation.
\newblock {\em Comm. Pure Appl. Math.}, 62(11):1501--1550, 2009.

\bibitem{LW22}
Yang Lan and Zhong Wang.
\newblock Strongly interacting multi-solitons for generalized benjamin-ono
  equations, 2022.

\bibitem{MM11}
Yvan Martel and Frank Merle.
\newblock Description of two soliton collision for the quartic g{K}d{V}
  equation.
\newblock {\em Ann. of Math. (2)}, 174(2):757--857, 2011.

\bibitem{MM11inelastic}
Yvan Martel and Frank Merle.
\newblock Inelastic interaction of nearly equal solitons for the quartic
  g{K}d{V} equation.
\newblock {\em Invent. Math.}, 183(3):563--648, 2011.

\bibitem{MMT02}
Yvan Martel, Frank Merle, and Tai-Peng Tsai.
\newblock Stability and asymptotic stability in the energy space of the sum of
  {$N$} solitons for subcritical g{K}d{V} equations.
\newblock {\em Comm. Math. Phys.}, 231(2):347--373, 2002.

\bibitem{MN20}
Yvan Martel and Ti\'{\^e}n-Vinh Nguy\~{\^e}n.
\newblock Construction of 2-solitons with logarithmic distance for the
  one-dimensional cubic {S}chr\"{o}dinger system.
\newblock {\em Discrete Contin. Dyn. Syst.}, 40(3):1595--1620, 2020.

\bibitem{MP17}
Yvan Martel and Didier Pilod.
\newblock Construction of a minimal mass blow up solution of the modified
  {B}enjamin-{O}no equation.
\newblock {\em Math. Ann.}, 369(1-2):153--245, 2017.

\bibitem{MR18}
Yvan Martel and Pierre Rapha\"{e}l.
\newblock Strongly interacting blow up bubbles for the mass critical nonlinear
  {S}chr\"{o}dinger equation.
\newblock {\em Ann. Sci. \'{E}c. Norm. Sup\'{e}r. (4)}, 51(3):701--737, 2018.

\bibitem{Mer01}
Frank Merle.
\newblock Existence of blow-up solutions in the energy space for the critical
  generalized {K}d{V} equation.
\newblock {\em J. Amer. Math. Soc.}, 14(3):555--578, 2001.

\bibitem{Miz03}
Tetsu Mizumachi.
\newblock Weak interaction between solitary waves of the generalized {K}d{V}
  equations.
\newblock {\em SIAM J. Math. Anal.}, 35(4):1042--1080, 2003.

\bibitem{MT22}
Luc Molinet and Tomoyuki Tanaka.
\newblock Unconditional well-posedness for some nonlinear periodic
  one-dimensional dispersive equations.
\newblock {\em J. Funct. Anal.}, 283(1):Paper No. 109490, 45, 2022.

\bibitem{Mun10}
Claudio Mu\~{n}oz.
\newblock On the inelastic two-soliton collision for g{K}d{V} equations with
  general nonlinearity.
\newblock {\em Int. Math. Res. Not. IMRN}, 9:1624--1719, 2010.

\bibitem{LNP22}
F\'{a}bio Natali, Uyen Le, and Dmitry~E. Pelinovsky.
\newblock Periodic waves in the fractional modified {K}orteweg--de {V}ries
  equation.
\newblock {\em J. Dynam. Differential Equations}, 34(2):1601--1640, 2022.

\bibitem{Ngu17}
Ti\'{\^e}n-Vinh Nguy\~{\^e}n.
\newblock Strongly interacting multi-solitons with logarithmic relative
  distance for the g{K}d{V} equation.
\newblock {\em Nonlinearity}, 30(12):4614--4648, 2017.

\bibitem{Ngu19}
Ti\'{\^e}n-Vinh Nguy\~{\^e}n.
\newblock Existence of multi-solitary waves with logarithmic relative distances
  for the {NLS} equation.
\newblock {\em C. R. Math. Acad. Sci. Paris}, 357(1):13--58, 2019.

\bibitem{OS98}
Yu.~N. Ovchinnikov and I.~M. Sigal.
\newblock Long-time behaviour of {G}inzburg-{L}andau vortices.
\newblock {\em Nonlinearity}, 11(5):1295--1309, 1998.

\bibitem{SW21}
Jean-Claude Saut and Yuexun Wang.
\newblock Global dynamics of small solutions to the modified fractional
  {K}orteweg--de {V}ries and fractional cubic nonlinear {S}chr\"{o}dinger
  equations.
\newblock {\em Comm. Partial Differential Equations}, 46(10):1851--1891, 2021.

\bibitem{SV96}
Victor~I. Shrira and Vyacheslav~V. Voronovich.
\newblock Nonlinear dynamics of vorticity waves in the coastal zone.
\newblock {\em J. Fluid Mech.}, 326:181--203, 1996.

\bibitem{WO82}
Miki Wadati and Kenji Ohkuma.
\newblock Multiple-pole solutions of the modified {K}orteweg-de {V}ries
  equation.
\newblock {\em J. Phys. Soc. Japan}, 51(6):2029--2035, 1982.

\bibitem{weinstein1985modulation}
Michael~I. Weinstein.
\newblock Modulational stability of ground states of nonlinear
  {S}chr\"{o}dinger equations.
\newblock {\em SIAM J. Math. Anal.}, 16(3):472--491, 1985.

\bibitem{weinstein1987existence}
Michael~I. Weinstein.
\newblock Existence and dynamic stability of solitary wave solutions of
  equations arising in long wave propagation.
\newblock {\em Communication in Partial Differential Equation},
  12(10):1133--1173, 1987.

\end{thebibliography}

\end{document}